\documentclass[11pt]{amsart}
\usepackage{geometry}
\newcommand{\A}{\mathbb{A}}

\newcommand{\C}{\mathbb{C}}

\newcommand{\Hbb}{\mathbb{H}}

\newcommand{\M}{\mathbb{M}}

\newcommand{\Q}{\mathbb{Q}}
\newcommand{\R}{\mathbb{R}}

\newcommand{\Z}{\mathbb{Z}}

\newcommand{\Ocal}{\mathcal{O}}

\newcommand{\Tcal}{\mathcal{T}}

\newcommand{\fbf}{\mathbf{f}}

\newcommand{\jbf}{\mathbf{j}}

\newcommand{\zbf}{\mathbf{z}}

\newcommand{\Hbf}{\mathbf{H}}

\newcommand{\Tbf}{\mathbf{T}}

\usepackage{imakeidx}
\makeindex[columns=3]
\makeindex[columns=2,options= -s index_style_xyz.ist]
\usepackage[utf8]{inputenc}
\usepackage{tgtermes}
\usepackage[english]{babel}
\usepackage{amsmath,
	mathtools,
	amssymb,
	graphicx,
	bbm,
	xcolor,
	amsthm,
	mathrsfs,
	enumitem,
	mathdots}
\usepackage[all]{xy}
\usepackage[backref=page]{hyperref}
\hypersetup{}
\usepackage{tikz-cd}
\usetikzlibrary{decorations.pathmorphing}
\usepackage{listings}
\usepackage[normalem]{ulem}

\theoremstyle{plain}

\newtheorem{theorem}{Theorem}[section]
\newtheorem{corollary}[theorem]{Corollary}
\newtheorem{definition}[theorem]{Definition}
\newtheorem{lemma}[theorem]{Lemma}
\newtheorem{proposition}[theorem]{Proposition}

\newtheorem{remark}[theorem]{Remark}

\newtheorem{examples}[theorem]{Examples}

\newtheorem{questions}[theorem]{Questions}

\newtheorem*{theorem-no-num}{Theorem}
\newtheorem*{proposition-no-num}{Proposition}
\newtheorem*{corollary-no-num}{Corollary}

\makeatletter
\@namedef{subjclassname@2020}{%
	\textup{2020} Mathematics Subject Classification}
\makeatother

\newtheorem{taggedtheoremx}{}
\newenvironment{taggedtheorem}[1]
{\renewcommand\thetaggedtheoremx{#1}\taggedtheoremx}
{\endtaggedtheoremx}

\newcommand{\GLmr}{\mathrm{GL}}
\newcommand{\Symmr}{\mathrm{Sym}}
\newcommand{\Spmr}{\mathrm{Sp}}

\DeclareFontFamily{U}{wncy}{}
\DeclareFontShape{U}{wncy}{m}{n}{<->wncyr10}{}
\DeclareSymbolFont{mcy}{U}{wncy}{m}{n}
\DeclareMathSymbol{\Sha}{\mathord}{mcy}{"58}

\usepackage{hyperref}

\begin{document}
	\title{Non-vanishing mod $p$ of theta lifts for $(\mathrm{O}_{2n+1},\mathrm{Mp}_{4n})$}
	\author{Xiaoyu Zhang}
	\address{Universit\"{a}t Duisburg-Essen,
		Fakult\"{a}t f\"{u}r Mathematik,
		Mathematikcarr\'{e}e
		Thea-Leymann-Straße 9,
		45127 Essen,
		Germany}
	\date{\today}
	\email{xiaoyu.zhang@uni-due.de}
	\maketitle

	\begin{abstract}
		We establish the non-vanishing mod $p$ of global theta lifts from
		an odd definite orthogonal group
		$\mathrm{O}_{2n+1}$
		over $\Q$
		to a metaplectic group
		$\mathrm{Mp}_{4n}$ over $\Q$
		under mild conditions. The problem is closely related to non-vanishing modulo $p$ of toric integrals on $\mathrm{O}_{2n+1}$. For this, we exploit the distribution properties of toric orbits of unipotent elements on $\mathrm{O}_{2n+1}(\Q_\ell)$ using Ratner's theorems on unipotent flows and we deduce that the toric integral of a $p$-primitive automorphic form on $\mathrm{O}_{2n+1}$ is non-zero modulo $p$ for infinitely many characters.
	\end{abstract}
	\tableofcontents

	\section{Introduction}
	It is a fundamental problem in representation theory and number theory to relate automorphic forms/representations on a reductive linear algebraic group $H_{/\Q}$ to automorphic forms/representations on another reductive linear algebraic group $G_{/\Q}$. One such instance is the Langlands functoriality. There are several tools of studying such problems, and in this article we shall work with one such tool, the theta correspondence. In the theory of theta correspondence, the non-vanishing of theta lifts of automorphic forms is a central question and has attracted much attention. Most of the work in the literature relies on the local-global criterion initiated by Rallis (\cite{RallisSchiffman1978,RallisSchiffman1980,RallisSchiffman1981,GanTakeda2011,GanQiuTakeda2012}). In this article, we will use a quite different method and relate the non-vanishing problem to toric integrals. More precisely, in this article we study modulo $p$ properties of 
	global theta lifts of automorphic forms
	between certain reductive dual pairs $(H,G)$ over $\Q$ using distribution results of toric orbits of unipotent elements on $H(\Q_\ell)$.
	In particular,
	we show that under mild conditions,
	a $p$-primitive automorphic form has $p$-primitive theta lift
	by a carefully chosen theta series.\footnote{This is a substantial revision of previous versions, including one titled `Non-vanishing mod $p$ of theta lifts' (\url{https://arxiv.org/abs/2203.05359}). In the previous versions, we worked over a totally real number field instead of the rational numbers $\Q$. This generality seems to have obscured the main idea of the work and to have burdened the notations. So to make our work more transparent, we work over $\Q$ and clarifies many places compared to the previous versions. Moreover, in previous versions, we also had applications to Bloch-Kato conjectures. However, we feel that the present article is already very long and is of independent interest, so it is better to put this application in a separate article.}

	\subsection{Non-vanishing modulo $p$ of theta lifts}
	Suppose $H$ and $G$, as above, are reductive algebraic groups over $\mathbb{Q}$, $\mathbf{f}$ an automorphic form on $H(\mathbb{Q})\backslash H(\mathbb{A})$ and $\Theta_{\phi}(\pi)$ its theta lift to $G(\mathbb{Q})\backslash G(\mathbb{A})$ by a Bruhat-Schwartz function $\phi$. Then we are interested in the following

	\begin{questions}
		Suppose
		$\mathbf{f}$ satisfies certain arithmetic properties
		(for example, algebraic, rational, integral,
		$p$-integral, $p$-primitive, etc),
		is it also true for $\Theta_{\phi,\mathbf{f}}$
		(for certain $\phi$)?
	\end{questions}

	These questions have important applications on non-vanishing (modulo $p$) of $L$-values and other important arithmetic invariants.
	For example
	\begin{enumerate}
		\item
		$(H,G)=(\mathrm{GU}_2,\mathrm{GU}_3)$
		in \cite{Finis2000} and \cite{Finis2006}
		where the $p$-adic valuation of certain Hecke $L$-value is deduced;

		\item
		$(H,G)=(\mathrm{GL}_2,\mathrm{GO}_B)$
		for indefinite quaternion algebra $B$ in
		\cite{Prasanna2006}
		where the integrality of quotients of Petersson inner products of
		modular forms is deduced
		and in
		\cite{Hida2010}
		where
		the $\mu$-invariant of certain Hecke $L$-functions is computed;

		\item
		$(H,G)=(\mathrm{GSO}_B,\mathrm{GSp}_4)$
		for definite quaternion algebra $B$ in
		\cite{HsiehNamikawa2017}
		and \cite{ChidaHsieh2016}
		where
		the $p$-primitivity of certain normalized $L$-values
		and the vanishing of certain $\mu$-invariants are deduced.
		This is also used
		(in combination with \cite{Prasanna2006})
		in a crucial way in
		\cite{TilouineUrban2018}
		for deducing integral relations of certain periods of modular forms.		
	\end{enumerate}

	In this article we want to study the $p$-primitivity problem of
	$\Theta_{\phi,\mathbf{f}}$ for
	the reductive dual pair $(H,G)$ over $\Z$ where $H$ is the orthogonal group of a non-degenerate quadratic module $(U,\langle-,-\rangle_{U})$ over $\Z$ of odd rank $2n+1\ge3$ and $G$ is the symplectic group of a non-degenerate symplectic module $(V,\langle-,-\rangle_V)$ over $\Z$ of rank $4n$. We will assume that under the standard basis of $U$, $\langle-,-\rangle_U$ is represented by a diagonal matrix $Q_U=\mathrm{diag}(\delta_1,\cdots,\delta_{2n+1})$ with $\delta_1,\cdots,\delta_{2n+1}\in\Z_{>0}$. In this case, the theta lift $\Theta_{\phi,\mathbf{f}}$ of an automorphic form $\mathbf{f}$ on $H(\mathbb{A})$ is an automorphic form on $\widetilde{G}(\mathbb{A})$ where $\widetilde{G}$ is the metaplectic group of $(V,\langle-,-\rangle)$.

	To state our main result, we need some more notations: we fix a rational prime $p$ and a field isomorphism $\C\simeq\C_p$, let $\mathcal{O}$ be the ring of integers of $\mathbb{C}_p$
	and $\mathfrak{P}$
	its maximal ideal, $\kappa$ the residue field. Fix then an irreducible algebraic representation $(\rho_{\lambda},\M_{\lambda})$
	of
	$H_{/\mathcal{O}}$
	of highest weight
	$\lambda\in\Z^n$ and an irreducible algebraic representation $(\rho_{\tau^\circ},\M_{\tau^\circ})$ of $\mathrm{GL}_{2n/\mathcal{O}}$ of highest weight $\tau^\circ\in\mathbb{Z}^{2n}$.\footnote{Under the theta correspondence between $H$ and $G$, these two highest weights $\rho,\tau^\circ$ are related to each other. See (\ref{tau^circ}).}

	For an automorphic form	$\mathbf{f}$
	on
	$H(\mathbb{Q})\backslash H(\mathbb{A})$
	of weight $\lambda$,
	we say it is \textit{$p$-integral}
	if its $p$-adic avatar
	\[
	\mathbf{f}_p\colon H(\mathbb{A}_f)\rightarrow\M_{\lambda}(\C_p)
	\]
	takes values in the lattice $\M_{\lambda}(\mathcal{O})$ of $\M_\lambda(\C_p)$.
	We say $\mathbf{f}$
	is \textit{$p$-primitive}
	if furthermore $\overline{\mathbf{f}_p}:=\mathbf{f}_p(\mathrm{mod}\,\mathfrak{P})$ is not identically zero.
	For a holomorphic automorphic form $\mathbf{F}$ on
	$G(\mathbb{Q})\backslash \widetilde{G}(\mathbb{A})$
	of weight $\tau^\circ$,
	we say it is \textit{$p$-integral} if
	the classical Siegel modular form associated to $\mathbf{F}$
	has Fourier coefficients in $\M_{\tau^\circ}(\mathcal{O})$ (a lattice in $\M_{\tau^\circ}(\C_p)\simeq\M_{\tau^\circ}(\C)$). We say $\mathbf{F}$ is \textit{$p$-primitive} if
	furthermore one of the Fourier coefficients is non-zero modulo
	$\mathfrak{P}$
	(see \S \ref{automorphic forms on H and G} for the precise definitions).

	We fix a finite set $\mathbb{S}$ of finite places of $\Q$ including those dividing $2\delta_1\cdots\delta_{2n+1}$ and a positive integer $r_\ell$ for each $\ell\in\mathbb{S}$ (satisfying some technical conditions, see the beginning of §\ref{Bruhat-Schwartz function}). Let 
	$K^{(\mathbb{S})}\subset H(\widehat{\mathbb{Z}})$ be a compact open subgroup (see the beginning of §\ref{Bruhat-Schwartz function}). We fix an isotropic decomposition $V=V^+\oplus V^-$ and write $W=U\otimes V$, $W^\pm=U\otimes V^\pm$. After fixing a basis for $U,V^+$, we identify $W^+$ with $\mathrm{M}_{2n+1,2n}(\Z)$.

	The main result of the article is as follows:
	\begin{theorem}\label{main result for special case}
		We construct a Bruhat-Schwartz function
		\[
		\phi_{\lambda}\in
		\mathcal{S}(W^+(\mathbb{A}))
		\otimes_\C
		(\M_\lambda(\C)\otimes\M_{\tau^\circ}(\C))
		\]
		such that the following holds:
		\begin{enumerate}
			\item
			Assume conditions \ref{condition on p} and \ref{condition p not dividing the level}.	Let	$\mathbf{f}\in\mathcal{A}_{\rho_{\lambda}}( H,K^{(\mathbb{S})})$ be a $p$-integral automorphic form on $H(\mathbb{A})$ of weight $\lambda$ and of level $K^{(\mathbb{S})}$, then the theta lift $\Theta_{\phi_{\lambda},\mathbf{f}}\in
			\mathcal{A}_{\rho_{\tau^\circ}}			(\widetilde{G},\Gamma_0(2,N_{\mathbb{S}}),\chi_{U}^\circ)$
			is a holomorphic cuspidal genuine automorphic form on $\widetilde{G}(\mathbb{A})$ and it is $p$-integral.

			\item
			Assume conditions \ref{condition on lambda}, \ref{condition on p}, \ref{condition p not dividing the level} and \ref{condition-delta are squares in Q}.
			If furthermore $\mathbf{f}$ is
			$p$-primitive, then so is $\Theta_{\phi_{\lambda},\mathbf{f}}$.
		\end{enumerate}		
	\end{theorem}
	See \S \ref{automorphic forms on H and G}
	for more details on the definitions of the notations.
	The first part is
	Theorem \ref{p-integrality of lift of f}
	and the second part is
	Theorem \ref{main result of the article}. In fact, we can show that there are infinitely many Fourier coefficients of $\Theta_{\phi_{\lambda},\mathbf{f}}$ which are non-zero modulo $\mathfrak{P}$.
	This result generalizes
	\cite[Theorem 5.3]{HsiehNamikawa2017}\footnote{Our theorem also recovers the main result (Theorem A) of \cite{HsiehNamikawa2017}: indeed, the orthogonal group $H'$ in \emph{loc.cit} is just $H\times H$ for $n=1$ and the automorphic form $\mathbf{f}'$ on $H'(\mathbb{A})$ is thus a pair $(\mathbf{f}_1,\mathbf{f}_2)$ where $\mathbf{f}_1,\mathbf{f}_2$ are automorphic forms on $H(\mathbb{A})$. Thus a necessary condition for the non-vanishing of the theta lift $\Theta_{\phi,\mathbf{f}'}$ to $G(\mathbb{A})$ is the local root condition (the condition (LR) in \emph{loc.cit}). If we assume that (LR) is satisfied, then applying our theorem to both $\mathbf{f}_1$ and $\mathbf{f}_2$ gives the non-vanishing modulo $p$ of their theta lifts $\Theta_{ \phi_{\lambda_1},\mathbf{f}_1}$ and $\Theta_{ \phi_{\lambda_2},\mathbf{f}_2}$, and thus the non-vanishing modulo $p$ of $\Theta_{\phi,\mathbf{f}'}=\Theta_{\phi_{\lambda_1},\mathbf{f}_1}\times\Theta_{\phi_{\lambda_2},\mathbf{f}_2}$.} and \cite[\S 6]{Jia2010}.
	Our work can be seen as an intermediate step in transferring the well established characteristic $0$
	global theta correspondences to
	modulo $p$ global theta correspondences.

	Our work gives also a proof of the non-vanishing (\textit{without} modulo $p$)
	of the corresponding theta lifts by a completely different method as in previous works using local-global criterion. Besides we can use the above theorem 
	to prove new cases of $p$-part Bloch-Kato conjectures using the Rallis inner product formula, as done in \cite{AgarwalKlosin2013,BochererDummiganSchulzePillot2012,Jia2010}, etc.
	We will discuss this in a separate work.

	We next explain the strategy of the proof for the second part of Theorem \ref{main result for special case}. There are two main technical steps,
	which can be summarized as follows:
	(1) we need to show the non-vanishing modulo $p$ of certain toric periods,
	(2) we need to show certain toric orbits satisfy some equidistribution property.
	In the following we explain the first point and refer the reader to the second part of this introduction
	for the second point.

	We fix an element $\mathbf{z}\in W^+$ (see (\ref{z})) and write $S_{\mathbf{z}}=\mathbf{z}^\mathrm{t}Q_U\mathbf{z}$. We first express the $S_{\mathbf{z}}$-th Fourier coefficients $a_{\Theta_{ \phi_{\lambda},\mathbf{f}}}(S_{\mathbf{z}})$ of the classical Siegel modular form $\Theta_{ \phi_{\lambda},\mathbf{f}}^\ast$ associated to $\Theta_{\phi_{\lambda},\mathbf{f}}$ using a sum of $\mathbf{f}(h)$	for $h$
	running through a finite set $[\mathcal{E}_{\mathbf{z},K^{(\mathbb{S})}}]$	
	\[
	a_{\Theta_{ \phi_{\lambda},\mathbf{f}}}(S_{\mathbf{z}})
	=
	[H(\widehat{\Z}):K^{(\mathbb{S})}]^{-1}
	\sum_{[h_f]\in[\mathcal{E}_{\mathbf{z},K^{(\mathbb{S})}}]}
	w_{\mathbf{z},h_f}
	\langle
	\Delta_{\lambda}(\mathcal{T}\mathbf{z}),
	\mathbf{f}(h_f)
	\rangle_{W,U}.
	\]
	(we refer to Proposition \ref{Fourier coefficient of classicla Siegel modular form} for the details and explanation of the notations).

	Similarly, we fix a maximal torus $T$ of $H$ and for a character of finite order $\psi\colon\colon T(\mathbb{Q})\backslash T(\mathbb{A}_f)\to\C^\times$, we have the toric period (see Definition \ref{toric period} for more details)
	\[
	\mathbf{B}_{S_{\mathbf{z}},\psi}(g)
	=
	\mathbf{e}_\infty(-i\mathrm{Tr}(S))
	\int\limits_{T(\Q)\backslash T(\mathbb{A}_f)}
	\mathbf{a}_{\Theta_{\phi_{\lambda},\mathbf{f}},S_{\mathbf{z}}}
	(\mathbf{j}(t)g)
	\psi(t)
	dt.
	\]
	Here
	\begin{enumerate}
		\item 
		$\mathbf{a}_{\Theta_{ \phi_{\lambda},\mathbf{f}},S_\mathbf{z}}\colon\widetilde{G}(\mathbb{A})\to\mathbb{M}_\tau(\C)$ is the $S_\mathbf{z}$-th Fourier coefficient of $\Theta_{ \phi_{\lambda},\mathbf{f}}$ (see Definition \ref{definition of Siegel modular forms}), it is related to $a_{\Theta_{ \phi_{\lambda},\mathbf{f}},S_\mathbf{z}}$ as in (\ref{a_{f,S} and a_f(S)});

		\item 
		$\mathbf{j}\colon T\to G$ is a morphism of algebraic groups from the maximal torus $T$ of $H$ to $G$, sending $t$ to $\mathrm{diag}(\jbf'(t),\jbf'(t)^{-\mathrm{t}})$ where $\jbf'(t)=(\zbf')^{-1}t\zbf'$ with $\zbf'$ the first $2n$ rows of $\zbf$ (see also (\ref{j,j'})) .
	\end{enumerate}	
	A simple argument shows that if $\mathbf{f}$ is $p$-integral, then
	$\mathbf{B}_{S_{\mathbf{z}},\psi}$ also takes values in $\M_{\tau^\circ}(\mathcal{O})$. Moreover if $\mathbf{B}_{S_{\mathbf{z}},\psi}$
	is non-zero modulo $\mathfrak{P}$, then there exists some $a_{\Theta_{ \phi_{\lambda},\mathbf{f}}}(S)$ which is non-zero modulo $\mathfrak{P}$
	($S$ is related to $S_{\mathbf{z}}$, see Theorem \ref{p-primitivity of theta lift of f}). So it suffices to show that there exists some character $\psi$
	such that $\mathbf{B}_{S_{\mathbf{z}},\psi}$ is non-zero modulo $\mathfrak{P}$:
	\begin{theorem}
		Assume \ref{condition on lambda}, \ref{condition on p}, \ref{condition p not dividing the level} and \ref{condition-delta are squares in Q},
		then there are infinitely many characters of finite order
		\[
		\psi
		\colon
		T(\mathbb{Q})\backslash
		T(\mathbb{A}_f)
		\rightarrow
		\mathbb{C}^\times
		\]
		such that $\mathbf{B}_{S_{\mathbf{z}},\psi}$ is non-zero modulo $\mathfrak{P}$.
	\end{theorem}
	See Theorem \ref{theorem on the non-vanishing mod p of toric periods} for a more precise formulation of the above result. To prove this theorem, we first express the toric period as a finite sum over $\mathcal{R}$ of the product of
	$\psi$ and certain function $\overline{F}_{\mathbf{f},\underline{r}',\underline{\alpha}'}$ (valued in a $\kappa$-vector space) related to $\mathbf{f}$ (see (\ref{widetilde{F}})), where $\mathcal{R}$
	is a finite subset of $T(\mathbb{A})$. More precisely, we write $F_\fbf$ for the function on $H(\A_f)$ sending $g$ to $\langle\Delta_\lambda(\Tcal\zbf),\fbf_p(g)\rangle_{W,U}\in\M_\tau(\C)$. One can show that for $p$-integral $\fbf$, $F_\fbf$ takes values in $\M_\tau(\Ocal)$. Then $\overline{F}_{\fbf,\underline{r}',\underline{\alpha}'}$ is a certain $\psi$-weighted sum of $\overline{F}_\fbf$. Then we consider those characters $\psi$ of $\ell$-power conductors for some carefully chosen fixed prime $\ell\ne p$ and choose $\mathcal{R}$ satisfying some non-commensurability condition at $\ell$ (see (\ref{pairwise distinct})). The assumptions \ref{condition on lambda} and \ref{condition on p} imply that the image of $\overline{\mathbf{f}_p}$ generates $\M_{\lambda}(\kappa)$ and $\M_{\lambda}$ is of dimension $>1$ imply that $\overline{F}_{\mathbf{f},\underline{r}',\underline{\alpha}'}$
	is \textit{not} $\widetilde{H}^1(\mathbb{A}_f)$-invariant. Here $\widetilde{H}^1(\mathbb{A}_f)$ is the image of a certain spin group $\mathbf{H}^1(\mathbb{A}_f)$ in $H(\mathbb{A}_f)$. Now we can apply Theorem \ref{main theorem-2} to conclude that for any sufficiently large conductor $\underline{r}$, there exists a character $\psi$ of conductor $\underline{r}$ such that
	\[
	\sum_{t'\in\mathcal{R}}
	\overline{F}_{\mathbf{f},\underline{r}',\underline{\alpha}'}(t'\varsigma_{\xi(\underline{r})})
	\psi(t')
	\neq0.
	\]
	Theorem \ref{main theorem-2} is a result on ergodic theory and is a generalization of \cite{CornutVatsal2005}. Since it is of quite different flavour from the main topic of this article, we put it in the appendix. See also the next subsection for ideas involved in the proof of this theorem.

	\begin{remark}
		If the ranks of $U$ and $V$ satisfy $\mathrm{rk}(U)\leq\mathrm{rk}(V)/2$, then we are in the \emph{stable range} for the theta correspondence of the pair $(H,G)$. It is well-known that for an irreducible automorphic representation $\pi$ of $H(\mathbb{A})$, the theta lift $\Theta(\pi)$ of $\pi$ to $G$ is always non-zero (\cite{JianshuLi1992}). In fact, in this case and under certain conditions on the level of $\mathbf{f}$, one can relate the Fourier coefficient $a_{\Theta_{ \phi_{\lambda},\mathbf{f}}}(S_z)$ to a particular value of $\mathbf{f}$ (instead of a finite sum!) for $z$ of maximal rank (equal to $\mathrm{rk}(U)$). Thus it follows immediately that $\Theta_{\phi_{\lambda},\mathbf{f}}$ is non-zero. By suitably choosing the Bruhat-Schwartz function, one can even show that this theta lift is $p$-primitive if $\mathbf{f}$ is so. However, as long as $\mathrm{rk}(U)>\mathrm{rk}(V)/2$, this method fails and our work can be seen as the first step toward the non-vanishing modulo $p$ of theta lifts in this case.
	\end{remark}

	\subsection{Toric orbits of unipotent elements}

	We first recall the following important result of Cornut-Vatsal (\cite[Corollary 2.10]{CornutVatsal2005}): let $B$ be a definite quaternion algebra over $\Q$ and write $\mathbf{H}=B^\times$, viewed as an algebraic group over $\Q$. Write $\mathbf{H}^1$ for the derived subgroup of $\mathbf{H}$. Fix a maximal torus $\mathbf{T}\subset\mathbf{H}$, a prime number $\ell$, a positive integer $r$ and a compact open subgroup $K$ of $\mathbf{H}(\mathbb{A}_f)$. We then fix $r$ element $g_1,\cdots,g_r$ in $\mathbf{T}(\mathbb{A}_f)$. Then we have the following maps
	\[
	\mathbf{T}(\Q)
	\backslash
	\mathbf{H}(\mathbb{A}_f)/K
	\xrightarrow{\mathfrak{R}_K^r}
	(\mathbf{H}(\Q)\backslash \mathbf{H}(\mathbb{A}_f)/K)^r
	\xrightarrow{\mathfrak{A}_K^r}
	(\mathbf{H}(\Q)\backslash \mathbf{H}(\mathbb{A}_f)/K\mathbf{H}^1(\mathbb{A}_f))^r,
	\]
	where
	$\mathfrak{R}_K^r(x)
	=
	(\mathbf{H}(\Q)g_ix)_{i=1}^r
	$
	and
	$\mathfrak{A}_K^r$
	is the natural projection map.

	\begin{theorem}\label{Cornut-Vatsal}
		Suppose that $B$ is split at $\ell$.
		Let $\mathcal{R}$
		be a finite subset of
		$\mathbf{T}(\mathbb{A}_f)$
		such that the $\ell$-component of its elements are pairwise distinct modulo
		$
		\mathbf{T}(\mathbb{Q})Z(\mathbf{H}(\mathbb{Q}_\ell))$
		and
		$S$ be the image of a $\mathbf{H}(\mathbb{Q}_\ell)$-orbit
		through the projection
		$\mathbf{T}(\Q)\backslash \mathbf{H}(\mathbb{A}_f)
		\rightarrow
		\mathbf{T}(\Q)\backslash \mathbf{H}(\mathbb{A}_f)/K$.
		Let $\mathcal{G}$
		be an open subgroup of
		$\mathbf{T}(\Q)\backslash \mathbf{T}(\mathbb{A}_f)$.
		Then for all but finitely many
		$s\in S$,
		one has
		\[
		\mathfrak{R}_K^r(\mathcal{G}s)
		=
		(\mathfrak{A}_K^r)^{-1}
		\left(
		\mathfrak{A}_K^r\circ\mathfrak{R}_K^r(\mathcal{G}s)
		\right).
		\]
	\end{theorem}
	
	\begin{remark}\rm
		\begin{enumerate}
			\item
			There are many applications of the above result to number theory: (1)
			Mazur's conjecture (\cite{Vatsal2002});
			(2) the non-vanishing modulo $p$
			(with $p\ne\ell$)
			of certain Yoshida lifts
			(3) the non-vanishing of certain Euler systems
			(in the domain of Iwasawa theory
			as in \cite{Howard2004}). So it is of extreme importance to generalize the above theorem to more other groups.

			\item 
			In applications,
			the set $S$ is often taken to be the image of a
			$\mathbf{H}(\mathbb{Q}_\ell)$-orbit
			containing Heegner points of $\ell$-power conductors.
			In other words,
			it suffices to consider those $S$
			that are the image of a subset of a
			$\mathbf{H}(\mathbb{Q}_\ell)$-orbit
			and are infinite.
			In fact,
			if $H$ is split at $\ell$,
			the Heegner points of conductor $\ell^k$
			are represented by
			matrices of the form
			$
			\begin{pmatrix}
				\ell^k & 1 \\
				0 & 1
			\end{pmatrix}
			\in
			\mathrm{GL}_2(\mathbb{Q}_\ell)
			\simeq
			\mathbf{H}(\mathbb{Q}_\ell)$.
			Multiplying this matrix on the left by
			$\mathrm{diag}(\ell^{-k},1)
			\in
			\mathbf{T}(\mathbb{Q}_\ell)
			\subset
			\mathbf{T}(\mathbb{A}_f)$,
			what one actually considers in the proof of
			the above theorem is
			that for
			$k\gg0$,
			the last identity in the above theorem
			holds with
			$s$ replaced by
			the element
			$
			\begin{pmatrix}
				1 & \ell^{-k} \\
				0 & 1
			\end{pmatrix}
			\in
			\mathrm{GL}_2(\mathbb{Q}_\ell)
			\simeq
			\mathbf{H}(\mathbb{Q}_\ell)
			\subset
			\mathbf{H}(\mathbb{A}_f)
			$.
		\end{enumerate}

	\end{remark}

	Our goal in the appendix is to provide a generalization of the above result to a certain class of algebraic groups, which include as examples compact forms of classical groups. For simplicity, in this introduction we take $\mathbf{H}$ to be a non-abelian general spin group over $\Q$ (for application in the first part of the introduction, we can simply take $\mathbf{H}$ to be the general spin group over $\mathbb{Q}$ associated to the quadratic space $(U,\langle-,-\rangle_{U})$). We assume that $\mathbf{H}$ and $\mathbf{T}$ are split at $\ell$. In these cases, the most important property that we need on $\mathbf{H}$ and $\mathbf{T}$ is that there are subgroups $\mathbf{H}_1,\mathbf{H}_2,\cdots,\mathbf{H}_n$ of $\mathbf{H}$
	such that each $\mathbf{H}_j$ is split at $\ell$ and is isomorphic to $B_j^\times$ where $B_j$ is a (definite) quaternion algebra over $\Q$. We fix isomorphisms $\mathbf{H}_j(\Q_\ell)\simeq\mathrm{GL}_2(\Q_\ell)$. Then we have $n$ rank one unipotent subgroups $\mathbf{U}_1,\cdots,\mathbf{U}_n$
	of $\mathbf{H}(\Q_\ell)$ with $\mathbf{U}_j\subset \mathbf{H}_j(\Q_\ell)$.
	Moreover we can choose these $\mathbf{H}_j$	such that	$\mathbf{H}_j^1(\Q_\ell)$ generate $\mathbf{H}^1(\Q_\ell)$, $\mathbf{T}(\Q_\ell)$ normalizes each $\mathbf{H}_j(\Q_\ell)$ and $\mathbf{U}_j$, and all $\mathbf{U}_j$ commute with each other.
	Let $\mathbf{T}(\Q_\ell)$ act on $\mathbf{U}_j$ by conjugation via a character $\chi_j$. We can choose these $\Hbf_j$ such that these characters $\chi_j$ ($j=1,\cdots,n$) are linearly independent in the character group $X^\ast(\Tbf)$ of $\Tbf$. Fix isomorphisms $\Q_\ell\simeq \mathbf{U}_j$ sending $t$ to
	$u_j(t)$. In this appendix, we prove the following (see also Theorem \ref{density of torus orbits}):
	\begin{theorem}\label{Cornut-Vatsal's theorem}
		Let $\mathcal{R}=\{g_1,\cdots,g_r\}\subset\mathbf{T}(\mathbb{A}_f)$	be a finite subset of $r$ elements such that
		\begin{equation}\label{pairwise distinct}
			(g_k)_\ell(g_i)_\ell^{-1}
			\notin
			\mathbf{T}(\Q_\ell)
			\cap
			\big(
			\cup_{j=1}^n
			\mathbf{H}_j(\Q)Z(\mathbf{H}_j(\Q_\ell))
			\big),
			\quad
			\forall
			i\ne k.
		\end{equation}
		Let $\mathcal{G}$ be an open subgroup of $\mathbf{T}(\Q)\backslash \mathbf{T}(\mathbb{A}_f)$. For an $n$-tuple of integers	$\underline{N}=(N_1,\cdots,N_n)$, we write $u(\ell^{-\underline{N}})
		=\prod_{i=1}^nu_i(\ell^{-N_i})$. Then for $N_1,N_2,\cdots,N_n\gg0$,
		we have
		\[
		\mathfrak{R}_K^r
		\left(
		\mathcal{G}
		u(\ell^{-\underline{N}})
		\right)
		=
		(\mathfrak{A}_K^r)^{-1}
		\left(
		\mathfrak{A}_K^r\circ
		\mathfrak{R}_K^r
		\left(
		\mathcal{G}
		u(\ell^{-\underline{N}})
		\right)
		\right).
		\]
	\end{theorem}

	Now let's briefly indicate the ideas behind of the proof of the above result: the proof borrows ideas from the proof of Theorem \ref{Cornut-Vatsal}. So maybe it is useful to first indicate the proof of this theorem for the case $\mathbf{H}$ split at $\ell$. Let $h\in \mathbf{T}(\mathbb{A}_f)$, and
	$\Delta\colon\mathbf{H}^1(\Q_\ell)\rightarrow\mathbf{H}^1(\Q_\ell)^r$
	be the diagonal map.
	\begin{enumerate}[label=(\alph*)]
		\item
		Ratner's orbit closure theorem says that the closure of $\prod_{i=1}^r\Gamma_{K}(h g_i)\Delta(U)$ in $\mathbf{H}^1(\Q_\ell)^r$
		contains $c\Delta(\mathbf{H}^1(\Q_\ell))c^{-1}$	for some $c\in U^r$,
		where $\Gamma_{K}(h g_i)$ is a lattice in $\mathbf{H}^1(\Q_\ell)$ (see §\ref{adelic formulation} for the definition). Moreover one can show that for any $h\in \mathbf{T}(\mathbb{A}_f)$, the lattices $\Gamma_K(h g_i)$
		and $\Gamma_K(h g_k)$ are \textit{not} $U$-commensurable
		for any $i\ne k$, which implies that the closure of $\prod_i\Gamma_K(h g_i)\Delta(\mathbf{U}_i)$ is $\mathbf{H}^1(\Q_\ell)^r$.

		\item
		Ratner's uniform distribution theorem on unipotent orbits says that
		for any locally constant function $f$ on $(\mathbf{H}(\Q)\backslash
		\mathbf{H}(\mathbb{A}_f)/K)^r$ and any compact open subset $\kappa$ of $\Q_\ell$, one has
		\[
		\lim\limits_{N\rightarrow+\infty}
		\frac{1}{\lambda(\kappa_N)}
		\int_{\kappa_N}
		f(\mathfrak{R}_K^r(h g_iu_i(t)))dt
		=
		\int_{(\mathbf{H}(\Q)\backslash \mathbf{H}(\mathbb{A}_f)/K)^r}
		f
		d\mu
		\]
		Here $\kappa_N=\kappa/\ell^N$. Now one integrates both sides over
		$h\in\mathcal{G}$ and then take $f$ to be certain characteristic functions,
		one concludes the proof.
	\end{enumerate}

	We next explain the substantial difficulties in these two steps when generalizing
	Theorem \ref{Cornut-Vatsal} to Theorem \ref{Cornut-Vatsal's theorem}. In (a) one should find suitable conditions on $\mathbf{H}$ such that $\prod_{i=1}^r\Gamma_K(h g_i)\Delta(\mathbf{U}_1\cdots \mathbf{U}_n)$
	is dense in $\mathbf{H}^1(\Q_\ell)^r$. First consider the case $r=1$. For the case $\mathbf{H}$ a compact general spin group over $\Q$, the closure
	$\Gamma_K(h g_i)\mathbf{U}_j$ inside $\mathbf{H}^1(\Q_\ell)$ should contain the closure of $(\Gamma_K(h g_i)\cap \mathbf{H}_j(\Q))\mathbf{U}_j$,
	the latter closure contains $\mathbf{H}_j^1(\Q_\ell)$. Since these	$\mathbf{H}_j^1(\Q_\ell)$ generate $\mathbf{H}^1(\Q_\ell)$, we conclude that
	$\Gamma_K(h g_i)\mathbf{U}_1\mathbf{U}_2\cdots \mathbf{U}_n$
	should be dense in the whole group $\mathbf{H}^1(\Q_\ell)$. For the case
	$r>1$, we use the arguments in \cite{CornutVatsal2005} to show that again
	$\Gamma_K(h g_i)\cap \mathbf{H}_j(\Q)$ and $\Gamma_K(h g_k)\cap \mathbf{H}_j(\Q)$ are not $\mathbf{U}_i$-commensurable for any $j$
	and $i\ne k$. Therefore the closure of $\prod_{i=1}^r\Gamma_K(h g_i)
	\Delta(\mathbf{U}_j)$ inside $\mathbf{H}^1(\Q_\ell)^r$ contains
	$\mathbf{H}_j^1(\Q_\ell)^r$. Now these $\mathbf{H}_j^1(\Q_\ell)^r$
	generate $\mathbf{H}^1(\Q_\ell)^r$, one sees immediately that
	$\prod_{i=1}^r\Gamma_K(h g_i)\Delta(\mathbf{U}_1\cdots \mathbf{U}_n)$
	is dense in $\mathbf{H}^1(\Q_\ell)^r$.

	In (b), one should choose
	a set $S$ as in
	Theorem \ref{Cornut-Vatsal}.
	As we are in the case $\mathbf{H}$ split at $\ell$,
	so it is natural to consider a set $S$ of unipotent elements.
	However even this is not enough.
	In fact, we should require that
	$S$ is contained in a unipotent subgroup whose rank is no greater than that of
	$\mathbf{T}(\Q_\ell)\cap \mathbf{H}^1(\Q_\ell)$.
	Indeed,
	in (2) for our case
	one needs an identity of the form
	\[
	\mathcal{G}\prod_{i=1}^nu_i(\kappa_{N_i})K
	=
	\mathcal{G}\prod_{i=1}^nu_i(\ell^{-N_i})K
	\]
	for a sufficiently small compact open subgroup
	$\kappa$ of $\Q_\ell^\times$
	(here 
	$\kappa_{N_i}=\kappa\ell^{-N_i}$).
	One way to absorb
	$\kappa_{N_i}$
	into the groups $\mathcal{G}$ and $K$
	is to use the non-trivial conjugate action of
	$\mathcal{G}$ on
	$\mathbf{U}_1,\cdots,\mathbf{U}_n$.
	For $\kappa$ sufficiently small,
	as long as
	the conjugate action of
	$\mathcal{G}$ on each $\mathbf{U}_i$
	is non-trivial
	(in particular, the rank of
	$\mathbf{T}(\Q_\ell)\cap \mathbf{H}^1(\Q_\ell)$
	should be no less than $n$),
	the above identity can be achieved.
	On the other hand,
	$S$ should not be too small,
	in fact,
	the limit in (2) for our case is over
	$N_1,N_2,\cdots,N_n\rightarrow+\infty$.
	Thus
	we take
	\[
	S=
	\left\{
	u_1(\ell^{-N_1})\cdots u_n(\ell^{-N_n})|N_1,\cdots,N_n\in\mathbb{N}
	\right\}.
	\]

	From this theorem, one deduces quite easily the following application to automorphic forms (see also Theorem \ref{automorphic application of appendix} for weaker assumptions on $\mathcal{G}$):
	\begin{theorem}\label{main theorem-2}
		Take $\mathcal{G}=\mathbf{T}(\mathbb{Q})\backslash\mathbf{T}(\mathbb{A}_f)$.
		Let $A$ be a ring and $M$ an $A$-module.
		Let $\{\beta_i\}_{i=1}^r$ be a set of elements in $A$ with $\beta_1\in A^\times$. Consider a map $f\colon \mathbf{H}(\mathbb{Q})\backslash\mathbf{H}(\mathbb{A}_f)/K\rightarrow M$ which is not $\mathbf{H}^1(\mathbb{A}_f)$-invariant (under right translation). Then for any $N_1,\cdots,N_n\gg0$, there exists $h=h_{\underline{N}}\in\mathbf{T}(\mathbb{A}_f)$ such that
		\[
		\sum_{i=1}^r\beta_if(hg_iu(\ell^{-\underline{N}}))\neq0.
		\]
	\end{theorem}

	Both Theorems \ref{main result for special case} and \ref{main theorem-2} can be generalized without difficulties to a totally real number field instead of $\Q$. However, we choose to work over $\Q$ to make the ideas in our work more transparent.

	\subsection*{Outline}
	In \S \ref{automorphic forms} and \S \ref{theta lift}
	we recollect the basic notions on automorphic forms and theta lifts that we will
	need in this article.
	Moreover we define the notions of
	$p$-integral and $p$-primitive automorphic forms on orthogonal group and metaplectic
	group,
	which are the basis for our work.
	Then in \S \ref{Bruhat-Schwartz function}
	we define carefully the Bruhat-Schwartz function 
	used in global theta correspondences
	and deduce some simple properties of the theta lifts of
	automorphic forms from $H$ to $\widetilde{G}$ by this
	particular Bruhat-Schwartz function.
	In \S \ref{toric periods and toric intergals}
	we consider toric periods and toric integrals
	and in
	\S \ref{non-vanishing of theta lift mod p}
	we deduce the $p$-primitivity of theta lifts from the non-vanishing modulo $\mathfrak{P}$ of
	toric periods.
	\S \ref{non-vanishing of Besse mod p} is the technical heart of this article.
	We make a digression in
	\S \ref{Equidistribution of CM points}
	and recall the main results from
	Appendix \ref{appendix}
	that will be used in the next subsection.
	In \S \ref{non-vanishing of Besse mod p-subsection},
	we prove the main technical result of this article on the non-vanishing modulo $\mathfrak{P}$
	of the toric periods.
	In Appendix \ref{appendix},
	we give a quite general treatment of the crucial ingredient mentioned in
	\S \ref{Equidistribution of CM points}
	and we believe it may have interest on its own.
	As the reader can see when reading this article,
	the works \cite{CornutVatsal2005,Jia2010,HsiehNamikawa2017}
	have a great influence on the ideas and presentations of
	this article.

	\subsection*{Notations}
	\begin{enumerate}
		\item 
		Write
		$\mathbb{A}$
		for the ring of adèles of $\Q$,
		$\mathbb{A}_f$ for the ring of finite adèles. We fix a prime number
		$p$ in this article
		and
		an isomorphism of fields $\mathbb{C}
		\simeq
		\mathbb{C}_p$,
		compatible with the embeddings
		$\overline{\mathbb{Q}}
		\hookrightarrow
		\mathbb{C}$
		and
		$\overline{\mathbb{Q}}
		\hookrightarrow
		\mathbb{C}_p$.
		We write
		$\mathcal{O}$
		for the ring of integers of
		$\mathbb{C}_p$,
		$\mathfrak{P}$
		the maximal ideal of
		$\mathcal{O}$
		and
		$\kappa=\mathcal{O}/\mathfrak{P}$
		the residue field.
		Using the isomorphism $\C\simeq\C_p$,
		we will identify
		$\mathbb{C}$ with $\mathbb{C}_p$
		when they appear as the values of certain maps
		or the base field of certain spaces of algebraic representations
		of a group.
		In particular we will view
		$\mathcal{O}$
		as a subring of $\mathbb{C}$.

		\item 
		We write
		$\mathrm{Sym}_{n}(R)$
		for the set of
		$n\times n$ symmetric matrices
		with entries in a ring $R$,
		$\mathrm{M}_{m,n}(R)$
		the set of
		$m\times n$-matrices with entries in $R$. If $m=n$, we simply write $\mathrm{M}_n(R)$ for $\mathrm{M}_{n,n}(R)$.
		For a matrix 
		$A$,
		we write
		$A^\mathrm{t}$
		for its transpose.
		For an abelian group
		$M$
		and an ring $R$,
		we write
		\[
		M(R)=M\otimes_{\Z}R.
		\]

		\item 
		For an algebraic group $G$
		over $\Q$,
		we write
		\[
		[G]=G(\Q)\backslash G(\mathbb{A}),
		\quad
		[G]_f=G(\Q)\backslash G(\mathbb{A}_f).
		\]

		\item 
		For a (abstract) group $G$ acting on a set $X$, we write $G_x$ for the stabilizer of an element $x\in X$ in $G$
		\[
		G_x=\{g\in G|gx=x\}.
		\]

		\item 
		For a function $f$ (on a set $X$) valued in an $\mathcal{O}$-module $M$, we write
		\[
		\overline{f}=f(\mathrm{mod}\,\mathfrak{P}).
		\]
	\end{enumerate}

	\section*{Acknowledgments}
	We would like to thank Massimo Bertolini, Christophe Cornut, Haruzo Hida, Ming-Lun Hsieh and Jacques Tilouine for a lot of useful and stimulating discussions concerning earlier drafts of this article. We thank Zheng Liu in particular for careful reading which uncovers an error in an early version of the article. At last we would like to express our thanks to the referee for his/her many useful comments, suggestions and questions, which help improve a lot the structure and the readability of this article.

	\section{Automorphic forms}\label{automorphic forms}
	In this section we define
	vector-valued automorphic forms
	on        
	orthogonal group
	and
	metaplectic group.
	This is well-known in the literature and we omit the proof of some facts.
	Along the way we shall also define the notions of $p$-integral and
	$p$-primitive automorphic forms,
	the basic objects of study in this article.

	\subsection{The groups $H$ and $G$}
	\subsubsection{Orthogonal group}\label{orthogonal group}
	Let
	$(U,\langle-,-\rangle_U)$ be a non-degenerate quadratic module over $\Z$ free of \emph{odd} rank $2n+1\ge3$.
	We fix a $\Z$-basis
	\[
	\mathfrak{B}=(E_1,E_2,\cdots,E_{2n+1})
	\]
	for
	$U$ and write
	$Q_U$ for the symmetric matrix
	associated to the quadratic form on $U$
	under this basis. We assume that $Q_U$ is of the form
	\[
	Q_U
	=
	\mathrm{diag}
	(\delta_1,\delta_2,\cdots,\delta_{2n+1})
	\]
	with
	$\delta_i$ positive integers.
	For $i=1,2,\cdots,2n+1$, we define the algebraic numbers
	\[
	d_i=\frac{1}{\sqrt{2(-1)^{i-1}\delta_i}}\in\overline{\Q}.
	\]

	In the following, we write
	\[
	R=\Z[d_1,\cdots,d_{2n+1}].
	\]
	Then the quadratic space
	$(U(R),\langle-,-\rangle_{U(R)})$
	is split and
	we have another basis for $U(R)$:
	\[
	\widetilde{\mathfrak{B}}
	=
	(\widetilde{E}_1,\cdots,\widetilde{E}_{2n+1}),
	\]
	where each basis vector is given by
	\[
	\widetilde{E}_i=
	\begin{cases*}
		d_{2i-1}E_{2i-1}+d_{2i}E_{2i},
		&
		$i=1,\cdots,n$;
		\\
		d_{2(i-n)-1}E_{2{i-n}-1}-d_{2(i-n)}E_{2(i-n)},
		&
		$i=n+1,\cdots,2n$;
		\\
		E_{2n+1},
		&
		$i=2n+1$.
	\end{cases*}
	\]
	Under this basis, the quadratic form $\langle-,-\rangle_{U(R)}$ is represented by the matrix
	\[
	\widetilde{Q}_U
	=
	\begin{pmatrix}
		0 & 1_{n} & 0 \\
		1_{n} & 0 & 0 \\
		0 & 0 & 1
	\end{pmatrix}.
	\]
	The transformation matrix from $\mathfrak{B}$ to $\widetilde{\mathfrak{B}}$ is given by
	\begin{equation}\label{T}
		\mathcal{T}
		=
		\begin{pmatrix}
			d_1 & d_2 & 0 & 0 & \cdots & 0\\
			0 & 0 & d_3 & d_4 & \cdots  & 0\\
			\vdots & \vdots & \vdots & \vdots & \ddots & 0\\
			d_1 & -d_2 & 0 & 0 & \cdots & 0\\
			0 & 0 & d_3 & -d_4  & \cdots & 0\\
			\vdots & \vdots & \vdots & \vdots & \ddots & 0 \\
			0 & \cdots & \cdots & \cdots & \cdots & 1 
		\end{pmatrix}.
		\index{T@$\mathcal{T}$}
	\end{equation}

	The orthogonal group scheme
	$H=\mathrm{O}_U$
	over $\Z$
	consists of
	$g\in\mathrm{GL}_{n}$
	such that
	$g^\mathrm{t}Q_Ug=Q_U$. We write
	\[
	H^1=\mathrm{SO}_U
	\]
	for the special orthogonal subgroup of $H$ and we have an isomorphism of group schemes over $\Z$
	\[
	H=H^1\times\{\pm1_{2n+1}\}
	\simeq H^1\times\{\pm1\}.
	\]
	We have the following torus subgroups of
	$ H$:
	\begin{align*}
		T_i
		&
		=
		\mathrm{SO}_{\Z(E_{2i-1},E_{2i})}
		\quad
		(i=1,\cdots,n)
		\\
		T
		&
		=
		\prod_{i=1}^{n}
		T_i,
	\end{align*}
	where
	$\Z(E_{2i-1},E_{2i})$
	is the quadratic submodule of $U$
	generated by the basis vectors
	$E_{2i-1},E_{2i}$.
	Then $T$ is a maximal torus of $ H$. We define group isomorphisms
	\[
	\mu_i\colon
	T_i(R)
	\rightarrow
	R^\times,
	\quad
	\begin{pmatrix}
		a & bd_{2i}^2/d_{2i-1}^2\\
		b & a
	\end{pmatrix}
	\mapsto
	a+bd_{2i}/d_{2i-1}.\index{m@$\mu_i$}
	\]

	Now we consider algebraic representations of $H^1_{/R}$.
	We write $B$ for the Borel subgroup of $H^1_{/R}$ containing $T_{/R}$, consisting of elements of the form $\mathcal{T}\begin{pmatrix}
		a & b & ac \\ 0 & a^{-\mathrm{t}} & 0 \\ 0 & -c^{\mathrm{t}} & 1
	\end{pmatrix}\mathcal{T}^{-1}$ with $a\in\GLmr_{n/R}$ upper triangular. We then write $B^-$ for the Borel subgroup of $H_{/R}$ opposite to $B$. Then there is a unique element $w_0$ of maximal length in the Weyl group of $T_{/R}$ such that $w_0(B)=B^-$.
	We identify the character group $X^\ast(T_{/R})$ with $\Z^n$:
	\[
	\Z^n\to X^\ast(T_{/R}),
	\quad
	\lambda=(\lambda_1,\cdots,\lambda_n)
	\mapsto
	\left(
	\prod_{i=1}^nT_i\ni t=t(t_1,\cdots,t_n)\mapsto\prod_{i=1}^n\mu_i(t_i)^{\lambda_i}
	\right).
	\]
	Note that the Weyl element $w_0$ also acts on $X^\ast(T_{/R})$, sending $(\lambda_1,\cdots,\lambda_n)$ to $(-\lambda_1,\cdots,-\lambda_n)$.
	A character $\lambda\in\Z^n$ is \emph{dominant} (with respect to $(T_{/R},B)$) if $\lambda_1\ge\lambda_2\ge\cdots\ge\lambda_n\ge0$.

	Each $\lambda\in\Z^n$ induces a character of $B$. Then we have the induction algebraic representation $(\rho_{\lambda},\M_{\lambda})$ of $H^1_{/R}$:
	\[
	\M_{\lambda}
	=
	\left\{
	f\in R[H^1_{/R}]\Big|
	f(bg)=(w_0\lambda)(b)f(g),\,
	\forall
	R\text{-algebra }R',\,g\in H(R'),b\in B(R')
	\right\}.\index{r@$(\rho_{\lambda},\M_{\lambda})$}
	\]
	Here $H^1_{/R}$ acts by right translation. It is well known that this module is non-zero if and only if $\lambda$ is dominant (\cite[II.2.6]{Jantzen2003}). If $\widetilde{R}$ is an $R$-algebra which is a field of characteristic $0$ and $\lambda$ is dominant, then $\M_{\lambda}(\widetilde{R})$ is the irreducible algebraic representation of $H(\widetilde{R})$ of highest weight $\lambda$.
	
	Since $H\simeq H^1\times\{\pm1\}$, we list the irreducible algebraic representations $(\rho,\M)$ of $H_{/R}$ by their highest weights
	\[
	\lambda=(\lambda_1,\cdots,\lambda_{n};\varepsilon)\index{l@$\lambda$}
	\]
	Here $(\lambda_1,\cdots,\lambda_n)$ is a dominant character of $T_{/R}$ and $\varepsilon=\pm1$. The representation is defined as follows:
	$(\rho,\M)$ is isomorphic to $(\rho_{(\lambda_1,\cdots,\lambda_n)},\M_{(\lambda_1,\cdots,\lambda_n)})$ as representations of $H^1_{/R}$. Moreover, $-1_{2n+1}\in H$ acts on $\M$ by the scalar $\varepsilon$.
	In the following, we again denote this representation $(\rho,\M)$ of $H_{/R}$ by $(\rho_{\lambda},\M_{\lambda})$.

	For applications in theta correspondence, it is often very useful to have an explicit description of the representation $\rho_{\lambda}$, which we will do now: we set
	\[
	R'=R
	\Big[
	\frac{1}{(n-1)!},\frac{1}{(\lambda_1-\lambda_2)!},\cdots,\frac{1}{(\lambda_{n-1}-\lambda_n)!},\frac{1}{\lambda_n!}
	\Big].\index{R@$R'$}
	\]
	Write $(\rho_{\mathrm{std}},\M_\mathrm{std}=U_{/R'})$
	for the standard representation
	of
	$H_{/R'}$
	\[
	\rho_{\mathrm{std}}
	\colon
	H_{/R'}\to\GLmr(\M_{\mathrm{std}})\simeq\mathrm{GL}_{n/R'}.
	\]
	Then we define the following representation of $H_{/R'}$
	\[
	\widetilde{\M}_{\lambda }
	:=
	\bigotimes_{i=1}^{n-1}
	\Symmr^{\lambda_i-\lambda_{i+1}}(\wedge^i\M_{\mathrm{std}}).
	\]
	So $\M_\lambda$ is the irreducible subrepresentation of $\widetilde{\M}_\lambda$ generated by the highest weight vector (of weight $\lambda$)
	\[
	v_{\lambda }
	:=
	\bigotimes_{i=1}^{n-1}(\widetilde{E}_1\wedge\cdots\wedge\widetilde{E}_i)^{\lambda_i-\lambda_{i+1}}.
	\]
	Here $v_1v_2:=
	\frac{1}{2!}(v_1\otimes v_2+v_2\otimes v_1)$,
	$v_1\wedge v_2:=
	\frac{1}{2!}(v_1\otimes v_2-v_2\otimes v_1)$
	and similarly for
	$v_1v_2\cdots v_k,v_1\wedge v_2\wedge\cdots\wedge v_k$ (due to our assumptions on $R'$, these vectors are well-defined over $R'$).
	This vector is an eigenvector for the maximal torus $T_{/R'}$
	\[
	\rho_{\lambda }(t)v_{\lambda }
	=
	\prod_{i=1}^{n}
	\mu_i(t_i)^{\lambda_i }
	v_{\lambda },
	\quad
	\forall
	t=\mathrm{diag}(t_1,t_2,\cdots,t_{n})\in
	T_{/R'}.
	\]

	We impose the following condition
	\begin{taggedtheorem}{(C1)}\label{condition on lambda}
		The dimension $\dim_{\C}\M_{\lambda}(\C)>1$\index{C@(C1)}

	\end{taggedtheorem}

	\begin{remark}\label{C1 ensure non-vanishing of theta lifting}\rm
		This condition will ensure that the theta lifting of the discrete series of highest weight $\lambda$ of $H(\R)$ to $\widetilde{G}(\R)$
		is non-zero (\cite[Theorem 6.9, p.25]{KashiwaraVergne1978}): indeed, according to \emph{loc.cit}, if either $\varepsilon=(-1)^{\sum_i\lambda_i}$ or $\varepsilon=-(-1)^{\sum_i\lambda_i}$ and $1\leq r\leq n$ (here $r$ is the maximal integer such that $\lambda_r>0$), then the theta lifting is non-zero (with an explicit formula for the highest weight of the theta lifting in terms of the weight $\lambda$). The condition $\dim_{\C}\mathbb{M}_{\lambda}(\C)>1$ implies in particular that $\lambda_1\neq0$ and thus $r\ge1$.

	\end{remark}

	\subsubsection{Metaplectic group}\label{symplectic group}
	Let $(V,\langle-,-\rangle_V)$ be a non-degenerate symplectic module over $\Z$ free of rank $4n$. Fix a symplectic $\Z$-basis for $V$
	\[
	V=\Z(e_1^+,\cdots,e_{2n}^+,e_1^-,\cdots,e_{2n}^-).
	\]
	We assume that $\langle e_i^+,e_j^-\rangle_V=\delta_{i,j}$ for all $i,j=1,\cdots,2n$. Then the matrix corresponding to $\langle-,-\rangle$ is given by
	\[
	Q_V
	=
	\begin{pmatrix}
		0 & 1_{2n} \\
		-1_{2n} & 0
	\end{pmatrix}.
	\]
	Then we write $G=\Spmr_V$ to be the symplectic group scheme	associated to
	$V$ over $\Z$, consisting of $g\in\mathrm{SL}_{4n}/\Z$ such that
	$g^\mathrm{t}Q_Vg=Q_V$.

	We also write
	\[
	V^\pm=\Z(e_1^\pm,\cdots,e_{2n}^\pm),
	\]
	both are maximal isotropic submodules of $V$.

	For each place $v$ of $\Q$, there is a unique non-trivial two-fold cover $\widetilde{G}(\Q_v)$ of $G(\Q_v)$, the metaplectic group associated to $G(\Q_v)$:
	\[
	1
	\rightarrow
	\{\pm1\}
	\rightarrow
	\widetilde{G}(\Q_v)
	\rightarrow
	G(\Q_v)
	\rightarrow
	1.
	\]
	We write elements in $\widetilde{G}(\Q_v)$ as pairs $(g,\varepsilon)$ with $g\in G(\Q_v)$ and $\varepsilon=\pm1$. This central extension splits over any Siegel parabolic subgroup of $G(\Q_v)$ and splits in a unique way over any unipotent subgroup of $G(\Q_v)$. It also splits over $G(\Z_\ell)$ if $v=\ell\neq 2$.
	If we view $\GLmr_{2n}(\Q_v)$ as a subgroup of $G(\Q_v)$ (via the map $g\mapsto\mathrm{diag}(g,g^{-\mathrm{t}})$), then we have
	\[
	(g_1,\varepsilon_1)(g_2,\varepsilon_2)
	=
	(g_1g_2,\varepsilon_1\varepsilon_2\cdot(\mathrm{det}(g_1),\mathrm{det}(g_2))_v).
	\]
	Here $(-,-)_v$ is the Hilbert symbol on $\Q_v$. In particular, this central extension splits over $\mathrm{SL}_{2n}(\Q_v)\subset\mathrm{GL}_{2n}(\Q_v)$.

	We define the adelic metaplectic group $\widetilde{G}(\mathbb{A})$ to be the quotient
	\[
	\widetilde{G}(\mathbb{A})
	=
	\prod_v\,'\widetilde{G}(\Q_v)\Big/
	\Big\{(\varepsilon_v=\pm1)_v\mid\prod_v\varepsilon_v=1
	\Big\}.
	\]
	In the following, the component $\varepsilon=(\varepsilon_v)_v$ in an element $(g_v,\varepsilon_v)_v\in\widetilde{G}(\mathbb{A})$ will be identified with the product $\prod_{v}\varepsilon_v=\pm1$ (so that we can and will view $\varepsilon$ as a scalar ($\pm1$)).
	It is well-known that the two-fold cover $\widetilde{G}(\mathbb{A})\rightarrow G(\mathbb{A})$ splits over $G(\Q)$. Moreover, for any algebraic subgroup $H$ of $G$, we write $\widetilde{H}(\mathbb{A})$ for the preimage of $H(\mathbb{A})$, resp. $H$, $H$ by this cover map; similarly, for any subgroup $H$ of $G(\Q_v)$, resp. $G(\A)$, we write $\widetilde{H}$ for the preimage of $H$ in $\widetilde{G}(\Q_v)$, resp. $\widetilde{G}(\A)$ by this cover map. We then write
	\[
	[\widetilde{G}]=G(\Q)\backslash\widetilde{G}(\mathbb{A}).
	\]
	The highest weight of an irreducible algebraic representation $(\rho_{\tau},\M_\tau)$ of $\GLmr_{2n/\Z}$ with respect to the subgroup of upper triangular matrices is given by an $2n$-tuple
	\[
	\tau=(\tau_1,\cdots,\tau_{2n})
	\in\mathbb{Z}^{2n}.
	\]
	The construction of $\M_{\tau}$ is similar to $\M_{\lambda}$ as in the case of special orthogonal group $H^1$.

	Recall that we have fixed $(\rho_{\lambda},\M_\lambda)$. We fix also an embedding $R'\hookrightarrow\C$ of rings. We assume moreover that $R'$ embeds into $\C$.
	In the following we will take (\cite[p.27]{KashiwaraVergne1978})
	\[
	\tau=
	\begin{cases*}
		(\lambda _1,\cdots,\lambda _{n},0,\cdots,0),
		&
		$\varepsilon=(-1)^{\sum_i\lambda_i}$;
		\\
		(\lambda_1,\cdots,\lambda_{r},1,1,\cdots,1,
		\overbrace{0,\cdots,0}^\text{$(r-1)$-terms}),
		&
		$\varepsilon=-(-1)^{\sum_i\lambda_i}$ and $1\leq r\leq n$.
	\end{cases*}
	\index{t@$\tau$}
	\]
	Here $r$ is the maximal integer such that $\lambda_r>0$. Note that the condition $r\ge1$ is always satisfied due to \ref{condition on lambda} (\emph{loc.cit} Remark \ref{C1 ensure non-vanishing of theta lifting}).
	We denote again $(\rho_{\tau},\M_{\tau})$ the base change $(\rho_{\tau},\M_{\tau})\otimes_{\Z}R'$, which is now an algebraic representation of $\GLmr_{2n/R'}$.
	We give an explicit description of $(\rho_{\tau},\M_{\tau})$ as follows. We first digress to recall the notion of pluri-harmonic polynomials from \cite[II.(5.2)]{KashiwaraVergne1978}: we write
	\begin{align*}
		W
		&
		=U\otimes V,
		\\
		W^\pm
		&
		=U\otimes V^\pm.
	\end{align*}
	We equip
	$W$
	with the natural symplectic pairing
	$\langle-,-\rangle_W$
	induced from
	$U$ and $V$:
	\[
	\langle u\otimes v,u'\otimes v'\rangle_{W}
	:=
	\langle u,u'\rangle_U\cdot\langle v,v'\rangle_{V},
	\quad
	u,u'\in U,\,
	v,v'\in V.
	\]
	Moreover, we identify $W$, resp. $W^\pm$ with $\mathrm{M}_{2n+1,4n}(\mathbb{Z})$, resp. $\mathrm{M}_{2n+1,2n}(\mathbb{Z})$ via our fixed basis and we will write elements $W,W^\pm$ as matrices.
	Then a polynomial function $f(w)$ on $w=(w_{i,j})_{i,j}\in W^+(\C)$ is \emph{pluri-harmonic} if
	\[
	\sum_{s=1}^{2n+1}\frac{\partial^2}{\delta_s\partial w_{s,i}\partial w_{s,j}}f(w)=0,
	\quad
	\forall
	1\leq i,j\leq 2n
	\]
	Write $\mathfrak{H} $ for the space of all pluri-harmonic polynomials on $W^+(\C)$. Then $H(\C)\times\GLmr_{2n}(\C)$ acts on $\mathfrak{H}$ in a natural way. Then we write $\mathfrak{H}_{\lambda}$ for its subspace of $\lambda $-isotypic component. According to II.(6.4) of \textit{loc.cit}, there is an isomorphism of representations of $H(\mathbb{C})\times\mathrm{GL}_{2n}(\mathbb{C})$:
	\begin{equation}\label{decomposition of pluri-harmonic}
		\mathfrak{H}_{\lambda } 
		\simeq
		\M_{\lambda }(\C)\otimes\M_{\tau }(\C).
	\end{equation}
	Write $\mathfrak{H} (\lambda )$ for the space of pluri-harmonic $\M_{\lambda }$-valued polynomials $f$ on $\mathrm{M}_{2n+1,2n}(\mathbb{C})$ such that $f(wh)=\rho_{\lambda }(h)f(w)$ for any $w$ and $h\in H(\mathbb{C})$. 
	Then we have an isomorphism of representations of $\GLmr_{2n}(\C)$:
	\[
	\M_\tau(\C)\simeq\mathfrak{H}(\lambda).
	\]	
	We fix one such isomorphism and identify these two spaces. We then define
	\[
	\M_\tau(R')
	=
	\left\{
	f\in\mathfrak{H}(\lambda)|f(\mathrm{M}_{2n+1,2n}(R'))\subset\M_{\lambda }(R')
	\right\}.
	\]
	This $R'$-module is indeed stable under the action of $\GLmr_{2n}(R')$.

	We write $(\rho_{\tau^\vee},\M_{\tau^\vee})$ for the dual representation of $(\rho_{\tau},\M_\tau)$.

	\subsection{Automorphic forms}\label{automorphic forms on H and G}
	We review the notion of automorphic forms on
	$H(\mathbb{A})$ and $\widetilde{G}(\mathbb{A})$
	and then
	we will define the notions of $p$-integral and
	$p$-primitive automorphic forms
	on these groups.
	
	Recall we have fixed an irreducible algebraic representation $(\rho_{\lambda},\M_\lambda)$ of $H_{/R'}$ and an irreducible algebraic representation $(\rho_{\tau},\M_\tau)$ of $\GLmr_{2n/R'}$.

	We consider the following conditions on $p$:
	\begin{taggedtheorem}{(C2)}\label{condition on p}
		$p\nmid2\delta_1\cdots\delta_{2n+1}$ and $p>\max(n-1,
		\lambda_1 -\lambda_2 ,\cdots,\lambda_{n-1} -\lambda_{n},
		\lambda_{n})$.\index{C@(C2)}
	\end{taggedtheorem}
	These conditions ensure that there is an embedding $R'\hookrightarrow\mathcal{O}$ (we can choose this embedding to be compatible with $R'\hookrightarrow\C$ via the isomorphism $\C\simeq\C_p$). In particular, $\M_{\lambda}(\mathcal{O})$ is an irreducible algebraic representation of $H(\mathcal{O})$, thus also an irreducible algebraic representation of $H(\Z_p)$. Similar for $\M_{\tau}(\mathcal{O})$.

	\subsubsection{Automorphic forms on $ H(\mathbb{A})$}

	\begin{definition}\label{automorphic forms on SO}
		\begin{enumerate}
			\item 
			For a compact open subgroup $K$ of
			$ H(\mathbb{A}_f)$,
			we write the space of
			automorphic forms on
			$ H(\mathbb{A}_{F})$
			of weight $\lambda$
			and of level $K$ to be
			$\mathcal{A}_{\rho_\lambda}( H,K)$
			\[
			\mathcal{A}_{\rho_\lambda}( H,K)
			=
			\left\{
			\mathbf{f}\colon
			H(\mathbb{A}_{F})
			\rightarrow
			\M_{\lambda}(\C)
			\Big|
			\mathbf{f}(z\gamma hk)
			=
			\rho_\lambda(h_\infty^{-1})\mathbf{f}(h_f),
			\begin{array}{l}
				\text{for any }
				h=h_\infty h_f
				\in
				H(\mathbb{R})\times
				H(\mathbb{A}_f),
				\\
				\text{and }
				(z,\gamma,k)
				\in
				Z_{ H}(\mathbb{A}_{F})
				\times
				H(\Q)
				\times
				K
			\end{array}
			\right\}
			\]
			Here
			$Z_{ H}$
			is the center of
			$ H$.

			\item 
			For
			$\mathbf{f}\in\mathcal{A}_{\rho_{\lambda}}( H,K)$,
			we define its \textbf{$p$-adic avatar} to be\footnote{Here we use the isomorphism $\C\simeq\C_p$ to identify $\M_{\lambda}(\C_p)$ with $\M_{\lambda}(\C)$}
			\[
			\mathbf{f}_p
			\colon
			H(\mathbb{A}_f)
			\rightarrow
			\M_\lambda(\C_p),
			\quad
			h
			\mapsto
			\rho_\lambda(h_p)^{-1}\mathbf{f}(h).\index{f@$\mathbf{f}_p$}
			\]
		\end{enumerate}
	\end{definition}            
	Thus for
	$\mathbf{f}\in
	\mathcal{A}_{\rho_{\lambda}}( H,K)$,
	one has
	\[
	\mathbf{f}_p(\gamma h u)
	=
	\rho_\lambda(u_p)^{-1}\mathbf{f}_p(h),
	\quad
	\forall
	\gamma\in H(\Q),
	h\in H(\mathbb{A}_f),
	u\in K.
	\]
	The values of
	$f_p$
	is determined by its values on a set of representatives of the finite double coset
	$ H(\mathbb{Q})\backslash
	H(\mathbb{A}_f)/
	K$.
	
	\begin{definition}\label{p-integral forms on special orthogonal group}
		Let
		$K=K_pK^p$ be a compact open subgroup of
		$ H(\mathbb{A}_f)$
		such that the
		$p$-th component
		$K_p\subset
		H(\Z_p)$.
		\begin{enumerate}
			\item
			We call
			$\mathbf{f}\in
			\mathcal{A}_{\rho_{\lambda}}( H,K)$
			\textbf{$p$-integral}
			if
			$\mathbf{f}_p(g)\in
			\M_\lambda
			(\mathcal{O})$
			for all $g\in H(\mathbb{A}_f)$.
			
			\item
			Let $\mathbf{f}\in\mathcal{A}_{\rho_{\lambda}}(H,K)$ be $p$-integral.
			We call
			$\mathbf{f}$
			\textbf{$p$-primitive} if $\overline{\mathbf{f}_p}\neq0$ (recall $\overline{\mathbf{f}_p}$ is reduction modulo $\mathfrak{P}$ of $\mathbf{f}_p$).
			
		\end{enumerate}
	\end{definition}

	\begin{remark}
		\rm
		In view of the transformation property of
		$\mathbf{f}_p$,
		the definition of
		$\M_\lambda(\mathcal{O})$
		and the assumption on
		$K_p$,
		$\mathbf{f}$ is $p$-integral if and only if
		for a (or any)
		set of representatives
		$S$ of
		$ H(\Q)
		\backslash
		H(\mathbb{A}_f)/K$,
		$\mathbf{f}_p(g)
		\in
		\M_\lambda(\mathcal{O})$
		for all $g\in S$.

	\end{remark}

	\subsubsection{Automorphic forms on $\widetilde{G}(\mathbb{A})$}

	Let $\mathbb{H}_{2n}$ be the Siegel upper half space of degree
	$2n$,
	\[
	\Hbb_{2n}
	=
	\left\{
	Z=X+iY\in\Symmr_{2n}(\C)|Y\text{ positive definite}
	\right\}.
	\]
	The group
	$G(\mathbb{R})$
	acts on
	$\mathbb{H}_{2n}$
	by
	fractional linear transformation
	\[
	\begin{pmatrix}
		A & B \\
		C & D
	\end{pmatrix}
	\cdot Z
	=
	(AZ+B)(CZ+D)^{-1}.
	\]
	The stabilizer of
	$i=i\cdot1_{2n}$
	in
	$G(\mathbb{R})$
	is
	$\mathbf{K}_\infty $,	
	the compact subgroup of
	$G(\mathbb{R})$
	consisting of $g$ such that
	$g^\mathrm{t}g=1$.
	Therefore we have an identification
	\[
	\mathbb{H}_{2n}\simeq G(\mathbb{R})/\mathbf{K}_\infty.
	\]
	The factor of automorphy is given by
	\[
	J 
	\colon
	G(\mathbb{R})
	\times
	\mathbb{H}_{2n}
	\rightarrow
	\mathrm{GL}_{2n}(\mathbb{C}),
	\quad
	(\begin{pmatrix}
		A & B \\
		C & D
	\end{pmatrix},Z)
	\mapsto
	CZ+D.
	\]    
	For a positive integer $N=\prod_\ell\ell^{e_{\ell}}$, we write $N'=\mathrm{l.c.m}(2,N)$ and define a congruence subgroup of $G(\Z)$:
	\begin{equation}\label{Gamma_0(2,N)}
		\Gamma_0(2,N)
		=
		\left\{
		\begin{pmatrix}
			A & B \\
			C & D
		\end{pmatrix}
		\in
		G(\widehat{\Z})\Big|
		\begin{array}{l}
			C_{\ell}
			\equiv
			0(\mathrm{mod}\,N'\Z_\ell)),
			B_{\ell}\equiv
			0
			(\mathrm{mod}\,2\Z_{\ell}),
			\\
			\mathrm{det}(D_{\ell})
			\equiv1
			(\mathrm{mod}\,4\Z_{\ell}),
			\forall \ell
		\end{array}		
		\right\}.
		\index{G@$\Gamma_0(2,N)$}
	\end{equation}
	The same argument as in \cite[Proposition 2.8, p.17] {Gelbart2006} or \cite[II.10, p.43]{MoeglinVignerasWaldspurger} shows that the subgroup $\Gamma_0(2,N)$ splits the two-fold cover $\widetilde{G}(\mathbb{A})\rightarrow G(\mathbb{A})$. For a compact open subgroup $K$ of $G(\mathbb{A}_f)$, by strong approximation, one has
	\begin{equation}\label{upper half space}
		(K\cap G(\Q))
		\backslash
		\mathbb{H}_{2n}
		=
		G(\Q)\backslash G(\mathbb{A})/(\mathbf{K}_\infty K)
		=
		G(\Q)\backslash \widetilde{G}(\mathbb{A})/\widetilde{\mathbf{K}_\infty K}.
	\end{equation}

	We fix a Cartan decomposition
	\[
	\mathrm{Lie}(G(\mathbb{R}))
	=\mathfrak{k}\oplus\mathfrak{p}
	\]
	Here $\mathfrak{k}$ is the Lie algebra of $\mathbf{K}_{\infty} $. Then the center of $\mathfrak{k}$ is of dimension one, generated by $Z=\frac{1}{2}\begin{pmatrix}
		0 & 1_{2n} \\ -1_{2n} & 0
	\end{pmatrix}$. Then the eigenvalues of $Z$ on $\mathfrak{p}^{\mathbb{C}}$ are $\pm i$. We write $\mathfrak{p}^\pm$ for the subspace of $\mathfrak{p}^{\mathbb{C}}$ consisting of $X$ such that $[Z,X]=\pm i X$. More explicitly,	write the Cayley matrix $\mathfrak{c}=\frac{1}{2}\begin{pmatrix}
		1_{2n} & i\cdot1_{2n} \\
		i\cdot1_{2n} & 1_{2n}
	\end{pmatrix}$, then
	\[
	\mathfrak{p}^+
	=
	\mathfrak{c}
	\begin{pmatrix}
		0 & \mathrm{Sym}_{2n}(\mathbb{C}) \\
		0 & 0
	\end{pmatrix}
	\mathfrak{c}^{-1},
	\quad
	\mathfrak{p}^-
	=
	\mathfrak{c}
	\begin{pmatrix}
		0 & 0 \\
		\mathrm{Sym}_{2n}(\mathbb{C}) & 0
	\end{pmatrix}
	\mathfrak{c}^{-1}.
	\]

	Write
	$U_G$ for the unipotent subgroup of
	$G$
	\[
	U_G=
	\left\{
	n(X)=
	\begin{pmatrix}
		1_{2n} & X \\
		0 & 1_{2n}
	\end{pmatrix}
	\Bigg|
	X=X^{\mathrm{t}}
	\right\}.\index{n@$n(X)$}
	\]
	We write $dn(X)$ for the normalized Haar measure on $[U_G]$ (the totally volume is $1$).

	We fix an additive character
	\[
	\mathbf{e}
	=
	\otimes_v\mathbf{e}_v
	\colon
	F\backslash\mathbb{A}
	\rightarrow
	\mathbb{C}^\times,
	\quad
	\mathbf{e}_v
	(x)
	=
	\begin{cases*}
		\exp(2i\pi x),
		&
		if
		$v\mid\infty$;
		\\
		\exp(-2i\pi\{x\}),
		&
		if $v$ a finite place.
	\end{cases*}\index{e@$\mathbf{e},\mathbf{e}_v$}
	\]
	Here
	$\{x\}$
	is the fractional part of
	$x\in\mathbb{Q}_v$.

	\begin{definition}\label{definition of Siegel modular forms}

		Fix a continuous character of finite order $\chi
		\colon
		(1+4\Z_{2})
		\times
		\prod_{\ell\neq2}
		\Z_{\ell}^\times
		\rightarrow
		\mathbb{C}^\times$,
		\begin{enumerate}
			\item 
			We denote the space of (genuine) metaplectic automorphic form of weight $\tau$, of level
			$\Gamma_0(2,N)$ and of character $\chi$ by
			\[
			\mathcal{A}_{\rho_{\tau}}
			(\widetilde{G},{\Gamma}_0(2,N),\chi)
			=
			\left\{
			\begin{array}{c}
				f\colon
				\widetilde{G}(\mathbb{A})
				\rightarrow
				\M_{\tau}(\C)
				\\
				\text{ smooth}
			\end{array}
			\Bigg|
			\begin{array}{l}
				f((\gamma gk_\infty k_f,\varepsilon))
				=
				\frac{\varepsilon\chi(\mathrm{det}(D))}{\mathrm{det}(J(k_{\infty},i))^{1/2}}
				\cdot
				\rho_{\tau}(J(k_\infty,i))^{-1}(f(g))
				\\
				\forall
				(\gamma,k_\infty,k_f)
				\in
				G(\Q)
				\times
				\mathbf{K}_\infty
				\times
				\Gamma_0(2,N)
				\text{ with }
				k_f=\begin{pmatrix}
					A & B \\
					C & D
				\end{pmatrix}
			\end{array}
			\right\}
			\]
			Here the square root $\mathrm{det}(J(k_{\infty},i))^{1/2}$ is defined using the principal argument of $\mathrm{det}(J(k_{\infty},i))$.
			We say that an automorphic form $f\in\mathcal{A}_{\rho_{\tau}}			(\widetilde{G},{\Gamma}_0(2,N),\chi)$ is holomorphic if moreover it is annihilated by $\mathfrak{k}^{\mathbb{C}}\oplus\mathfrak{p}^-$.

			\item 			
			For any
			$S\in\mathrm{Sym}_{2n}(\Q)$,
			the $S$-th Fourier coefficient
			of a metaplectic automorphic forms
			$f\in
			\mathcal{A}_{\rho_{\tau}}
			(\widetilde{G},\Gamma_0(2,N),\chi)$
			is given by
			\[
			\mathbf{a}_{f,S}
			\colon
			\widetilde{G}(\mathbb{A})
			\rightarrow
			\M_{\tau}(\C),
			\quad
			g
			\mapsto
			\int\limits_{[U_G]}
			f(n(X)g)
			\mathbf{e}(-\mathrm{Tr}(SX))
			dn(X).
			\index{a@$\mathbf{a}_{f,S}$}
			\]
			
		\end{enumerate}
	\end{definition}

	It follows from the definition that
	\begin{align}\label{Fourier coefficient, transformation property}
		\begin{split}
			\mathbf{a}_{f,S}
			(n(X)g)
			&
			=
			\mathbf{e}(-\mathrm{Tr}(SX))
			\mathbf{a}_{f,S}(g),
			\\
			\mathbf{a}_{f,S}
			(
			(\mathrm{diag}(A,A^{-\mathrm{t}}),1)
			g)
			&
			=
			\mathbf{a}_{f,A^\mathrm{t}SA}(g),
			\quad
			\forall
			A\in
			\mathrm{GL}_{2n}(\Q).
		\end{split}
	\end{align}
	Moreover,
	we have the Fourier expansion for $f$
	\begin{equation*}
		f(g)
		=
		\sum_{S\in\mathrm{Sym}_{2n}(\Q)}
		\mathbf{a}_{f,S}(g),\quad
		\forall
		g\in\widetilde{G}(\mathbb{A}).
	\end{equation*}

	We next relate metaplectic automorphic forms to
	classical Siegel modular forms:
	\begin{definition}
		For any
		$Z=X+iY\in\mathbb{H}_{2n}$,
		we choose
		$g\in G(\R)$
		such that
		$g(i)=Z$.
		Then
		for a metaplectic automorphic form
		$f\in
		\mathcal{A}_{\rho_{\tau}}(\widetilde{G},{\Gamma}_0(2,N),\chi)$,
		its associated classical
		Siegel modular form
		$f^\ast
		\colon
		\mathbb{H}_{2n}
		\rightarrow
		\M_{\tau}(\C)
		$
		is given by
		\[
		f^\ast
		(Z)
		=
		\mathrm{det}(J(g,i))^{1/2}
		\cdot
		\rho_{\tau}(J(g,i))
		f((g,1)).\index{f@$f^\ast$}
		\]
	\end{definition}
	It follows from the definition that
	$f^\ast$
	is independent of the choice $g$. Moreover due to (\ref{upper half space}) and the transformation property of $f$ under the element $(1,-1)\in\widetilde{G}(\mathbb{A})$, $f$ is uniquely determined by $f^\ast$. For any
	$\gamma=\begin{pmatrix}
		1_{2n} & b \\ 0 & 1_{2n}
	\end{pmatrix}\in\Gamma_0(2,N)$,
	we have
	\[
	f^\ast(Z+b)
	=
	f^\ast(\gamma(Z))=
	f^\ast(Z).
	\]
	We assume in the following that $f$ is holomorphic.
	Then the above identity gives us the Fourier expansion
	for the classical Siegel modular form
	\begin{equation}\label{Fourier expansion}
		f^\ast(Z)
		=
		\sum_{S\in\mathrm{Sym}_{2n}^\circ(\Z)}a_f(S)q^S
		\index{a@$a_f(S)$}
	\end{equation}
	with
	\begin{align*}
		q^S
		&
		=
		\mathbf{e}_\infty(\mathrm{Tr}(SZ )),
		\\
		\mathrm{Sym}_{2n}^\circ(\Z)
		&
		=
		\{
		S\in
		\mathrm{Sym}_{2n}(\Q)|
		2S\in
		\mathrm{Sym}_{2n}(\Z)
		\text{ and }
		S_{1,1},S_{2,2},\cdots,S_{2n,2n}
		\in\Z
		\}.\index{S@$\mathrm{Sym}_{2n}^\circ(\Z)$}
	\end{align*}
	By definition, we have
	\[
	a_f(S)
	=
	\int_{\mathrm{Sym}_{2n}(\R)/\mathrm{Sym}_{2n}(\Z)}
	f^\ast(X+iY)\exp(-2i\pi\mathrm{Tr}(S(X+iY)))dX.
	\]
	Thus it follows from the definition of $\mathbf{a}_{f,S}$ that for any matrix $A\in\mathrm{GL}_{2n}(\R)$, we have
	\begin{equation}\label{a_{f,S} and a_f(S)}
		a_f(S)
		=
		\mathrm{det}(A)^{-1/4}\exp(2\pi\mathrm{Tr}(A^{\mathrm{t}}SA))
		\rho_{\tau}(A^{-\mathrm{t}})
		\mathbf{a}_{f,S}(\mathrm{diag}(A,A^{-\mathrm{t}})).
	\end{equation}
	For this, it is enough to put $Y=AA^\mathrm{t}$ in the expression for $a_f(S)$ and use the relation between $f$ and $f^\ast$.

	\begin{definition}\label{p-integral Siegel modular forms}
		Let
		$f\in
		\mathcal{A}_{\rho_{\tau}}(\widetilde{G},\Gamma_0(2,N),\chi)$
		be as above.
		We call
		$f$
		\textbf{$p$-integral}
		if
		$a_f(S)\in
		\M_{\tau}
		(\mathcal{O})$
		for all
		$S\in\mathrm{Sym}_{2n}^\circ(\Z)$.
		We call
		$f$
		\textbf{$p$-primitive}
		if
		furthermore
		$a_f(S_0)\not\equiv
		0(\mathrm{mod}\,\mathfrak{P})$
		for some
		$S_0
		\in
		\mathrm{Sym}_{2n}^\circ(\Z)$.
	\end{definition}

	\begin{remark}
		\rm 
		Our notion of $p$-integral metaplectic automorphic forms
		coincides with the one defined using automorphic sheaves on
		Siegel modular varieties over
		$\Z[1/(2N)]$,
		since both of which are characterized by the $p$-integrality of Fourier expansions.
	\end{remark}

	\subsection{The pairings}

	Recall we have fixed $\tau=(\tau_1,\cdots,\tau_{2n})$. We write
	\begin{equation}\label{tau^circ}
		\tau^\circ
		=
		\tau+\underline{n}
		=
		(\tau_1+n,\cdots,\tau_{2n}+n).\index{t@$\tau^\circ$}
	\end{equation}

	Under condition \ref{condition on p}, we have an embedding $R'\hookrightarrow\mathcal{O}$ of rings.

	For any positive integer $r$, we extend the symmetric pairing
	$\langle-,-\rangle_U$
	on
	$U$ to
	a non-degenerate
	$H$-equivariant
	pairing	
	on
	$U^{\otimes r}$
	in a natural way (still denoted by $\langle-,-\rangle_U$):
	\[
	\left\langle
	v_1\otimes v_2\otimes\cdots\otimes v_r,
	v_1'\otimes v_2'\otimes\cdots\otimes v_r'
	\right\rangle_U
	:=
	\prod_{i=1}^r
	\langle v_i,v_i'\rangle_U.
	\]	
	Similarly this extends to pairings $\langle-,-\rangle_U$ on $\widetilde{\M}_\lambda$ and $\M_\lambda$.
	Under
	condition
	\ref{condition on p},
	we have a \textit{perfect}
	$ H(\mathcal{O})$-equivariant pairing
	\[
	\langle
	-,-
	\rangle_U
	\colon
	\M_\lambda(\mathcal{O})
	\otimes_{\mathcal{O}}
	\M_\lambda(\mathcal{O})
	\rightarrow
	\mathcal{O}.
	\]
	This induces
	$ H(\mathcal{O})$-equivariant 
	pairings
	\begin{align}\label{pairing W,U-complex case}
		\begin{split}
			\langle
			-,-
			\rangle_{W,U}
			\colon
			&
			\left(
			\M_\lambda(\mathcal{O})
			\otimes_{\mathcal{O}}
			\M_{\tau}(\mathcal{O})
			\right)
			\otimes
			\M_\lambda(\mathcal{O})
			\rightarrow
			\M_{\tau}
			(\mathcal{O}),
			\\
			\langle
			-,-
			\rangle_{W,U}
			\colon
			&
			\left(
			\M_\lambda(\mathcal{O})
			\otimes_{\mathcal{O}}
			\M_{\tau^\circ}(\mathcal{O})
			\right)
			\otimes
			\M_\lambda(\mathcal{O})
			\rightarrow
			\M_{\tau^\circ}
			(\mathcal{O}).
		\end{split}
	\end{align}

	Recall we have an isotropic decomposition $V=V^+\oplus V^-$. The embedding $\GLmr_{2n}\hookrightarrow G$ sending $g$ to $\mathrm{diag}(g,g^{-\mathrm{t}})$ induces actions of $\GLmr_{2n}$ on $V^+$ and $V^-$.
	The symplectic form
	on
	$V$
	induces a perfect $\GLmr_{2n}(\Z)$-equivariant pairing
	$\langle-,-\rangle_V
	\colon
	V^+\times V^-
	\rightarrow
	\Z$,
	and therefore a perfect pairing on the tensor products
	$\langle-,-\rangle_V
	\colon
	V^+(\mathcal{O})^{\otimes k}
	\otimes
	V^-(\mathcal{O})^{\otimes k}
	\rightarrow
	\mathcal{O}$
	for $k\ge0$.
	So
	we get the following
	$\GLmr_{2n}(\mathcal{O})$-equivariant pairings
	\begin{align*}
		\langle-,-\rangle_V
		&
		\colon
		\M_{\tau}(\mathcal{O})
		\otimes
		\M_{\tau^\vee}(\mathcal{O})
		\rightarrow
		\mathcal{O},
		\\
		\langle-,-\rangle_{W,V}
		&
		\colon
		\left(
		\M_\lambda(\mathcal{O})
		\otimes
		\M_{\tau}(\mathcal{O})
		\right)
		\otimes
		\M_{\tau^\vee}(\mathcal{O})
		\rightarrow
		\M_\lambda(\mathcal{O}).
	\end{align*}
	Under condition \ref{condition on p},
	the first pairing is
	perfect.

	\begin{lemma}
		For any
		$w\in
		\M_\lambda(\mathcal{O})
		\otimes
		\M_{\tau}(\mathcal{O})$,
		$u\in\M_\lambda(\mathcal{O})$,
		$v\in\M_{\tau^\vee}(\mathcal{O})$,
		we have
		\[
		\langle
		\langle
		(g,h)w,hv
		\rangle_{W,V},
		gu
		\rangle_U
		=
		\langle
		\langle
		(g,h)w,gu
		\rangle_{W,U},
		hv
		\rangle_V,
		\quad
		\forall
		g\in H(\mathcal{O}),
		\,
		h\in\GLmr_{2n}(\mathcal{O})
		\]		
	\end{lemma}
	\begin{proof}
		We have the following commutative diagram
		\[
		\begin{tikzcd}
			(\M_{\lambda}(\mathcal{O})
			\otimes
			\M_{\tau}(\mathcal{O}))
			\otimes
			(\M_{\lambda}(\mathcal{O})
			\otimes
			\M_{\tau^\vee}(\mathcal{O}))
			\arrow[r,"{\langle-,-\rangle_V}"]
			\arrow[d,"{\langle-,-\rangle_U}"]
			\arrow[rd,"{\langle-,-\rangle}"]
			&
			\M_{\lambda}(\mathcal{O})
			\otimes
			\M_{\lambda}(\mathcal{O})
			\arrow[d,"{\langle-,-\rangle_U}"]
			\\
			\M_{\tau}(\mathcal{O})
			\otimes
			\M_{\tau^\vee}(\mathcal{O})
			\arrow[r,"{\langle-,-\rangle_V}"]
			&
			\mathcal{O}
		\end{tikzcd}
		\]
		We denote the diagonal map/pairing by
		$\langle-,-\rangle$.
		Note that the horizontal pairings are
		$\GLmr_{2n}(\mathcal{O})$-equivariant while the vertical pairings
		are	$H(\mathcal{O})$-equivariant.
		Moreover we have
		\[
		\langle\langle w,v\rangle_{ W,V},u\rangle_U
		=
		\langle w,v\otimes u\rangle
		=
		\langle\langle w,u\rangle_{W,U},v\rangle_{V}.
		\]
		This concludes the proof of the lemma.
	\end{proof}
	
	\begin{remark}
		In this section, we can replace $\mathcal{O}$ by $R'$ everywhere and the results still hold. The main reason of working with $\mathcal{O}$ comes from the consideration of $p$-integrality and $p$-primitivity of theta lifts.
	\end{remark}

	\section{Theta lifts}\label{theta lift}
	We recall the notions of Weil representations and theta lifts.
	For later computation,
	we fix Haar measures on
	the spaces
	$V(\Q_v)$,
	the groups
	$G(\Q_v)$,
	etc.
	as follows
	\begin{enumerate}[label=(\alph*)]
		\item 
		The measure on
		$U(\Q_{\ell})$
		is such that
		the volume of
		$U(\Z_{\ell})$
		is equal to
		$|\delta_1\delta_2\cdots\delta_{2n+1}|_\ell$.
		The measure on
		$U(\mathbb{R})$
		is
		$\delta_1\cdots\delta_{2n+1}$
		times the standard Lebesgue measure with respect to the
		$\mathbb{R}$-basis
		$(E_1,\cdots,E_{2n+1})$.
		The measure on
		$ H(\Q_{\ell})$
		is such that
		$ H(\Z_{\ell})$
		has volume $1$.
		The measure on
		$ H(\mathbb{R})$
		is chosen such that it has total volume equal to $1$.
		The measure on
		$T(\mathbb{R})$
		is such that the total volume is $1$
		and the measure on
		$T(\Q_{\ell})$
		is such that
		$T(\Z_{\ell})$
		has volume $1$.

		\item 
		For a finite place
		$v=\ell$,
		fix the measure on
		$V^\pm(\Q_{\ell})$
		such that
		the volume of
		$V^\pm(\Z_{\ell})$
		is equal to $1$.
		we fix the measure
		on $V^\pm(\R)$
		to be the standard Lebesgue measure
		with respect to the $\mathbb{R}$-basis
		$(e_1^\pm,\cdots,e_{2n}^\pm)$.
		The measure on
		$G(\Q_{\ell})$
		is such that
		the volume of
		$G(\Z_{\ell})$
		is equal to $1$.
		The measure on
		$G(\mathbb{R})$
		is induced from
		the measures on
		the Siegel upper half space
		$\mathbb{H}_{2n}$
		and on
		$\mathbf{K}_\infty $,
		the former is given by
		\[
		\mathrm{det}(Y)^{-2n-1}
		\prod_{1\le i\le j\le 2n}
		dX_{i,j}dY_{i,j}
		\]
		where
		$Z=X+iY\in\mathbb{H}_{2n}$,
		and the latter is chosen such that
		$\mathbf{K}_\infty $
		has total volume $1$.

	\end{enumerate}

	For each place $v$ of $\Q$
	and
	$X=V,V^\pm,U,W,W^\pm$,
	we also write
	\[
	X_v=X(\Q_v).
	\]
	Recall we have identified
	$W^+$ with
	$\mathrm{M}_{2n+1,2n}(\Z)$
	using the basis vectors
	$(E_1,\cdots,E_{2n+1})$ and
	$(e_1^+,\cdots,e_{2n}^+)$.
	Thus $H(\Z)\times \mathrm{GL}_{2n}(\Z)$ acts on
	$W^+$ as follows
	\[
	(h,g)v=hvg^{-1}.
	\]

	Let $V$ be a $\C$-vector space of finite dimension.	Write $\mathcal{S}(W^+_v)$\index{S@$\mathcal{S}(W^+_v)$}
	for the space of $\mathbb{C}$-valued Bruhat-Schwartz functions on
	$W^+_v$ and $\mathcal{S}(W_v^+)\otimes_\C V$ for the space of  $V$-valued Bruhat-Schwartz functions on $W^+_v$.
	Then the local Weil representation $\omega_{W_v^+}$
	of $ H(\Q_v)\times\widetilde{G}(\Q_v)$ on $\mathcal{S}(W^+_v)\otimes V$
	is given by the Schrödinger model as follows: let
	$\chi_{U_v}\colon\Q_v^\times\rightarrow\{\pm1\}$ denote the quadratic character attached to the quadratic space $U_v$, given by
	\[
	\chi_{U_v}(x)=
	\left(
	x,(-1)^{n(2n+1)}\mathrm{det}(U_v)
	\right)_v\index{k@$\chi_{U_v}$}
	\]
	where $(-,-)_v$ is the Hilbert symbol on $\Q_v^\times$. Write $q_{U,v}$
	for the quadratic form associated to the bilinear form $\langle-,-\rangle_U$
	on $U_v$ and write $\gamma(q_{U,v})$ for the Weil index of the character of second degree $\Q_v\rightarrow\mathbb{C}^\times$ sending $x$ to
	$\mathbf{e}_v(q_{U,v}(x))$ (\cite[Theorem A.1]{Rao1993}). This is an $8$-th root of unity. For $a\in \Q_v^\times$, write $aq_{U,v}$ for the product of $a$ with $q_{U,v}$ and set $\gamma(a,q_{U,v})=\gamma(aq_{U,v})/\gamma(q_{U,v})$
	(\cite[A.3]{Rao1993}). Then we have (\cite[Lemma 4.1, p.17]{Kudla1996}):
	\[
	\gamma(ab,q_{U,v})=\gamma(a,q_{U,v})\gamma(b,q_{U,v})
	(a,b)_v,
	\quad
	\forall
	a,b\in \Q_v^\times.
	\]
	In particular we have the following
	\begin{equation}\label{gamma is a character}
		\gamma(ab,q_{U,v})
		=
		\gamma(a,q_{U,v})\gamma(b,q_{U,v})
		\text{ if }
		\begin{cases*}
			\text{either }
			v=\ell\nmid 2\infty
			\text{ and }
			a,b\in\Z_\ell^\times;
			\\
			\text{or }
			v=2
			\text{ and }
			a,b\in1+4\Z_2;
			\\
			\text{or }
			v=\infty
			\text{ and }
			a,b
			\text{ totally positive}.
		\end{cases*}
	\end{equation}
	The computations of these values $\gamma(q_{U,v})$ and $\gamma(a,q_{U,v})$ can be found in \cite[A.4]{Rao1993}. In particular, we have
	\begin{equation}\label{value for gamma}
		\gamma(q_{U,v})
		=
		1,
		\quad
		\forall
		v\nmid2\delta_1\cdots\delta_{2n+1}.
	\end{equation}
	We write
	\[
	\gamma(q)=\prod_v\gamma(q_{U,v}).
	\]
	With these preparations we can now give the following well-known formulas for the Weil representation of $H(\Q_v)\times \widetilde{G}(\Q_v)$ on the space
	$\mathcal{S}(W_v^+)\otimes V$ (\cite[pp.38-39]{Kudla1996}): for any $f\in\mathcal{S}(W_v^+)\otimes V$, one has (recall $\varepsilon=\pm1$)
	\begin{align}\label{Weil representation}
		\begin{split}
			\omega_{W_v^+}
			\left(
			\begin{pmatrix}
				A & 0 \\
				0 & A^\mathrm{-t}
			\end{pmatrix},\varepsilon
			\right)
			f(x)
			&
			=\varepsilon
			\gamma(\mathrm{det}(A),q_{U,v})^{-1}
			\chi_{U_v}(\mathrm{det}(A))
			|\mathrm{det}(A)|_v^{(2n+1)/2}
			f(xA),
			\\
			\omega_{W_v^+}
			\left(
			\begin{pmatrix}
				1_{2n} & B \\
				0 & 1_{2n}
			\end{pmatrix},\varepsilon
			\right)    
			f(x)
			&
			=\varepsilon
			\mathbf{e}_v
			(
			\frac{1}{2}
			\mathrm{Tr}(S_xB)
			)
			f(x),
			\\
			\omega_{W_v^+}
			\left(
			\begin{pmatrix}
				0 & 1_{2n} \\
				-1_{2n} & 0
			\end{pmatrix},\varepsilon
			\right)
			f(x)
			&
			=\varepsilon
			\gamma(q_{U,v})^{-2n}
			\cdot
			\int\limits_{W_v^+}
			f(y)
			\mathbf{e}_v
			(\langle
			Q_V
			y,x\rangle_W)
			dy,
			\\
			\omega_{W_v^+}
			(h)f(x)
			&
			=
			f(h^{-1}x),
			\quad
			h\in
			H(\Q_v).
		\end{split}        
	\end{align}
	Here we set
	\[
	S_x=x^\mathrm{t}Q_Ux\index{S@$S_x$}.
	\]

	We denote the restricted tensor product of all Bruhat-Schwartz functions on $W_v^+$ by
	\[
	\mathcal{S}(W^+(\mathbb{A}))\otimes V
	=
	\bigotimes_v\,'
	(\mathcal{S}(W^+_v)\otimes V).
	\index{S@$\mathcal{S}(W^+(\mathbb{A}))$}
	\]
	The global Weil representation $\omega_{W^+}$ of $ H(\mathbb{A})\times
	\widetilde{G}(\mathbb{A})$ on $\mathcal{S}(W^+(\mathbb{A}))\otimes V$	
	is the tensor product of local ones
	\[
	\omega_{W^+}
	=
	\bigotimes_v\omega_{W_v^+}.\index{o@$\omega_{W^+},\omega_{W_v^+}$}
	\]
	For any
	$\phi
	\in
	\mathcal{S}(W^+(\mathbb{A}))\otimes V$,
	we define the associated theta series
	$\Theta_\phi$
	to be
	a function on
	$ H(\mathbb{A})
	\times
	\widetilde{G}(\mathbb{A})$
	given by
	\[
	\Theta_\phi
	(h,g)
	:=
	\sum_{w\in W^+(\Q)}
	\omega_{W^+}(h,g)\phi(w),
	\quad
	(h,g)
	\in
	H(\mathbb{A})\times\widetilde{G}(\mathbb{A}).
	\]
	\begin{definition}\label{definition of theta lift}
		We take $V=\mathbb{M}_\lambda(\C)
		\otimes_\C\mathbb{M}_{\tau^\circ}(\C)$. For an automorphic form		$\mathbf{f}\in\mathcal{A}_{\rho_{\lambda}}( H,K)$, its \textbf{theta lift} to
		$\widetilde{G}(\mathbb{A})$ by $\phi\in	\mathcal{S}(W^+(\mathbb{A}))\otimes V$ is the $\mathbb{M}_{\tau^\circ}(\C)$-valued metaplectic automorphic form $\Theta_{\phi,\mathbf{f}}$ on $\widetilde{G}(\mathbb{A})$ given by
		\[
		\Theta_{\phi,\mathbf{f}}(g)
		:=
		\int\limits_{[ H]}
		\langle
		\mathbf{f}(h),\Theta_\phi(h,g)
		\rangle_{W,U}
		dh
		\in\mathbb{M}_{\tau^\circ}(\C),
		\quad
		g\in
		\widetilde{G}(\mathbb{A}).
		\]
	\end{definition}
	We will be more precise about the weight, level and character of $\Theta_{\phi,\mathbf{f}}$ later on.

	\section{Bruhat-Schwartz function}\label{Bruhat-Schwartz function}
	In this section we pick one particular Bruhat-Schwartz function
	and then study some basic properties of
	theta lifts by the theta series associated to
	this function.

	In the rest of this article, we fix
	\begin{enumerate}
		\item 
		a finite set $\mathbb{S}$\index{S@$\mathbb{S}$} of finite places of $\Q$, including those dividing $2\delta_1\cdots\delta_{2n+1}$ and
		
		\item 
		a positive integer $r_\ell$\index{r@$r_\ell$} for each $\ell\in\mathbb{S}$ such that
		\[
		\mathrm{ord}_\ell(\delta_{2n+1})+r_\ell
		>\max_{i=1,\cdots,2n}\mathrm{ord}_\ell(\delta_i).
		\]
		If $\ell=2$, we require $r_\ell\ge2$.
	\end{enumerate} 
	Then for each finite place $\ell$, we put
	\[
	K_\ell
	=
	\begin{cases*}
		\left\{
		\begin{pmatrix}
			A & B \\
			C & D
		\end{pmatrix}\in H(\Z_\ell)
		\mid
		A\in\mathrm{GL}_{2n}(\Z_\ell),
		\,
		C\in\ell^{r_\ell}\mathrm{M}_{1,2n}(\Z_\ell)
		\right\},
		&
		$\ell\in\mathbb{S}$;
		\\
		H(\Z_\ell),
		&
		$\ell\notin\mathbb{S}$.
	\end{cases*}
	\]
	We put $K^{(\mathbb{S})}=\prod_{\ell}K_\ell$\index{K@$K^{(\mathbb{S})}$}, a compact open subgroup of $H(\mathbb{A}_f)$ which will serve as the level subgroup of automorphic forms on $H(\mathbb{A})$.

	\subsection{A choice of Bruhat-Schwartz function}
	We choose one particular
	$\phi=\otimes_v'\phi_v
	\in
	\mathcal{S}(W^+(\mathbb{A}))$
	with which we study the non-vanishing modulo $p$ 
	of theta lifts.

	For each $\ell\in\mathbb{S}$, we consider the following compact open subset of $W^+(\Q_\ell)\simeq\mathrm{M}_{n,2n}(\Q_\ell)$:
	\[
	\mathrm{Supp}_\ell(r_\ell)
	=
	\begin{cases*}
		\begin{pmatrix}
			\mathrm{GL}_{2n}(\Z_\ell)
			\\
			\ell^{r_\ell}\mathrm{M}_{1,2n}(\Z_\ell)
		\end{pmatrix}	
		\subset
		\mathrm{M}_{n,2n}(\Q_\ell),
		&
		$\ell\in\mathbb{S}$;
		\\
		\mathrm{M}_{n,2n}(\Z_\ell),
		&
		$\ell\notin\mathbb{S}$.
	\end{cases*}\index{S@$\mathrm{Supp}_\ell(r_\ell)$}
	\]
	Then for each $\ell\in\mathbb{S}$, we define the Bruhat-Schwartz function at $\ell$ to be the characteristic function
	\[
	\phi_\ell
	=
	\mathbb{I}_{\mathrm{Supp}_\ell(r_\ell)}.
	\]

	For
	$v=\infty$,
	the construction is more involved. We define $\phi_{\infty}\in\mathcal{S}(W^+_\infty)$ by
	\[
	\phi_{\infty}(w)
	=
	\mathbf{e}_\infty
	\left(
	i\langle w,Q_Vw\rangle_{W}
	\right)
	=
	\mathbf{e}_\infty
	(i\mathrm{Tr}(w^\mathrm{t}Q_Uw)),
	\quad
	\forall
	w\in W^+(\R).
	\]

	We construct a highest weight vector in $\mathfrak{H}_{\lambda}$ as follows: 
	for a matrix
	$A$ of size
	$s\times t$
	and integers $0\le s'\le s,0\le t'\le t$,
	we write
	$A_{(s',t')}^{\mathrm{ul}}$
	for the upper-left
	block matrix of $A$ of size
	$s'\times t'$.
	Following
	\cite[(II.6.10)]{KashiwaraVergne1978},
	for an element
	$z =\begin{pmatrix}
		x \\
		y \\
		t
	\end{pmatrix}\in
	W^+(\C)$
	with
	$x=(x_{i,j}),y=(y_{i,j})\in
	\mathrm{M}_{n,2n}(\mathbb{C})$
	and
	$t=(t_{1,i})\in\mathrm{M}_{1,2n}(\mathbb{C})$,
	we define the following quantities associated to $z $:
	\begin{align}\label{vectors to define Delta polynomial}
		\begin{split}
			z_{(r)} 
			&
			=
			\mathrm{det}
			(x_{(r,r)}^{\mathrm{ul}}),
			\\
			\widetilde{z}_{(r)}
			&
			=
			\det
			\begin{pmatrix}
				x_{1,1} & \cdots & x_{1,2n+1-r}
				\\
				\vdots & \ddots & \vdots
				\\
				x_{n,1} & \cdots & x_{n,2n+1-r}
				\\
				y_{r+1,1} & \cdots & y_{j+1,2n+1-r}
				\\
				\vdots & \ddots & \vdots
				\\
				y_{n,1} & \cdots & y_{n,2n+1-r}
				\\
				t_{1,1} & \cdots & t_{1,2n+1-r}
			\end{pmatrix}
			\text{ with }
			r\leq n,
			\\
			z_{\lambda} 
			&
			=
			\begin{cases*}
				z_{(1)}^{\lambda_1 -\lambda_2 }
				z_{(2)}^{\lambda_2 -\lambda_3 }
				\cdots
				z_{(n)}^{\lambda_{n} },
				&
				$\varepsilon=(-1)^{\sum_i\lambda_i}$,
				\\
				z_{(1)}^{\lambda_1 -\lambda_2 }
				z_{(2)}^{\lambda_2 -\lambda_3 }
				\cdots
				z_{(n)}^{\lambda_{n} }
				\widetilde{z}_{(r)}
				&
				$\varepsilon=-(-1)^{\sum_i\lambda_i}$ and $1\leq r\leq n$.
			\end{cases*}
		\end{split}\index{z@$z_\lambda$}
	\end{align}
	Here $r$ is the maximal integer such that $\lambda_r\neq0$.
	Then by Proposition II.(6.11) of \textit{loc.cit}, $(\mathcal{T}z)_{\lambda}$, viewed as a function on $z$, lies in $\mathfrak{H}_\lambda$, and is a highest weight vector of weight $(\lambda,\tau)$ under the natural action of
	$ H(\C)\times\mathrm{GL}_{2n}(\C)$. Here $\mathcal{T}$ is given in (\ref{T}) and $\mathcal{T}z$ is the product of these two matrices $\mathcal{T}$ and $z$.

	\begin{remark}\rm
		Note that in \textit{loc.cit}, $H(\C)$ is defined using the matrix $J$ as in p.21 and p.25 of \textit{loc.cit} while our $H(\C)$ is defined using the diagonal matrix $Q_U$, that is why we use $(\mathcal{T}z)_\lambda$ instead of $z_\lambda$ in the above.
	\end{remark}

	We choose a scalar product on $\M_{\tau}$ such that $\rho_{\tau}(g^{\mathrm{t}})=\rho_{\tau}(g)^{\mathrm{t}}$ for $g\in\mathrm{GL}_{2n}(\C)$.
	We fix $\mathcal{O}$-basis for $\M_\lambda(\mathcal{O})$ and $\M_\tau(\mathcal{O})$:
	\begin{align*}
		\M_\lambda(\mathcal{O})
		&
		=\mathcal{O}(\mathbf{e}_1,\cdots,\mathbf{e}_{d_{\lambda}}),
		\\
		\M_\tau(\mathcal{O})
		&
		=\mathcal{O}(\mathbf{E}_1,\cdots,\mathbf{E}_{d_{\tau}}).
	\end{align*}
	We will take $\mathbf{e}_1=v_{\lambda}$ and $\mathbf{E}_1$ to be respectively the highest weight vectors. We have the dual basis
	\[
	\M_{\lambda}(\mathcal{O})\simeq
	\M_{\lambda^\vee}(\mathcal{O})
	=
	\mathcal{O}(\mathbf{e}_1^\vee,\cdots,\mathbf{e}_{d_\lambda}^\vee).
	\]
	Recall from (\ref{decomposition of pluri-harmonic}), we have an isomorphism
	\[
	\M_{\lambda}(\C)\otimes\M_{\tau}(\C)
	\simeq
	\mathfrak{H}_{\lambda},
	\quad
	\mathbf{e}_i\otimes\mathbf{E}_j
	\mapsto
	z_{\lambda;i,j}.
	\]
	Here each $z_{\lambda;i,j}$ is a pluri-harmonic polynomial on $z$ and $z_{\lambda;1,1}=(\mathcal{T}z)_\lambda$. We write
	\[
	\Delta_{\lambda}(z)
	=
	\sum_{i,j}
	z_{\lambda;i,j}
	\mathbf{e}_i^\vee\otimes\mathbf{E}_j,\index{D@$\Delta_{\lambda}$}
	\]
	which clearly satisfies the following formula
	\begin{equation}\label{joint action of O(U) and GL on harmonic polynomial}
		\Delta_{\lambda}
		(\mathcal{T}hzg)
		=
		\left(
		\rho_{\lambda}(h)
		\otimes
		\rho_{\tau}(g^{\mathrm{t}})
		\right)
		\Delta_{\lambda}(\mathcal{T}z),
		\quad
		\forall
		h\in H(\C),
		\,
		g\in
		\mathrm{GL}_{2n}(\C).
	\end{equation}

	We define our archimedean
	(vector-valued)
	Bruhat-Schwartz function
	$\phi_{\lambda,\infty}$ as follows
	(\textit{cf.}
	\cite[p.27]{KashiwaraVergne1978})
	\[
	\phi_{\lambda,\infty}
	\colon
	W^+(\R)
	\rightarrow\M_{\lambda}(\C)
	\otimes
	\M_{\tau^\circ}(\C),
	\quad
	z
	\mapsto
	\phi_\infty(z)
	\cdot
	\Delta_{\lambda}(\mathcal{T}z),
	\]
	which satisfies the following
	\begin{equation}\label{action of O(U) and GL_n}
		\phi_{\lambda,\infty}(hzg)
		=
		\left(       
		\rho_{\lambda}
		(h)
		\otimes
		\rho_{\tau^\circ}
		(g^{\mathrm{t}})
		\right)
		\phi_{\lambda,\infty}(z),
		\quad
		\forall
		h\in H(\R),
		\,
		g\in
		\mathrm{GL}_{2n}(\R),
		\,
		z\in
		W^+(\R).
	\end{equation}

	\begin{lemma}\label{action of O(U) and compact of Sp(V)}
		For any
		$h\in H(\R)$,
		$k\in
		\mathbf{K}_\infty$
		and
		$v\in W^+_\infty$, we have
		\[
		\omega_{W_{\infty}^+}(h,(k,\varepsilon))
		\phi_{\lambda,\infty}(v)
		=
		\varepsilon\cdot\mathrm{det}(J(k,i))^{1/2}\cdot
		\left(
		\rho_{\lambda}(h^{-1})
		\otimes
		\rho_{\tau^\circ}(k^{-1})
		\right)
		\phi_{\lambda,\infty}(v).
		\]
	\end{lemma}	
	\begin{proof}
		This follows from our definition of Weil representation for $H(\R)\times\widetilde{G}(\R)$ and the fact $k^{-1}=k^{\mathrm{t}}$. See also \cite[p.16, Proposition 4.4]{KashiwaraVergne1978}.
	\end{proof}

	We put these local factors
	$\phi_v$
	together and set
	\begin{equation*}
		\phi_{\lambda}
		=
		\phi_{\lambda,\infty}
		\bigotimes
		\Big(
		\bigotimes_{v\neq\infty}
		\phi_v
		\Big)
		\colon
		W^+(\mathbb{A})
		\rightarrow\M_{\lambda}(\C)
		\otimes
		\M_{\tau^\circ}(\C).\index{f@$\phi_\lambda$}
	\end{equation*}

	\subsection{Fourier coefficients of theta lifts}
	For
	$\mathbf{f}\in
	\mathcal{A}_{\rho_\lambda}( H,K)$
	and
	$\phi_{\lambda}$
	as above,
	the theta lift of
	$\mathbf{f}$ by
	$\phi_\lambda$
	is by definition given by
	\[
	\Theta_{\phi_\lambda,\mathbf{f}}(g)
	=
	\int\limits_{[ H]}
	\langle\Theta_{\phi_{\lambda}}(h,g),\mathbf{f}(h)\rangle_{W,U}
	dh
	\in\M_{\tau^\circ},
	\quad
	g\in \widetilde{G}(\mathbb{A}).
	\]

	Moreover we have fixed a finite set $\mathbb{S}$ of places of $\Q$ and an integer $r_\ell$ for each $\ell\in\mathbb{S}$. We write
	\[
	N_{\mathbb{S}}=\prod_{\ell\in\mathbb{S}}\ell^{2r_{\ell}}.\index{N@$N_{\mathbb{S}}$}
	\]
	We then set
	\[
	\chi_{U}=\otimes_v'\chi_{U_v}
	\colon
	\Q^\times\backslash\mathbb{A}^\times\to\mathbb{C}^\times.
	\]
	We define a product \textit{map}
	$\chi_{U}^\circ
	\colon
	\mathbb{A}^\times
	\rightarrow
	\mathbb{C}^\times$
	sending
	$x$ to
	$\gamma(x)\chi_{U}(x)$.
	In general this is not a group homomorphism.
	However, by (\ref{gamma is a character}),
	its restriction to the following compact open subgroup is
	indeed a group homomorphism
	\[
	\chi_{U}^\circ
	\colon
	(1+4\Z_2)
	\times
	\prod_{\ell\ne2}
	\Z_{\ell}^\times
	\rightarrow
	\mathbb{C}^\times,
	\quad
	x
	\mapsto
	\gamma(x)
	\chi_{U}(x).\index{k@$\chi_{U},\chi_{U}^\circ$}
	\]
	\begin{lemma}
		For any
		$g=\begin{pmatrix}
			A & B \\
			C & D
		\end{pmatrix}
		\in
		\Gamma_0(2,N_{\mathbb{S}})$
		with
		$A,B,C,D$
		block matrices of size
		$(2n)\times(2n)$,
		we have
		\[
		\omega_{W^+}((g,\varepsilon))
		\phi_{\lambda}
		=
		\varepsilon
		\chi_{U}^\circ
		(\mathrm{det}\,D)
		\phi_{\lambda}.
		\]
	\end{lemma}
	\begin{proof}
		We will drop the component $\varepsilon$ since $\Gamma_0(2,N_{\mathbb{S}})$ splits the cover $\widetilde{G}(\mathbb{A})\rightarrow G(\mathbb{A})$.
		For a finite place
		$\ell\nmid2\delta_U$,
		we know that
		$g_\ell$
		can be expressed as a product of matrices in $G(\Z_\ell)$
		of the forms
		$\begin{pmatrix}
			a & 0 \\
			0 & a^{-\mathrm{t}}
		\end{pmatrix}$,
		$\begin{pmatrix}
			0 & 1_{2n} \\
			-1_{2n} & 0
		\end{pmatrix}$
		and
		$\begin{pmatrix}
			1_{2n} & b \\
			0  & 1_{2n}
		\end{pmatrix}$.
		Now the formulas for Weil representation
		$\omega_{W_\ell^+}$ in
		(\ref{Weil representation})
		as well as
		the value for $\gamma(q_{U,\ell})$ in
		(\ref{value for gamma})
		show that
		all these three types of matrices act as identity on
		$\phi_q=\mathbb{I}_{W^+(\Z_\ell)}$
		and thus
		\[
		\omega_{W_\ell^+}(g_\ell)\phi_\ell
		=
		\phi_\ell.
		\]

		For a place
		$\ell\mid2\delta_1\cdots\delta_{2n+1}$,
		we use the following decomposition
		(for ease of notation,we write $A,B,C,D$ for $A_\ell,B_\ell,C_\ell,D_\ell$)
		\[
		\begin{pmatrix}
			A & B \\
			C & D
		\end{pmatrix}
		=	    
		\begin{pmatrix}
			A & 0 \\
			0 & A^{-\mathrm{t}}
		\end{pmatrix}
		\begin{pmatrix}
			0 & 1_{2n} \\
			-1_{2n} & 0
		\end{pmatrix}
		\begin{pmatrix}
			1 & -A^{\mathrm{t}}C \\
			0 & 1
		\end{pmatrix}
		\begin{pmatrix}
			0 & 1_{2n} \\
			-1_{2n} & 0
		\end{pmatrix}^{-1}
		\begin{pmatrix}
			1 & A^{-1}B \\
			0 & 1
		\end{pmatrix}
		=:
		a_1a_2a_3a_4a_5.
		\]

		Assume first
		$\ell\neq2$.
		Then we have the following step by step computations:
		\begin{align*}
			\omega_{W_\ell^+}(a_5)\phi_\ell(x)
			&
			=
			\mathbb{I}_{\mathrm{Supp}_\ell(r_\ell)}(x)
			\mathbf{e}_\ell(\frac{1}{2}\mathrm{Tr}(S_xA^{-1}B))
			=
			\phi_\ell(x),
			\,
			\text{ since }
			x\in W^+(\Z_\ell),
			\\
			\omega_{W_\ell^+}(a_4a_5)\phi_\ell(x)
			&
			=
			\omega_{W_\ell^+}(a_4)\phi_\ell(x)
			=
			\gamma(q_{U,\ell})^{-1}
			\mathcal{F}(\mathbb{I}_{\mathrm{Supp}_\ell(r_\ell)})(x),
			\\
			\omega_{W_\ell^+}(a_3a_4a_5)\phi_\ell(x)
			&
			=
			\gamma(q_{U,\ell})^{-1}\mathcal{F}(\mathbb{I}_{\mathrm{Supp}_\ell(r_\ell)})(x)
			\mathbf{e}_\ell
			(-\frac{1}{2}\mathrm{Tr}(S_xA^{\mathrm{t}}C))
			=
			\gamma(q_{U,\ell})^{-1}\mathcal{F}(\mathbb{I}_{\mathrm{Supp}_\ell(r_\ell)})(x),
			\\
			\omega_{W_\ell^+}(a_2a_3a_4a_5)\phi_\ell(x)
			&
			=
			\gamma(q_{U,\ell})\gamma(q_{U,\ell})^{-1}
			\mathcal{F}^{-1}\mathcal{F}(\mathbb{I}_{\mathrm{Supp}_\ell(r_\ell)})(x)
			=\mathbb{I}_{\mathrm{Supp}_\ell(r_\ell)}(x),
			\\
			\omega_{W_\ell^+}(a_1a_2a_3a_4a_5)\phi_\ell(x)
			&
			=
			\chi_{U}^\circ(\mathrm{det}(D))
			\mathbb{I}_{\mathrm{Supp}_\ell(r_\ell)}(xA)
			=
			\chi_{U}^\circ(\mathrm{det}(D))\mathbb{I}_{\mathrm{Supp}_\ell(r_\ell)}(x),
			\text{ since }
			A\in
			\mathrm{GL}_{2n}(\Z_\ell).
		\end{align*}
		In the second step $\mathcal{F}(-)$ is the Fourier transform and it is easy to see the support of $\mathcal{F}(\mathbb{I}_{\mathrm{Supp}_\ell(r_\ell)})$ is contained in $\ell^{-2r_\ell}W^+(\Z_\ell)$.
		This proves the lemma for the place $\ell\mid\delta_U$.

		For the case $\ell=2$,
		the computation is the same as above,
		except for the first and the last steps:
		in
		$\omega_{W_\ell^+}(a_5)$,
		we have
		$B\equiv0(\mathrm{mod}\,2)$
		since $g\in\Gamma_0(2,\delta_U)$,
		thus we again have
		$\omega_{W_\ell^+}(a_5)\phi_\ell(x)=\phi_\ell(x)$;
		in
		$\omega_{W_\ell^+}(a_1a_2a_3a_4a_5)\phi_\ell(x)$,
		we have
		$\mathrm{det}(A)\equiv1(\mathrm{mod}\,4)$,
		which lies in the domain of definition of the \textit{character}
		$\chi_{U}^\circ$.
	\end{proof}

	\begin{corollary}\label{weight, level of theta lift}
		For any
		$\mathbf{f}\in\mathcal{A}_{\rho_{\lambda}}
		( H,K^{(\mathbb{S})})$,
		the theta lift 
		$\Theta_{\phi_\lambda,\mathbf{f}}$
		is a cuspidal Siegel modular form
		of weight
		$\tau^\circ$,
		of level
		$\Gamma_0(2,N_{\mathbb{S}})$
		and of character
		$\chi_{U}^\circ$.
	\end{corollary}
	\begin{proof}
		Cuspidality comes from our choice of Bruhat-Schwartz function $\phi_{\lambda}$ at places $\mathbb{S}$. It remains to deal with the level and character.
		For any
		$k
		=
		\begin{pmatrix}
			A & B \\
			C & D
		\end{pmatrix}
		\in
		\Gamma_0(2,N_{\mathbb{S}})$
		as above,
		\begin{align*}
			\Theta_{\phi_\lambda,\mathbf{f}}(gk,\varepsilon)
			&
			=
			\int\limits_{[ H]}
			\sum_{v\in W^+(\Q)}
			\langle
			\omega_{W^+}(h,(gk,\varepsilon))\phi_{\lambda}(v),\mathbf{f}(h)
			\rangle_{W,U}
			dh
			\\
			&
			=
			\varepsilon
			\int\limits_{[ H]}
			\sum_{v\in W^+(\Q)}
			\langle
			\omega_{W^+}(h,g)
			(\omega_{W^+}(1,k)\phi_{\lambda})(v),\mathbf{f}(h)
			\rangle_{W,U}
			dh
			\\
			&
			=
			\varepsilon
			\chi_{U}^\circ(\mathrm{det}\,D)
			\int\limits_{[ H]}
			\sum_{v\in W^+(\Q)}
			\langle
			\omega_{W^+}(h,g)\phi_{\lambda}(v),\mathbf{f}(h)
			\rangle_{W,U}
			dh
			\\
			&
			=
			\varepsilon
			\chi_{U}^\circ(\mathrm{det}\,D)
			\Theta_{\phi_\lambda,\mathbf{f}}(g).
		\end{align*}
	\end{proof}

	By Witt's theorem,
	one has
	\[
	H(\Q)z
	=
	\left\{
	x\in W^+(\Q)|
	S_x=S_z
	\right\}.
	\]
	We consider the following subsets
	of $H(\Q_\ell)$
	\begin{align*}
		\begin{split}
			\mathcal{E}_{z;\ell}
			&
			=
			\left\{
			h\in
			H(\Q_\ell)|
			h^{-1}z\in  \mathrm{Supp}_\ell(r_\ell)
			\right\},
			\quad
			\mathcal{E}_z
			=
			\prod_{\ell\neq\infty}
			\mathcal{E}_{z;\ell}.
		\end{split}\index{E@$\mathcal{E}_z$}
	\end{align*}
	Recall that $H(\mathbb{A}_f)_z$ is the stabilizer of $z$ in $H(\mathbb{A}_f)$.
	It it easy to see $H(\mathbb{A}_f)_z
	\cdot
	\mathcal{E}_z
	\cdot
	K^{(\mathbb{S})}
	=
	\mathcal{E}_z$ and if $\mathcal{E}_z$ is non-empty, then $z$ is of maximal rank.
	Moreover
	$H(\mathbb{A}_f)_z
	\backslash
	\mathcal{E}_z/
	K^{(\mathbb{S})}$
	is finite
	(\textit{cf.}
	\cite[Proposition 1.5]{Yoshida1984})
	and thus so is the double quotient
	\begin{equation*}
		[\mathcal{E}_{z,K^{(\mathbb{S})}}]
		:=
		H(\Q)_z
		\backslash
		\mathcal{E}_z/
		K^{(\mathbb{S})}.
	\end{equation*}
	For $S\in\mathrm{Sym}_{2n}(\Q)$
	and
	$g\in
	\widetilde{G}(\mathbb{A})$,
	one has
	\begin{align*}
		\mathbf{a}_{\Theta_{ \phi_{\lambda},\mathbf{f}},S}(g)
		&
		=
		\int\limits_{[U_G]}
		\Theta_{ \phi_{\lambda,\mathbf{f}}}(g)
		\mathbf{e}(-\mathrm{Tr}(SX))
		dn(X)
		\\
		&
		=		
		\int\limits_{[ H]}
		\sum_{v\in  W^+(\Q)}
		\langle
		\omega_{ W^+}(h,g) \phi_{\lambda}(v), \mathbf{f}(h)
		\rangle_{W,U}
		dh
		\int\limits_{[U_G]}
		\mathbf{e}(\mathrm{Tr}(S_vX)-\mathrm{Tr}(SX))dn(X)
		\\
		&
		=
		\int\limits_{[ H]}
		\sum_{v\in  W^+(\Q),\,S_v=S}
		\langle
		\omega_{ W^+}(h,g) \phi_{\lambda}(v), \mathbf{f}(h)
		\rangle_{W,U}
		dh.
	\end{align*}
	Thus if
	$\mathbf{a}_{\Theta_{ \phi_{\lambda},\mathbf{f}},S}(g)\ne0$,
	then
	$S=S_z$
	for some
	$z\in  W^+(\Q)$.
	In this case,
	one has
	\begin{align*}
		\mathbf{a}_{\Theta_{ \phi_{\lambda},\mathbf{f}},S_z}(g)
		&
		=
		\int\limits_{[ H]}
		\sum_{\gamma
			\in H(\Q)_z
			\backslash
			H(\Q)}
		\langle
		\omega_{ W^+}(g)
		\phi_{\lambda}(h^{-1}\gamma^{-1}z),\mathbf{f}(h)
		\rangle_{W,U}
		dh
		\\
		&
		=
		\int\limits_{ H(\Q)_z
			\backslash
			H(\mathbb{A})}
		\langle
		\omega_{ W^+}(g) \phi_{\lambda}(h^{-1}z),
		\mathbf{f}(h)
		\rangle_{W,U}
		dh.
	\end{align*}

	For $\xi\in\GLmr_{2n}(\mathbb{A})$, we write
	\[
	g_{\xi}:=
	\begin{pmatrix}
		\xi & 0 \\
		0 & \xi^\mathrm{-t}
	\end{pmatrix}\in G(\mathbb{A}).\index{g@$g_\xi$}
	\]
	So one gets
	\begin{equation}\label{Fourier coefficients-1}
		\mathbf{a}_{\Theta_{ \phi_{\lambda},\mathbf{f}},S_z}
		(g_\xi,\varepsilon)
		=\varepsilon
		\gamma(\mathrm{det}(\xi),q)
		\chi_{U}(\mathrm{det}\,\xi)
		|\mathrm{det}\,\xi|_{\mathbb{A}}^{(2n+1)/2}
		\cdot
		\int\limits_{ H(\Q)_{z}\backslash H(\mathbb{A})}
		\langle
		\phi_{\lambda}(h^{-1}z\xi),
		\mathbf{f}(h)
		\rangle_{W,U}
		dh.
	\end{equation}
	We write $\xi=\xi_\infty\xi_f$ and require
	$\xi_f=1$
	and
	$\mathrm{det}\,\xi_\infty>0$.
	By the definitions of
	the character
	$\chi_{U}$ and
	$\gamma$
	as well as
	the
	$H(\R)$-equivariance of the pairing
	$\langle-,-\rangle_{W,U}$,
	one gets
	\begin{align*}
		\mathbf{a}_{\Theta_{ \phi_{\lambda},\mathbf{f}},S_z}
		(g_\xi,\varepsilon)
		&
		=\varepsilon
		(\mathrm{det}\,\xi_\infty)^{(2n+1)/2}
		\int\limits_{
			H(\Q)_z
			\backslash
			H(\mathbb{A}_f)}
		\phi_f(h_f^{-1}z)
		\langle
		\phi_{\lambda,\infty}(z\xi_\infty),
		\mathbf{f}(h_f)
		\rangle_{W,U}
		dh_f
		\\
		&
		=\varepsilon
		(\mathrm{det}\,\xi_\infty)^{(2n+1)/2}
		\int\limits_{
			H(\Q)_{z}
			\backslash H(\mathbb{A}_f)}
		\mathbb{I}_{\mathcal{E}_z}(h_f)
		\langle
		\phi_{\lambda,\infty}(z\xi_\infty),
		\mathbf{f}(h_f)
		\rangle_{W,U}
		dh_f.
	\end{align*}
	By
	(\ref{action of O(U) and GL_n})
	and the definition of
	$\phi_{\lambda,\infty}$,
	we get furthermore
	\begin{align*}
		&
		\mathbf{a}_{\Theta_{ \phi_{\lambda},\mathbf{f}},S_z}
		(g_\xi,\varepsilon)
		\\
		=
		&\varepsilon
		(\mathrm{det}\,\xi_\infty)^{(2n+1)/2}
		\int
		\limits_{
			H(\Q)_{z}
			\backslash
			H(\mathbb{A}_f)}
		\mathbb{I}_{\mathcal{E}_z}(h_f)
		\rho_{\tau^\circ}(\xi_\infty^\mathrm{t})
		\left(
		\langle
		\Delta_{\lambda}(\mathcal{T}z),
		\mathbf{f}(h_f)
		\rangle_{W,U}
		\right)
		\mathbf{e}_\infty
		(i\mathrm{Tr}(S_z\xi_\infty\xi_\infty^\mathrm{t}))
		dh_f
		\\
		=
		&\varepsilon
		(\mathrm{det}\,\xi_\infty)^{(2n+1)/2}
		\mathrm{vol}(K^{(\mathbb{S})})
		\times
		\rho_{\tau^\circ}(\xi_\infty^\mathrm{t})
		\left(
		\sum_{[h_f]\in[\mathcal{E}_{z,K^{(\mathbb{S})}}]}
		w_{z,h_f}
		\langle
		\Delta_\lambda(\mathcal{T}z), \mathbf{f}(h_f)
		\rangle_{W,U}
		\right)
		\mathbf{e}_\infty
		(i\mathrm{Tr}(S_z\xi_\infty\xi_\infty^\mathrm{t})),
	\end{align*}
	where
	\begin{equation}\label{index of stailibize of z}
		w_{z,h_f}
		=
		\frac{1}{\sharp
			(
			H(\Q)_{z}
			\cap
			h_fK^{(\mathbb{S})}h_f^{-1}
			)}.\index{w@$w_{z,h_f}$}
	\end{equation}
	Expressing the above result using
	classical Siegel modular forms,
	we get
	the following
	\begin{proposition}\label{Fourier coefficient of classicla Siegel modular form}
		Let
		$f\in
		\mathcal{A}_{\rho_{\lambda}}( H,K^{(\mathbb{S})})$.
		Then the
		Fourier expansion
		for the classical Siegel modular form
		$\Theta_{\phi_{\lambda},\mathbf{f}}^\ast$
		is given by
		\[
		\Theta_{ \phi_{\lambda},\mathbf{f}}^\ast
		(Z)
		=
		\sum_{S\in\mathrm{Sym}_{2n}^\circ(\Z)}
		a_{\Theta_{ \phi_{\lambda},\mathbf{f}}}(S)
		q^S,
		\]
		where
		$S=S_z$
		for some
		$z\in W^+(\Q)$
		and
		\[
		a_{\Theta_{ \phi_{\lambda},\mathbf{f}}}(S_{z})
		=
		\frac{1}{[H(\widehat{\Z}):K^{(\mathbb{S})}]}
		\sum_{[h_f]\in[\mathcal{E}_{z,K^{(\mathbb{S})}}]}
		w_{z,h_f}
		\langle
		\Delta_{\lambda}(\mathcal{T}z),
		\mathbf{f}(h_f)
		\rangle_{W,U}.
		\]
	\end{proposition}        
	\begin{proof}
		The proof follows from
		(\ref{joint action of O(U) and GL on harmonic polynomial})
		and the relation between
		$\tau$
		and
		$\tau^\circ$.        	
	\end{proof}

	In the following we will take
	$\mathbf{z}
	\in
	\mathrm{M}_{2n+1,2n}(\Q)$
	\begin{equation}\label{z}
		\mathbf{z}
		=
		\begin{pmatrix}
			\mathbf{z}'
			\\
			0_{1\times 2n}		
		\end{pmatrix}
		\text{ with }
		\mathbf{z}'
		=
		\begin{pmatrix}
			1 & 0 & 0 & \cdots & 0 & 0 & \cdots \\
			0 & 0 & 0 & \cdots & 1 & 0 & \cdots \\
			0 & 1 & 0 & \cdots & 0 & 0 & \cdots \\
			0 & 0 & 0 & \cdots & 0 & 1 & \cdots \\
			0 & 0 & 1 & \cdots & 0 & 0 & \cdots \\
			\vdots & \vdots & \vdots & \ddots & \vdots & \vdots & \ddots
		\end{pmatrix}
		\in
		\mathrm{SL}_{2n}(\Z).\index{Z@$\mathbf{z},\mathbf{z}'$}
	\end{equation}
	One checks easily
	\begin{equation}\label{Q_Uz}
		\mathcal{T}
		\mathbf{z}
		=
		\begin{pmatrix}
			\mathbf{z}_1 & \mathbf{z}_2 \\
			\mathbf{z}_1 & -\mathbf{z}_2 \\
			0_{1\times n} & 0_{1\times n}
		\end{pmatrix},
	\end{equation}
	where
	\begin{align*}
		\mathbf{z}_1
		&
		=
		\mathrm{diag}
		(d_1,d_3,\cdots,d_{2n-1}),
		\\
		\mathbf{z}_2
		&
		=
		\mathrm{diag}
		(d_2,d_4,\cdots,d_{2n}).
	\end{align*}
	A simple computation using the definition in
	(\ref{vectors to define Delta polynomial})
	then gives:
	\begin{equation}\label{Q_Uz, computation}
		(\mathcal{T}\mathbf{z})_{\lambda} 
		=
		d_1^{\lambda_1 -\lambda_2 }
		(d_1d_3)^{\lambda_2 -\lambda_3 }
		\cdots
		(d_1d_3\cdots d_{2n-1})^{\lambda_{n} }.
	\end{equation}
	
	We define a subgroup scheme $H_{\mathbf{z}}$ of $H$ over $\Z$
	Moreover one has
	\[
	H_{\mathbf{z}}=\left\{\mathrm{diag}(1_{2n},\pm1)\right\}\simeq\{\pm1\}.
	\index{H@$H_{\mathbf{z}}$}
	\]
	In particular, for $L=\Q,\Z_v,\Q_v,\mathbb{A},\mathbb{A}_f$, we have
	\[
	H_{\mathbf{z}}(L)=H(L)_{\mathbf{z}}.
	\]

	\section{Toric periods and toric integrals}\label{toric periods and toric intergals}
	In this section we will define toric periods
	of the Fourier coefficients of
	the theta lift
	$\Theta_{\phi_{\lambda},\mathbf{f}}$
	and then relate them to certain toric periods of
	the corresponding automorphic form
	$\mathbf{f}$
	under certain conditions.

	Recall we have a maximal torus
	$T=T_1\times T_2\times\cdots\times T_{n}$ of
	$ H$.
	We define the following morphisms of algebraic groups
	\begin{align}\label{j,j'}
		\begin{split}
			\mathbf{j}'
			&
			\colon
			T
			\rightarrow
			\mathrm{SL}_{2n},
			\quad
			\mathrm{diag}(t_1,\cdots,t_{n})
			\mapsto
			(\mathbf{z}')^{-1}
			\mathrm{diag}(t_1,\cdots,t_{n})
			\mathbf{z}',
			\\
			\mathbf{j}
			&
			\colon
			T
			\rightarrow
			G,
			\quad
			t
			\mapsto
			\mathrm{diag}(\mathbf{j}'(t),\mathbf{j}'(t)^\mathrm{-t}).
		\end{split}\index{j@$\mathbf{j},\mathbf{j}'$}
	\end{align}
	So we have
	$t\mathbf{z}=\mathbf{z}\mathbf{j}'(t)$. Note that $\mathbf{j}(T(\mathbb{A}))$ splits the cover $\widetilde{G}(\mathbb{A})\rightarrow G(\mathbb{A})$.
	We write the Haar measure
	on
	$T(\mathbb{A})=T(\R)\times T(\mathbb{A}_f)$
	by
	$dt=dt_\infty dt_f$.
	Write
	$d\overline{t}$
	for the quotient measure on
	$[T]$
	induced from $dt$ on $T(\mathbb{A})$
	and similarly for
	$d\overline{t}_f$
	on
	$[T]_f$.
	We fix Haar measure
	$da=da_\infty da_f$
	on
	$ H(\mathbb{A})_{\mathbf{z}}$
	such that the volume of
	$H(\R)_{\mathbf{z}}$, resp. $H(\mathbb{A}_f)_{\mathbf{z}}$ is $1$

	\begin{definition}\label{toric period}
		For a Siegel modular form $\mathcal{F}$	of weight $\rho_{\tau^\circ}$
		(of certain level and character), a symmetric matrix $S\in\mathrm{Sym}_{2n}(\Q)$	and a continuous character of finite order
		$\psi\colon[T]\rightarrow\mathbb{C}^\times$, we define the
		\textbf{toric period} associated to the triple
		$(\mathcal{F},S,\psi)$ to be the following map
		\begin{align*}
			\mathbf{B}_{\mathcal{F},S,\psi}
			\colon
			\widetilde{G}(\mathbb{A})
			&
			\rightarrow\M_{\tau^\circ},
			\\
			g
			&			
			\mapsto
			\mathbf{e}_\infty(-i\mathrm{Tr}(S))
			\int\limits_{[T]}
			\mathbf{a}_{\mathcal{F},S}(\mathbf{j}(t)g)
			\psi(t)
			dt.
		\end{align*}
	\end{definition}
	Now we fix $\mathbf{f}\in\mathcal{A}_{\rho_{\lambda}}( H,K^{(\mathbb{S})})$.
	The toric period associated to the triple $(\Theta_{\phi_{\lambda},\mathbf{f}},S,\psi)$ can be computed for certain kind of $g\in\widetilde{G}(\mathbb{A})$ as follows: for $\xi\in\mathrm{GL}_{2n}(\mathbb{A}_f)$, as in the case of
	$\mathcal{E}_{\mathbf{z};v}$, we set
	\[
	\mathcal{E}_{\mathbf{z},\xi;v}
	=
	\left\{
	h\in H(\Q_v)|
	h^{-1}\mathbf{z}\xi_v
	\in
	\mathrm{Supp}_\ell(r_\ell)
	\right\},
	\quad
	\mathcal{E}_{\mathbf{z},\xi}
	=
	\prod_{v\neq\infty}
	\mathcal{E}_{\mathbf{z},\xi;\ell}
	\subset
	H(\mathbb{A}_f).
	\]
	By construction, we have $ H(\mathbb{A}_f)_{\mathbf{z}}	\mathcal{E}_{\mathbf{z},\xi}K^{(\mathbb{S})}=\mathcal{E}_{\mathbf{z},\xi}$.
	Moreover for any $h_f\in H({\mathbb{A}_f})$, one has
	\[
	\phi_{\lambda,f}
	(h_f^{-1}\mathbf{z}\xi)
	=
	\mathbb{I}_{\mathcal{E}_{\mathbf{z},\xi}}
	(h_f).
	\]

	For ease of notations, we will write in the following
	\[
	\mathbf{B}_{S,\psi}
	=
	\mathbf{B}_{\Theta_{\phi_{\lambda},\mathbf{f}},S,\psi}.\index{B@$\mathbf{B}_{S,\psi}$}
	\]

	Following \cite[Proposition 4.2]{HsiehNamikawa2017}, we have
	\begin{proposition}\label{expansion of toric period}
		Let
		$\mathbf{f}\in
		\mathcal{A}_{\rho_{\lambda}}( H,K^{(\mathbb{S})})$,
		$\psi
		\colon
		[T]
		\rightarrow
		\mathbb{C}^\times$
		a continuous character of finite order,
		then
		\[
		\mathbf{B}_{S_{\mathbf{z}},\psi}
		(g_{\xi},\varepsilon)
		=\varepsilon
		\gamma(\mathrm{det}(\xi),q_U)
		\chi_{U}(\mathrm{det}\,\xi)
		|\mathrm{det}\,\xi|^{(2n+1)/2}_{\mathbb{A}_f}
		\int\limits_{ H(\mathbb{A}_f)_{\mathbf{z}}
			\backslash
			\mathcal{E}_{\mathbf{z},\xi}}
		\int\limits_{[T\times H_{\mathbf{z}}]_f}
		\langle
		\Delta_{\lambda}(\mathcal{T}\mathbf{z}),
		\mathbf{f}(t_fh_f)
		\rangle_{W,U}
		\psi(t_f)
		dt_f
		dh_f.
		\]
	\end{proposition}
	\begin{proof}		
		By (\ref{Fourier coefficients-1}),
		we have
		\begin{align*}
			&
			\mathbf{e}_\infty(i\mathrm{Tr}(S_{\mathbf{z}}))
			\mathbf{B}_{S_{\mathbf{z}},\psi}
			(g_{\xi},\varepsilon)
			\\
			=
			&\varepsilon
			\int\limits_{[T]}
			\int\limits_{ H(\Q)_{\mathbf{z}}\backslash H(\mathbb{A})}
			\langle
			\omega_{ W^+}(\mathbf{j}(t)g_{\xi})
			\phi_{\lambda}(t^{-1}h^{-1}\mathbf{z}),
			\mathbf{f}(ht)
			\rangle_{W,U}
			\psi(t_f)	
			dh
			d\overline{t}
			\\
			=
			&\varepsilon
			\gamma(\mathrm{det}(\xi),q_U)
			\chi_{U}(\mathrm{det}\,\xi)
			|\mathrm{det}\,\xi|^{(2n+1)/2}_{\mathbb{A}_f}
			\int\limits_{ H(\Q)_{\mathbf{z}}\backslash H(\mathbb{A})}
			\int\limits_{[T]}
			\langle
			\phi_{\lambda}(t^{-1}h^{-1}\mathbf{z}\mathbf{j}'(t)\xi),
			\mathbf{f}(ht)
			\rangle_{W,U}
			\psi(t_f)
			d\overline{t}
			dh
			\\
			=
			&\varepsilon
			\gamma(\mathrm{det}(\xi),q_U)
			\chi_{U}(\mathrm{det}\,\xi)
			|\mathrm{det}\,\xi|^{(2n+1)/2}_{\mathbb{A}_f}
			\int\limits_{
				H(\Q)_{\mathbf{z}}
				\backslash
				H(\mathbb{A})}
			\int\limits_{[T]}
			\langle
			\phi_{\lambda}(t^{-1}h^{-1}t\mathbf{z}\xi),
			\mathbf{f}(ht)
			\rangle_{W,U}
			\psi(t_f)
			d\overline{t}
			dh.
		\end{align*}
		Now making a change of variables
		$t^{-1}ht\mapsto h$,
		noting that
		$ \phi_{\lambda,f}(h_f^{-1}\mathbf{z}\xi)
		=
		\mathbb{I}_{\mathcal{E}_{\mathbf{z},\xi}}(h_f)$
		and
		using
		(\ref{action of O(U) and GL_n}),
		Lemma
		\ref{action of O(U) and compact of Sp(V)},
		the transformation property of
		$\mathbf{f}$
		and the $H(\mathbb{R})$-equivariance of the pairing
		$\langle-,-\rangle_{W,U}$
		as in
		(\ref{pairing W,U-complex case}),
		we get
		\begin{align*}
			&
			\mathbf{e}_\infty(i\mathrm{Tr}(S_{\mathbf{z}}))
			\mathbf{B}_{S_{\mathbf{z}},\psi}(g_\xi,\varepsilon)
			\\
			=
			&\varepsilon
			\gamma(\mathrm{det}(\xi),q_U)
			\chi_{U}(\mathrm{det}\,\xi)
			|\mathrm{det}\,\xi|^{(2n+1)/2}_{\mathbb{A}_f}
			\int\limits_{ H(\Q)_{\mathbf{z}}\backslash H(\mathbb{A})}
			\int\limits_{[T]}
			\langle
			\phi_{\lambda}(h^{-1}\mathbf{z}\xi),
			\mathbf{f}(th)
			\rangle_{W,U}
			\psi(t_f)
			d\overline{t}
			dh
			\\
			=
			&\varepsilon
			\gamma(\mathrm{det}(\xi),q_U)
			\chi_{U}(\mathrm{det}\,\xi)
			|\mathrm{det}\,\xi|^{(2n+1)/2}_{\mathbb{A}_f}
			\int\limits_{ H(\Q)_{\mathbf{z}}
				\backslash
				H(\mathbb{A}_f)}
			\int\limits_{[T]_f}
			\phi_{\lambda,f}(h_f^{-1}\mathbf{z}\xi)
			\langle
			\phi_{\lambda,\infty}(\mathbf{z}),
			\mathbf{f}(t_fh_f)
			\rangle_{W,U}
			\psi(t_f)
			dt_f
			dh_f.
		\end{align*}
		Since $\phi_{\lambda,\infty}(\mathbf{z})=\mathbf{e}_\infty(i\mathrm{Tr}(S_{\mathbf{z}}))\Delta_{\lambda}(\mathcal{T}\mathbf{z})$, we get the desired formula.
	\end{proof}

	We denote the inner \emph{toric integral} in the above proposition by the following
	\[
	P_{\mathbf{f}}(\psi,h_f)
	:=
	\int\limits_{[T\times H_{\mathbf{z}}]_f}
	\langle
	\Delta_{\lambda}(\mathcal{T}\mathbf{z}),
	\mathbf{f}(t_fh_f)
	\rangle_{W,U}
	\psi(t_f)
	dt_f.\index{P@$P_{\mathbf{f}}(\psi,h_f)$}
	\]
	Note that $K^{(\mathbb{S})}$ contains $H(\mathbb{A}_f)_{\mathbf{z}}$. If $h_v$ commutes with $H(\Q_v)_{\mathbf{z}}$ for all finite places except possibly one, then one has
	\[
	P_{\mathbf{f}}(\psi,h_f)
	=
	\int\limits_{[T]_f}
	\langle
	\Delta_{\lambda}(\mathcal{T}\mathbf{z}),
	\mathbf{f}(t_fh_f)
	\rangle_{W,U}
	\psi(t_f)
	dt_f.
	\]

	We choose an auxiliary prime $\ell$
	of $\Q$ over a rational prime $\ell$
	such that
	\begin{taggedtheorem}{(A1)}\label{condition on ell}
		\begin{enumerate}
			\item
			${\ell}\nmid
			2pN_{\mathbb{S}}[T(\mathbb{A}_f):T(\Q)T(\widehat{\Z})]$
			and $\ell\equiv1(\mathrm{mod}\,4)$;

			\item
			$p\nmid(\ell-1)$;

			\item
			each rational prime factor of
			$\delta_1\cdots\delta_{2n+1}$ as well as $-1$ and $2$
			are squares in
			$\Z_\ell^\times$
			(so that the entries of
			$\mathcal{T}$
			lie in
			$\Z_\ell$
			via the embedding
			$\Q\hookrightarrow
			\Q_\ell$).	\index{A@(A1)}
		\end{enumerate}
	\end{taggedtheorem}
	
	So
	$H$ and the torus $T$ are split at $\ell$.
	We require furthermore that
	the $\ell$-th component of $K$
	satisfies
	$K_\ell= H(\Z_\ell)$.	
	For $i=1,\cdots,n$,
	recall the isomorphisms
	\[
	\mu_i
	\colon
	T_i(\Q_\ell)
	\simeq
	\Q_\ell^\times,
	\quad
	\begin{pmatrix}
		a_i & -b_id_{2i}^2/d_{2i-1}^2 \\
		b_i & a_i
	\end{pmatrix}
	\mapsto
	a_i+b_id_{2i}/d_{2i-1},
	\]
	which maps
	$T_i(\Z_\ell)$
	onto
	$\Z_\ell^{\times}$.
	Then we have an isomorphism
	\[
	\mu\colon
	T(\Q_\ell)
	\simeq
	(\Q_\ell^\times)^{n},\quad
	t=(t_1,\cdots,t_{n})
	\mapsto
	(\mu_1(t_1),\cdots,\mu_{n}(t_{n})),\index{m@$\mu$}
	\]
	mapping
	$T(\Z_\ell)$
	onto
	$(\Z_\ell^\times)^{n}$.
	One checks easily that for any
	$t\in T(\Q_\ell)$:
	\begin{equation*}
		\mathcal{T}
		t
		\mathcal{T}^{-1}
		=
		\begin{pmatrix}
			\mathrm{diag}(\mu_1(t_1),\cdots,\mu_{n}(t_{n})) & & \\
			& \mathrm{diag}(\mu_1(t_1),\cdots,\mu_{n}(t_{n}))^{-\mathrm{t}} & \\
			& & 1_{1\times1}
		\end{pmatrix}.
	\end{equation*}
	We define a subgroup of $T(\Z_\ell)$
	\[
	T(\Z_\ell)^\circ
	=
	\left\{
	t=(t_1,\cdots,t_n)\in T(\Z_\ell)|
	\mu_i(t_i)\equiv1(\mathrm{mod}\,\ell)\,\forall i
	\right\}.\index{T@$T(\Z_\ell)^\circ$}
	\]
	Then
	$\mu(T(\Z_\ell)^\circ)
	=
	(1+\ell\Z_\ell)^{n}$.
	We consider the following automorphism of
	$(\Z_\ell^\times)^{n}$
	\begin{align}\label{sigma automorphism}
		\begin{split}
			\sigma'
			\colon
			(\Z_\ell^\times)^{n}
			&
			\rightarrow
			(\Z_\ell^\times)^{n}
			\\
			t'=(t_1',\cdots,t_{n}')
			&
			\mapsto
			\sigma'(t')
			=
			(\sigma'_1(t'),\cdots,\sigma'_{n}(t'))
			=
			(t_1'/t_{n}',\cdots,t_{n-1}'/t_{n}',
			t_{n}'),
		\end{split}		
	\end{align}
	which induces the following automorphism of
	$T(\Z_\ell)$:
	\begin{equation}\label{sigma}
		\sigma
		\colon
		T(\Z_\ell)
		\xrightarrow{\mu}
		(\Z_\ell^\times)^{n}
		\xrightarrow{\sigma'}
		(\Z_\ell^\times)^{n}
		\xrightarrow{\mu^{-1}}
		T(\Z_\ell).\index{s@$\sigma,\sigma'$}
	\end{equation}

	\begin{remark}\rm
		The particular form of the automorphism
		$\sigma'$
		comes from the adjoin action of
		$T( \Q_\ell)$
		on certain unipotent subgroups of
		$H( \Q_\ell)$
		that we will consider later on.
	\end{remark}

	For an integer $r\ge0$,
	we put
	\begin{equation}\label{sigma_i(T)_f}
		T(\Z_\ell)_{i,r}
		=
		(\sigma'\circ\mu)^{-1}
		\{
		(
		\underbrace{1,\cdots,1}_\text{$i-1$},
		1+\ell^r\Z_\ell,
		\underbrace{1,\cdots,1}_\text{$n-i$}
		)
		\},
		\index{T@$T(\Z_\ell)_{i,r},T(\Z_\ell)_{\underline{r}},T(\widehat{\Z})_{\underline{r}}$}
	\end{equation}
	where for $r=0$,  we replace $1+\ell^0\Z_\ell$ by $\Z_\ell^\times$.
	We have the following induced isomorphism:
	\begin{equation}
		\widetilde{\mu}_i
		\colon
		T(\Z_\ell)_{i,r}
		\simeq
		1+\ell^r\Z_\ell,
		\quad
		t
		\mapsto
		\sigma'_i(\mu(t)).
	\end{equation}	
	For an $n$-tuple
	$\underline{r}=(r_1,r_2,\cdots,r_{n})
	\in
	\mathbb{N}^{n}$,
	we define the following compact open subgroup of
	$T(\Z_\ell)$ and 
	$T(\widehat{\Z})$
	respectively,
	\begin{equation*}
		T(\Z_\ell)_{\underline{r}}
		=
		\prod_{i=1}^{n}
		T(\Z_\ell)_{i,r_i},
		\quad
		T(\widehat{\Z})_{\underline{r}}
		=
		\{
		t
		\in
		T(\widehat{\Z})
		\mid
		t_\ell
		\in
		T(\Z_\ell)_{\underline{r}}
		\}.
	\end{equation*}

	We define the following quotient groups
	of
	$T(\mathbb{A}_f)$:
	\[
	\mathcal{G}(\infty)
	=
	T(\Q)\backslash
	T(\mathbb{A}_f)
	/    
	T(\widehat{\Z}^\ell),
	\quad
	\mathcal{G}(\infty)'
	=
	T(\mathbb{A}_f^\ell)
	/
	T(\widehat{\Z}^\ell).\index{G@$\mathcal{G}(\infty),\mathcal{G}(\infty)'$}
	\]
	Here $T(\mathbb{A}_f^\ell)$ means that the components at the place $\ell$
	are trivial. Note that $\mathcal{G}(\infty)$ is a compact group. We view
	$\mathcal{G}(\infty)'$ as a subgroup of $\mathcal{G}(\infty)$
	via the inclusion $T(\mathbb{A}_f^\ell)\hookrightarrow T(\mathbb{A}_f)$.
	We next consider the torsion subgroups of
	$\mathcal{G}(\infty)$ and $\mathcal{G}(\infty)'$:
	\[
	\mathcal{G}_1
	=
	\mathcal{G}(\infty)'_{\mathrm{tor}}
	\subset
	\mathcal{G}_0
	=
	\mathcal{G}(\infty)_{\mathrm{tor}}.\index{G@$\mathcal{G}_0,\mathcal{G}_1$}
	\]
	It is not hard to see $\mathcal{G}_1=\{1\}$ because $(T(\Q_q)/T(\Z_q))_{\mathrm{tor}}=\{1\}$ for any prime $q$. The theory of topological groups shows that the quotient 
	$\mathcal{G}(\infty)/\mathcal{G}_0$ is	isomorphic to $\Z_\ell^{n}$. The assumption
	$\ell\nmid[T(\mathbb{A}_f):T(\Q)T(\widehat{\Z})]$ implies that the following exact sequence of compact groups splits (because the last term has cardinal prime to $\ell$ while the first term is a pro-$\ell$ group):
	\begin{equation}\label{splitting of compact groups}
		1
		\rightarrow
		T(\Z_\ell)_{(1,\cdots,1)}
		\rightarrow
		\mathcal{G}(\infty)
		\rightarrow
		T(\mathbb{A}_f)/T(\Q)T(\widehat{\Z})_{(1,\cdots,1)}
		\rightarrow
		1.
	\end{equation}	
	Then it follows immediately
	\begin{lemma}\label{torsion subgroup isomorphic to quotient group}
		The composition map
		$\mathcal{G}_0
		\rightarrow
		\mathcal{G}(\infty)
		\rightarrow
		T(\Q)\backslash
		T(\mathbb{A}_f)/T(\widehat{\Z})_{(1,\cdots,1)}$
		is an isomorphism.
	\end{lemma}
	\begin{proof}
		Clearly this map is injective since
		$T(\Z_\ell)_{(1,\cdots,1)}
		\simeq
		\Z_\ell^{n}$
		is torsion-free.
		However
		the sequence
		(\ref{splitting of compact groups})
		splits,
		we may view
		$T(\Q)\backslash
		T(\mathbb{A}_f)/T(\widehat{\Z})_{(1,\cdots,1)}$
		as a subgroup of
		$\mathcal{G}(\infty)$,
		which is a finite subgroup, thus a torsion subgroup
		and therefore contained in
		$\mathcal{G}_0$.
		Moreover the composition map
		\[
		T(\Q)\backslash T(\mathbb{A}_f)/T(\widehat{\Z})_{(1,\cdots,1)}
		\rightarrow
		\mathcal{G}(\infty)
		\rightarrow
		T(\Q)\backslash
		T(\mathbb{A}_f)/T(\widehat{\Z})_{(1,\cdots,1)}
		\]
		is the identity map.
		We conclude that
		the composition map in the lemma is an isomorphism.
	\end{proof}

	It follows from the above description of
	$\mathcal{G}(\infty)$
	and
	$\mathcal{G}(\infty)'$
	that
	$\mathcal{G}_0$
	is a finite subgroup of
	$\mathcal{G}(\infty)$
	whose cardinal is prime to $\ell$
	and the composition map
	\[
	T(\Z_\ell)_{(1,\cdots,1)}
	\hookrightarrow
	\mathcal{G}(\infty)
	\rightarrow
	\Gamma
	:=
	\mathcal{G}(\infty)/\mathcal{G}_0
	\]
	is an isomorphism.
	We summarize the relations among the above mentioned groups
	in the following commutative diagram,
	where the middle vertical and horizontal sequences are exact:
	\begin{equation*}
		\begin{tikzcd}
			\mathcal{G}_1=1
			\arrow[r , hookrightarrow]
			\arrow[d , hookrightarrow]
			&
			\mathcal{G}_0
			\arrow[d , hookrightarrow]
			\arrow[dr , "(\ast)", hookrightarrow]
			&
			\\
			\mathcal{G}(\infty)'
			\arrow[r , hookrightarrow]
			\arrow[d, "="]
			&
			\mathcal{G}(\infty)
			\arrow[r , twoheadrightarrow]
			\arrow[d , twoheadrightarrow]
			&
			T( \Q_\ell)/T(\Q)
			\\
			\mathcal{G}(\infty)'
			\arrow[r , hookrightarrow]
			&
			\Gamma
			&
		\end{tikzcd}
	\end{equation*}
	The injectivity of the diagonal arrow $(\ast)$ follows from
	the fact
	$\mathcal{G}(\infty)'\cap\mathcal{G}_0=\mathcal{G}_1=\{1\}$.
	We can thus view
	$\mathcal{G}_0$
	as a subgroup of
	$T( \Q_\ell)/T(\Q)$.	
	We fix then a (non-canonical)
	decomposition
	\[
	\mathcal{G}(\infty)
	\simeq
	\mathcal{G}_0
	\times
	\Gamma.
	\]
	For
	an $n$-tuple of positive integers
	$\underline{r}
	=
	(r_1,r_2,\cdots,r_{n})$,
	we define the following finite groups:
	\[
	\Gamma_i(r_i)
	=
	T(\Z_\ell)_{i,1}/
	T(\Z_\ell)_{i,r_i},
	\quad
	\Gamma(\underline{r})
	=
	\prod_{i=1}^{n}
	\Gamma_i(r_i).\index{G@$\Gamma_i(r_i),\Gamma(\underline{r})$}
	\]
	Then the decomposition
	$\mathcal{G}(\infty)
	\simeq
	\mathcal{G}_0\times\Gamma$
	induces an isomorphism
	\begin{equation}\label{splitting of G(infty)}
		\mathcal{G}(\infty)/T(\Z_\ell)_{\underline{r}}
		\simeq
		\mathcal{G}_0
		\times
		\Gamma(\underline{r}).
	\end{equation}

	We consider the following unipotent elements in
	$H(\Q_\ell)$:
	for any $v\in \Q_\ell$
	and $i=1,\cdots,n$,
	we define
	$u_i(v)\in H(\Q_\ell)$ to be
	\begin{equation}\label{Q_Uu_i(t)Q_U^-1}
		u_i(v)
		=
		\begin{cases*}
			\mathcal{T}^{-1}
			\begin{pmatrix}
				1_{n}-vE_{n,i} & 0 & 0\\
				0 & 1_{n}+vE_{i,n} & 0\\
				0& 0 & 1_{1\times1}
			\end{pmatrix}
			\mathcal{T},
			&
			$i=1,\cdots,n-1$;
			\\
			\mathcal{T}^{-1}
			\begin{pmatrix}
				1_{n} & \frac{v^2}{2}E_{n,n}  &-vE_{n} \\
				0 & 1_{n} & 0 \\
				0 & -vE_{n}^{\mathrm{t}}& 1_{1\times1}
			\end{pmatrix}
			\mathcal{T},
			&
			$i=n$.
		\end{cases*}\index{u@$u_i(v)$}
	\end{equation}
	Here $E_{i,j}$
	are the elementary matrices of size $n\times n$, $E_{n}$ the column matrix of size $n\times1$ with $1$ on the last entry and $0$ elsewhere.
	One checks easily that these unipotent elements commute with each other:
	\[
	u_i(v)u_j(w)=u_j(w)u_i(v),
	\quad
	\forall
	i,j\in\{1,\cdots,n\},
	\,
	v,w\in
	\Q_\ell.
	\]
	The conjugate action $\mathrm{Ad}$ of
	$T( \Q_\ell)$ on
	these $u_i(v)$
	is given as follows:
	write
	$t=(t_1,\cdots,t_{n})
	\in
	\prod_{i=1}^{n}T_i
	( \Q_\ell)$, then
	\begin{equation*}		
		\mathrm{Ad}_t(u_i(v))=
		tu_i(v)t^{-1}
		=
		\begin{cases*}
			u_i(v\mu_{n}(t_{n})/\mu_i(t_i)),
			&
			$i=1,\cdots,n-1$;
			\\
			u_i(v\mu_{n}(t_{n})),
			&
			$i=n$.
		\end{cases*}
	\end{equation*}
	For
	$t=(t_1,\cdots,t_{n})
	\in
	\prod_{i=1}^{n}
	\sigma_i(T(\Z_\ell))_1$,
	we write
	$(\sigma'\circ\mu)(t)_i$
	for the $i$-th component of
	$(\sigma'\circ\mu)(t)$.
	Then we have
	\begin{equation}\label{adjoint action of T on u, another basis}
		tu_i(v)t^{-1}
		=
		u_i((\sigma'\circ\mu)(t)_iv),
		\quad
		i=1,\cdots,n.
	\end{equation}

	We consider an element
	$\xi=\xi(\underline{r})
	\in
	\mathrm{GL}_{2n}(\mathbb{A}_f)$\index{x@$\xi=\xi(\underline{r})$}
	such that the matrix
	$\mathbf{z}\xi=((\mathbf{z}\xi)_q)_q
	\in
	\mathrm{M}_{2n+1,2n}
	(\mathbb{A}_f)$
	is of the form:
	\begin{equation}\label{L=(L_1,...L_2n)}
		(\mathbf{z}\xi)_{\ell}
		=
		\begin{pmatrix}
			\widetilde{\mathbf{z}\xi}_{\ell}
			\\
			0_{1\times2n}
		\end{pmatrix}
		=
		\begin{cases*}
			\mathbf{z}_{\ell},
			&
			${q}\ne\ell$;
			\\
			\mathcal{T}^{-1}
			\prod_{i=1}^{n-1}
			\begin{pmatrix}
				1_{n}+\ell^{-r_i}E_{n,i}  & 0 \\
				0 & 1_{n}-\ell^{-r_i}E_{i,n} \\
				0_{1\times n} & 0_{1\times n}
			\end{pmatrix}
			\begin{pmatrix}
				1_{n} & -\frac{\ell^{-r_{n}}}{2}E_{n,n} \\
				0 & 1_{n-1} \\
				0 & \ell^{r_{n}} \\
				0_{1\times n} & 0_{1\times n}
			\end{pmatrix},
			&
			${q}=\ell$.
		\end{cases*}
	\end{equation}
	Note that $\xi$ is uniquely determined by the expression for $\mathbf{z}\xi$.

	We fix $\underline{r}$ and assume the following condition
	\begin{taggedtheorem}{(A2)}\label{condition on the character psi}
		The character
		$\psi
		\colon
		[T]
		\rightarrow
		\mathbb{C}^\times$
		is trivial on $T(\mathbb{R})$,
		and is the product of $n$ characters
		$\psi_i\colon
		[T_i]
		\rightarrow
		\mathbb{C}^\times$
		for $i=1,\cdots,n$
		where the conductors of
		$\psi_1,\psi_2,\cdots,\psi_{n}$
		are		
		$\ell^{r_1},\ell^{r_2},
		\cdots,\ell^{r_{n}}$,
		respectively.\index{A@(A2)}
	\end{taggedtheorem}

	So we view $\psi$ also as a character of
	$[T]_f$.
	We then write
	\[
	N(\underline{r})
	=
	[
	T(\Q)T(\widehat{\Z}):
	T(\Q)T(\widehat{\Z})_{\underline{r}}
	].\index{N@$N(\underline{r})$}
	\]
	For any element
	$h_f\in H(\mathbb{A}_f)$
	with
	$h^{-1}_f
	T(\widehat{\Z})_{\underline{r}}
	h_f
	\subset
	K^{(\mathbb{S})}$,
	we put
	\[
	\mathbf{P}_{\mathbf{f}}(\psi,h_f)
	=
	\sum_{
		t
		\in
		[T]_f
		/T(\widehat{\Z})_{\underline{r}}}
	\langle
	\Delta_{\lambda}(\mathcal{T}\mathbf{z}),
	\mathbf{f}(t_fh_f)
	\rangle_{W,U}
	\psi(t_f).\index{P@$\mathbf{P}_{\mathbf{f}}(\psi,h_f)$}
	\]
	It follows by definition
	\[
	P_{\mathbf{f}}(\psi,h_f)
	=
	N(\underline{r})^{-1}
	\mathbf{P}_{\mathbf{f}}(\psi,h_f).
	\]

	By the assumptions
	$\ell\neq p$
	and
	$p\nmid(\ell-1)$,
	we have
	\begin{equation}\label{p not dividing t_r}
		p\nmid  N(\underline{r}).
	\end{equation}

	For $k=0,1,\cdots,r_n$, we define an element $\varsigma_{\underline{r},k}=((\varsigma_{\underline{r},k})_q)_q\in H(\mathbb{A}_f)$\index{s@$\varsigma_{\underline{r},k}$} as follows
	\begin{equation}\label{varsigma}
		(\varsigma_{\underline{r},k})_q
		=
		\begin{cases*}
			1,
			&
			$q\neq\ell$;
			\\
			\mathcal{T}^{-1}
			\prod_{i=1}^{n}u_i(- q^{-r_i})
			\mathrm{diag}(1_{n-1},\ell^{-k},1_{n-1},\ell^{k},1)
			\mathcal{T},
			&
			$q=\ell$.
		\end{cases*}
	\end{equation}

	For a prime number $q$, we consider the following sets:
	\begin{align*}
		\begin{split}
			[\mathcal{E}_{\mathbf{z},\xi;{q}}]
			&
			=
			H(\Q_q)_{\mathbf{z}}
			\backslash
			\mathcal{E}_{\mathbf{z},\xi;{q}}/K^{(q)},
			\\
			\widetilde{\mathcal{E}}_{\mathbf{z},\xi;{q}}
			&
			=
			\{(\varsigma_{\underline{r},k})_q\mid k=0,\cdots,r_{n}\}.
		\end{split}
	\end{align*}
	
	Then we set
	\[
	\widetilde{\mathcal{E}}_{\mathbf{z},\xi}
	:=
	\prod_{q}
	\widetilde{\mathcal{E}}_{\mathbf{z},\xi;{q}}
	=
	\{\varsigma_{\underline{r},k}\mid k=0,\cdots,r_{n}\}.\index{E@$\widetilde{\mathcal{E}}_{\mathbf{z},\xi},\widetilde{\mathcal{E}}_{\mathbf{z},\xi;q},[\mathcal{E}_{\mathbf{z},\xi;{q}}]$}
	\]
	\begin{proposition}\label{E_z set of representatives}
		\begin{enumerate}[label=(\alph*)]
			\item
			The set
			$\widetilde{\mathcal{E}}_{\mathbf{z},\xi;{q}}$
			as above
			is a complete set of representatives for
			$[\mathcal{E}_{\mathbf{z},\xi;{q}}]$
			for any prime
			$q$.
			Therefore
			the set
			$\widetilde{\mathcal{E}}_{\mathbf{z},\xi}$
			is a complete
			set of representatives for the double quotient
			$
			H(\mathbb{A}_f)_{\mathbf{z}}
			\backslash
			\mathcal{E}_{\mathbf{z},\xi}/
			K^{(\mathbb{S})}$.

			\item
			Fix $k\in\{0,\cdots,r_{n}\}$ and define an $n$-tuple of integers
			\[
			\underline{r}(k)=(\max(r_1-k,0),\cdots,\max(r_{n}-k,0)).
			\]
			Then for
			$t\in
			T
			(\widehat{\Z})_{\underline{r}(k)}$,
			we have
			\[
			\varsigma_{\underline{r},k}^{-1}t\varsigma_{\underline{r},k}
			\in
			H(\Z_\ell)
			\times
			\prod_{q\ne\ell}
			T(\Z_{\ell}).
			\]
		\end{enumerate}

	\end{proposition}
	\begin{proof}
		Our argument is inspired from \cite[§4.3]{HsiehNamikawa2017}, in particular Lemma 4.3 and Proposition 4.6 of \emph{loc.cit}.

		For $(a)$,
		we consider the following cases separately:
		\begin{enumerate}
			\item 
			$
			{q}\notin\mathbb{S}\cup\{\ell\}$.
			In this case,
			$\mathcal{E}_{\mathbf{z},\xi;{q}}\neq\emptyset$
			and moreover
			$[\mathcal{E}_{\mathbf{z},\xi;{q}}]$
			has only one element
			by
			\cite[Proposition 1.3]{Yoshida1984}.
			We can just take
			$\widetilde{\mathcal{E}}_{\mathbf{z},\xi;q}=\{1\}$ as the set of representatives for $[\mathcal{E}_{\mathbf{z},\xi;{q}}]$.

			\item
			$q\in\mathbb{S}$. We write $\widetilde{\Lambda}=\mathrm{diag}(\delta_1,\cdots,\delta_{2n})$. So $\Lambda=\mathrm{diag}(\widetilde{\Lambda},\delta_{2n+1})$. For an element $g\in H(\Q_q)$, write $g^{-1}=\begin{pmatrix}
				A & B \\
				C & D
			\end{pmatrix}$ in block matrix with $A$ of size $2n\times 2n$. Then $g^{-1}\mathbf{z}=\begin{pmatrix}
				A\mathbf{z}' \\
				C\mathbf{z}'
			\end{pmatrix}\in\mathrm{Supp}_q(r_{q})$ implies that $A\in\mathrm{GL}_{2n}(\Z_q)$ and $C\in q^{r_{q}}\mathrm{M}_{1,2n}(\Z_q)$. On the other hand, $g^{-1}\in H(\Q_v)$ implies
			\[
			B=-(A^{\mathrm{t}}\widetilde{\Lambda})^{-1}C^{\mathrm{t}}\delta_{2n+1}D,
			\quad
			B^{\mathrm{t}}\widetilde{\Lambda}B+D^{\mathrm{t}}\delta_{2n+1}D=\delta_{2n+1}.
			\]
			Then our assumption on $r_{q}$ shows that
			\[
			B^{\mathrm{t}}\widetilde{\Lambda}B
			=
			D^2\delta_{2n+1}(\delta_{2n+1}(CA^{-1})\widetilde{\Lambda}^{-1}(CA^{-1})^{\mathrm{t}})\in qD^2\delta_{2n+1}\Z_q.
			\]
			Therefore we must have $D\in\Z_q^\times$ and $B\in q\mathrm{M}_{2n,1}(\Z_q)$. Therefore $g^{-1}\in K^{(q)}$.

			\item 
			${q}=\ell$. We first show that the problem can be reduced to the case $n=3$.
			Write
			\[
			\widetilde{\xi}(r_{n})
			=\prod_{i=1}^{n-1}u_i(\ell^{-r_i})\mathbf{z}\xi_\ell
			=
			\mathcal{T}^{-1}
			\begin{pmatrix}
				1_{n} & -\frac{\ell^{-r_{n}}}{2}E_{n,n} \\
				0 & \ell^{r_{n}}\cdot1_{n} \\
				0_{1\times n} & 0_{1\times n}
			\end{pmatrix}
			\]
			and for an element $g\in H(\Q_\ell)$, write
			\begin{equation}\label{widetilde{g}}
				\widetilde{g}=\prod_{i=1}^{n-1}u_i(\ell^{-r_i})g.
			\end{equation}
			Then for an element $g\in H(\Q_\ell)$, $g^{-1}\mathbf{z}\xi_\ell\in W^+(\Z_\ell)$ if and only if $\widetilde{g}^{-1}\widetilde{\xi}(r_{n})
			\in W^+(\Z_\ell)$. So
			\[
			\left\{
			g\in H(\Q_\ell)\mid g^{-1}\mathbf{z}\xi_\ell\in W^+(\Z_\ell)
			\right\}
			=
			\prod_{i=1}^{n-1}u_i(-\ell^{-r_i})
			\left\{
			\widetilde{g}\in H(\Q_\ell)\mid\widetilde{g}^{-1}\widetilde{\xi}(r_{n})\in W^+(\Z_\ell)
			\right\}.
			\]
			We rearrange the order of the elements in the basis $\widetilde{\mathfrak{B}}$ and put
			\[
			\widetilde{\mathfrak{B}}'
			=
			(\widetilde{E}_1,\cdots,\widetilde{E}_{n-1},\widetilde{E}_{n+1},\cdots,\widetilde{E}_{2n-1},\widetilde{E}_{n},\widetilde{E}_{2n},\widetilde{E}_{2n+1}).
			\]
			Under this basis, the quadratic form $\langle-,-\rangle_{U}$ is represented by the matrix
			\[
			\widetilde{\Lambda}':=\begin{pmatrix}
				& 1_{n-1} & & & \\
				1_{n-1} & & & & \\
				& & & 1 & \\
				& & 1 & & \\
				& & & & 1
			\end{pmatrix}
			\]
			We write $\mathcal{T}'$ for the transformation matrix from the basis $\mathfrak{B}$ to $\widetilde{\mathfrak{B}}'$. Then one gets
			\[
			\mathcal{T}'\widetilde{\xi}(r_{n})
			=\begin{pmatrix}
				1_{2n-2} & &  \\
				& 1 & -\frac{\ell^{-r_{n}}}{2} \\
				& 0 & \ell^{r_{n}} \\
				& 0 & 0
			\end{pmatrix}.
			\]
			In the following we write $\mathrm{O}_\Lambda(\Q_\ell)$ for the orthogonal group consisting of $g\in\mathrm{GL}_n(\Q_\ell)$ such that $g^{\mathrm{t}}\Lambda g=\Lambda$. Then
			\[
			\left\{g\in \mathrm{O}_\Lambda(\Q_\ell)\mid g^{-1}\widetilde{\xi}(r_{n})\in W^+(\Z_\ell)
			\right\}
			=
			\mathcal{T}'
			\left\{g\in
			\mathrm{O}_{\widetilde{\Lambda}'}(\Q_\ell)\mid g^{-1}\mathcal{T}'\widetilde{\xi}(r_{n})\in W^+(\Z_\ell)
			\right\}
			(\mathcal{T}')^{-1}.
			\]
			As a result we are reduced to finding those $g\in\mathrm{O}_{\widetilde{\Lambda}'}(\Q_\ell)$ such that $g^{-1}\mathcal{T}'\widetilde{\xi}(r_{n})\in W^+(\Z_\ell)$. For this, write $g^{-1}=\begin{pmatrix}
				A & B \\
				C & D
			\end{pmatrix}$ with $A$ of size $(2n-2)\times(2n-2)$. Then $A,C$ both have entries in $\Z_\ell$. Write $L_1$ for the $\mathbb{Z}_\ell$-submodule of $U(\Z_\ell)$ generated by the columns of $\begin{pmatrix}
				A \\ C
			\end{pmatrix}$. Then by \cite[82:15]{Omeara1963}, $L_1$ splits $U(\Z_\ell)$ (say, by another lattice $L_2$). Now we have two splittings of $U(\Z_\ell)$:
			\[
			U(\Z_\ell)
			=L_1\oplus L_2
			=\Z_\ell(\widetilde{E}_1,\cdots,\widetilde{E}_{n-1},\widetilde{E}_{n+1},\cdots,\widetilde{E}_{2n-1})\oplus\Z_\ell(\widetilde{E}_{n},\widetilde{E}_{2n},\widetilde{E}_n)
			\]
			and the the first factors in both splittings are isometric (by the map represented by the matrix $\begin{pmatrix}
				A \\ C
			\end{pmatrix}$). So by \cite[92:3]{Omeara1963}, the second factors $L_2$ and $\Z_\ell(\widetilde{E}_{n}^+,\widetilde{E}_{n}^-,E_n)$ are also isometric (say, by a map represented by matrix $\begin{pmatrix}
				B' \\ D'
			\end{pmatrix}$). Then the matrix $\begin{pmatrix}
				A & B' \\
				C & D'
			\end{pmatrix}$ lies in $\mathcal{O}(\widetilde{\Lambda}')(\Q_\ell)$.
			Write 
			\[
			\widehat{\Lambda}=\begin{pmatrix}
				0 & 1 &0 \\ 1 &0 &0 \\ 0& 0& 1
			\end{pmatrix}
			\]
			and $\mathrm{O}_{\widehat{\Lambda}}$ for the orthogonal matrix associated to $\widehat{\Lambda}$. Then there exists $h\in\mathrm{O}_{\widehat{\Lambda}}(\Q_\ell)$ such that
			\[
			\begin{pmatrix}
				A & B' \\ C & D'
			\end{pmatrix}^{-1}
			\begin{pmatrix}
				A & B \\ C & D
			\end{pmatrix}=
			\begin{pmatrix}
				1_{2n-2} &0  \\0 & h
			\end{pmatrix}
			\]
			Now we put
			\[
			\widehat{\mathbf{z}}=\begin{pmatrix}
				1 & -\ell^{-r_{n}}/2 \\
				0 & \ell^{r_{n}} \\
				0 & 0
			\end{pmatrix}
			\]
			We write
			\[
			\mathcal{E}_{\widehat{\mathbf{z}}}
			=
			\left\{
			\begin{pmatrix}
				1_{2n-2} &0 \\
				0 & h
			\end{pmatrix}
			\mid
			h\in\mathrm{O}_{\widehat{\Lambda}}(\Q_\ell),\,
			h^{-1}\widehat{\mathbf{z}}\in\mathrm{M}_{3,2}(\Z_\ell)
			\right\}.
			\]
			Then it is easy to check
			\[
			\left\{
			g\in\mathrm{O}_{\widetilde{\Lambda}'}(\Q_\ell)\mid
			g^{-1}\mathcal{T}'\widetilde{\xi}(r_{n})\in W^+(\Z_\ell)
			\right\}
			=
			\mathrm{O}_{\widetilde{\Lambda}'}(\Z_\ell)
			\mathcal{E}_{\widehat{\mathbf{z}}}.
			\]
			So the problem is reduced to find a set of representatives for the double quotient 
			\[
			\mathrm{O}_{\widehat{\Lambda}}(\Q_\ell)_{\widehat{\mathbf{z}}}
			\backslash
			\mathcal{E}_{\widehat{\mathbf{z}}}/\mathrm{O}_{\widehat{\Lambda}}(\Z_\ell).
			\]
			As a result we are reduced to the case $n=3$.

			So in the following we will assume $n=3$. In this case one has
			\[
			u_{n}
			(\ell^{-r_{n}})\mathbf{z}\xi_\ell
			=\begin{pmatrix}
				1 & 0 \\ 0 & \ell^{r_{n}} \\ 0 & 1
			\end{pmatrix}.
			\]

			For $g\in\mathrm{O}_{\widehat{\Lambda}}(\Q_\ell)$ such that $g^{-1}\in\mathcal{E}_{\mathbf{z},\xi;\ell}$, we write $(a_{i,j})_{i,j=1}^3=g^{-1}$. Then up to multiplication on the left by the element $\widehat{\Lambda}\in H(\Z_\ell)$, we can assume $0\le\mathrm{ord}_\ell(a_{1,1})\le\mathrm{ord}_\ell(a_{2,1})$. We write
			\begin{align*}
				k
				&
				=\min(\mathrm{ord}_\ell(a_{1,1}),r_{n}),
				\\
				\widehat{g}
				&
				=\mathrm{diag}(\ell^{k},\ell^{-k},1)g,
				\\
				(\widehat{a}_{i,j})_{i,j=1}^3
				&
				=\widehat{g}^{-1}.
			\end{align*}
			Note that  we have
			\begin{align}\label{relation among widehat{a}_{i,j}}
				\begin{split}
					2\widehat{a}_{1,1}\widehat{a}_{1,2}+\widehat{a}_{1,3}^2=0
					\\
					\widehat{a}_{1,1}\widehat{a}_{2,2}+\widehat{a}_{1,2}\widehat{a}_{2,1}+\widehat{a}_{1,3}\widehat{a}_{2,3}
					=1
				\end{split}
			\end{align}

			We next show $\mathrm{ord}_\ell(a_{1,1})\le r_{n}$. It follows from $\widehat{g}^{-1}\begin{pmatrix}
				\ell^{r_{n}} & 0 \\ 0 & 1 \\ 0 & 1
			\end{pmatrix}\in W^+(\Z_\ell)$ that $\widehat{a}_{1,2}+\widehat{a}_{1,3}\in\Z_\ell$. If $\mathrm{ord}_\ell(a_{1,1})>r_{n}$, then $\widehat{a}_{1,1}=a_{1,1}\ell^{-r_{n}}\in\ell\Z_\ell$, thus we have $\widehat{a}_{1,2}\in\Z_\ell$ and $\widehat{a}_{1,3}\in\ell\Z_\ell$. The same holds for $\widehat{a}_{2,1},\widehat{a}_{2,2}$ and $\widehat{a}_{2,3}$, which contradicts (\ref{relation among widehat{a}_{i,j}}). So we must have
			\[
			\mathrm{ord}_\ell(a_{1,1})\le r_{n}.
			\]
			We deduce $\widehat{a}_{1,1}\in\Z_\ell^\times$. Note that the stabilizer of $\mathbf{z}\xi_\ell$ in $H(\Q_\ell)$ is $\{\mathrm{diag}(1,1,\pm1)\}$. Therefore we have
			\[
			H(\Q_\ell)_{\mathbf{z}\xi_\ell}
			=
			\left\{
			1_3,
			\begin{pmatrix}
				1 & 2\ell^{-2r_{n}+2k} & -2\ell^{-r_{n}+k} \\
				0 & 1 & 0 \\
				0 & 2\ell^{-r_{n}+k} & -1
			\end{pmatrix}
			\right\}
			\]
			So up to multiplying $\widehat{g}^{-1}$ on the right by the non-trivial element in $H(\Q_\ell)_{\mathbf{z}\xi_\ell}$, we can assume that $\widehat{a}_{1,3}\in\Z_\ell$. From (\ref{relation among widehat{a}_{i,j}}) we deduce that $\widehat{a}_{1,2}\in\Z_\ell$. In other words, we can assume that the first row of $\widehat{g}^{-1}$ has entries all in $\Z_\ell$.

			We then proceed to show that the remaining entries of $\widehat{g}^{-1}$ are also in $\Z_\ell$. We know already $\widehat{a}_{2,1}\in\Z_\ell$. Now if $\mathrm{ord}_\ell(\widehat{a}_{2,3})=r<0$, then we deuce from $\ell^{r_{n}-k}\widehat{a}_{2,2}+\widehat{a}_{2,3}\in\Z_\ell$ and $2\widehat{a}_{2,1}\widehat{a}_{2,2}+\widehat{a}_{2,3}^2=0$ that
			\[
			\mathrm{ord}_\ell(\widehat{a}_{2,2})=2r<0,
			\]
			which contradicts (\ref{relation among widehat{a}_{i,j}}). Therefore the second row of $\widehat{g}^{-1}$ also has entries in $\Z_\ell$. The same argument applies to the last row of it and one has $\widehat{g}^{-1}\in H(\Z_\ell)$.

			As a summary, for the case $n=3$, each element in $\mathcal{E}_{\mathbf{z},\xi;\ell}$ is of the form
			\[
			h_1u_1(-\ell^{-r_{n}})\mathrm{diag}(\ell^{-k},\ell^{k},1)h_2
			\]
			where $k\in\{0,1,\cdots,r_{n}\}$, $h_1=\mathrm{diag}(1,1,\pm1)$ and $h_2\in H(\Z_\ell)$.

			For general $n$, each $\mathcal{E}_{\mathbf{z},\xi;\ell}$ is of the form
			\[
			h_1\prod_{i=1}^{n}u_i(-\ell^{-r_i})\mathrm{diag}(1_{n-1},\ell^{-k},1_{n-1},\ell^{k},1)h_2
			\]
			where $k\in\{0,1,\cdots,r_{n}\}$, $h_1=\mathrm{diag}(1_{2n},\pm1)$ and $h_2\in H(\Z_\ell)$.		    
			
		\end{enumerate}

		Now we proceed to $(b)$. We treat the case $n=5$ and assume $r_1\ge k$, the general case is similar.
		Clearly for the places outside $\ell$,
		this is true
		since
		$\varsigma_{\underline{r},k}$
		has trivial components outside the place $\ell$.
		Write $\mathcal{T}t\mathcal{T}^{-1}=\mathrm{diag}(t_1,t_2,1/t_1,1/t_2,1)$, then one has
		\[
		\mathcal{T}
		\varsigma_{\underline{r},k}^{-1}t_\ell\varsigma_{\underline{r},k}
		\mathcal{T}^{-1}
		=
		\begin{pmatrix}
			t_1 & 0 & 0 & 0 &0 \\
			\ell^{-r_1+k}(t_1-t_2) & t_2 & 0 & \ell^{-2r_2+2k}(t_2+t_2^{-1}-2)/2 & \ell^{-r_2+k}(t_2-1) \\
			0 & 0 & t_1^{-1} & \ell^{-r_1+k}(t_1^{-1}-t_2^{-1}) & 0 \\
			0 & 0 & 0 & t_2^{-1} &0 \\
			0 & 0 & 0 & \ell^{-r_2+k}(1-t_2^{-1}) & 1
		\end{pmatrix}
		\]
		So if $t_2-1\in\ell^{r_2-k}\Z_\ell$ and $t_1t_2^{-1}-1\in\ell^{r_1-k}\Z_\ell$, then the above matrix lies in $H(\Z_\ell)$.

	\end{proof}

	\begin{lemma}
		We have the following
		\begin{equation}\label{toric period as a sum}			
			\frac{\mathbf{B}_{S,\psi}
				(g_{\xi(\underline{r})},\varepsilon)}{\gamma(\mathrm{det}(\xi(\underline{r})),q_U)
				\chi_{U}(\mathrm{det}\,\xi)
				|\mathrm{det}\,\xi|^{(2n+1)/2}_{\mathbb{A}_f}}
			=
			\frac{
				\varepsilon\mathbf{P}_{\mathbf{f}}(\psi,\varsigma_{\underline{r},0})
			}{
				[ H(\widehat{\Z}):K^{(\mathbb{S})}]
			}.
		\end{equation}        
	\end{lemma}    
	\begin{proof}
		We claim that the following natural map is in fact a bijection
		\[
		H(\Q)_{\mathbf{z}}\backslash\mathcal{E}_{\mathbf{z},\xi}/K^{(\mathbb{S})}
		\to
		H(\mathbb{A}_f)_{\mathbf{z}}\backslash\mathcal{E}_{\mathbf{z},\xi}/K^{(\mathbb{S})}.
		\]
		Indeed, for any $g\in\mathcal{E}_{\mathbf{z},\xi}$, we can write
		\[
		g=h_{\mathbf{z}}\varsigma_{\underline{r},k}h^{(\mathbb{S})}
		\]
		with $h_{\mathbf{z}}\in H(\mathbb{A}_f)_{\mathbf{z}}$, $k\in\{0,\cdots,r_{n}\}$ and $h^{(\mathbb{S})}\in K^{(\mathbb{S})}$. If the $\ell$-th component $(h_\mathbf{z})_\ell=1_{2n+1}$, then $g=\varsigma_{\underline{r},k}h_{\mathbf{z}}h^{(\mathbb{S})}$ and $h_\mathbf{z}h^{(\mathbb{S})}\in K^{(\mathbb{S})}$. If $(h_{\mathbf{z}})_\ell=\mathrm{diag}(1_{2n},-1)$, then take $h_{\Q}=\mathrm{diag}(1_{2n},-1)\in H(\Q)_{\mathbf{z}}$ and we have
		\[
		g=h_\Q(h_\Q h_{\mathbf{z}})\varsigma_{\underline{r},k}h^{(\mathbb{S})}=h_\Q\varsigma_{\underline{r},k}(h_\Q h_{\mathbf{z}}h^{(\mathbb{S})}).
		\]
		This proves the claim.

		Now by definition, the LHS in the lemma is equal to $[H(\widehat{\Z}):K^{(\mathbb{S})}]^{-1}
		\sum_{k=0}^{r_{n}}\mathbf{P}_{\mathbf{f}}(\psi,\varsigma_{\underline{r},k})$. It suffices to show that for $k>0$, we must have
		\[
		\mathbf{P}_{\mathbf{f}}(\psi,\varsigma_{\underline{r},k})
		=
		\int_{[T]_f}
		\langle
		\Delta_{\lambda}(\mathcal{T}\mathbf{z}),\mathbf{f}(t_f\varsigma_{\underline{r},k})
		\rangle_{W,U}
		\psi(t_f)dt_f=0.
		\]
		By Proposition \ref{E_z set of representatives}(b), for $t\in T(\widehat{\Z})_{\underline{r}(k)}$, we have $\mathbf{f}(t\varsigma_{\underline{r},k})=\mathbf{f}(\varsigma_{\underline{r},k})$. However, condition \ref{condition on the character psi} says $\psi_{n}$ is of conductor equal to $\ell^{r_{n}}$, thus as long as $k>0$, one has $\mathbf{P}_{\mathbf{f}}(\psi,\varsigma_{\underline{r},k})=0$.
	\end{proof}

	\section{Non-vanishing modulo $p$ of theta lift}\label{non-vanishing of theta lift mod p}
	In this section, we study the $p$-integrality of $\Theta_{ \phi_{\lambda,\mathbf{f}}}$ and assume the non-vanishing modulo $\mathfrak{P}$, we show the $p$-primitivity of $\Theta_{ \phi_{\lambda},\mathbf{f}}$ for $p$-primitive $\mathbf{f}$.
	For
	$\mathbf{f}
	\in
	\mathcal{A}_{\rho_{\lambda}}(H,K^{(\mathbb{S})})$,
	we have the Fourier expansion
	of the theta lift
	as in
	(\ref{Fourier expansion}):
	\[
	\Theta_{\phi_{\lambda},\mathbf{f}}^\ast
	=
	\sum_{S\in\mathrm{Sym}_{2n}^\circ(\Z)}
	a_{\Theta_{ \phi_{\lambda},\mathbf{f}}}(S)q^S.
	\]

	First consider the $p$-integrality of theta lifts. We consider the following condition:
	\begin{taggedtheorem}{(C3)}\label{condition p not dividing the level}
		$p\nmid[H(\widehat{\Z}):K^{(\mathbb{S})}]$.\index{C@(C3)}
	\end{taggedtheorem}

	As in \cite[Proposition 5.1]{HsiehNamikawa2017}, we have the following $p$-integrality result for theta lifts:
	\begin{theorem}\label{p-integrality of lift of f}
		Assume conditions \ref{condition on p} and \ref{condition p not dividing the level}.
		If
		$\mathbf{f}
		\in
		\mathcal{A}_{\rho_{\lambda}}
		( H,K^{(\mathbb{S})})$
		is
		$p$-integral,
		then
		$\Theta_{ \phi_{\lambda},\mathbf{f}}
		\in
		\mathcal{A}_{\rho_{\tau^\circ}}
		(\widetilde{G},\Gamma_0(2,N_{\mathbb{S}}),\chi_{U}^\circ)$
		is also
		$p$-integral.
	\end{theorem}
	\begin{proof}
		By
		Proposition
		\ref{Fourier coefficient of classicla Siegel modular form}
		and the definition of
		the $p$-adic avatar
		$\overline{\mathbf{f}_p}$
		of
		$\mathbf{f}$,
		we have
		(for some $z\in W^+(\Q)$):
		\begin{equation}\label{a(S), as a sum}
			a_{\Theta_{ \phi_{\lambda},\mathbf{f}}}(S_{z})
			=
			\frac{1}{[H(\widehat{\Z}):K^{(\mathbb{S})}]}
			\cdot
			\sum_{[h_f]\in[\mathcal{E}_{z,K^{(\mathbb{S})}}]}
			w_{z,h_f}
			\langle
			\Delta_{\lambda}(\mathcal{T}h_p^{-1}z),
			\overline{\mathbf{f}_p}(h_f)
			\rangle_{W,U}.
		\end{equation}
		By definition of $h_f\in\mathcal{E}_{z,K^{(\mathbb{S})}}$,
		$h_p^{-1}z
		\in
		W^+(\mathcal{O})$
		and thus
		$\Delta_{\lambda}(\mathcal{T}h_p^{-1}z)
		\in
		\M_{\lambda}(\mathcal{O})
		\otimes_{\mathcal{O}}
		\M_{\tau}(\mathcal{O})$.

		Write $z=(z_1,z_2,\cdots,z_{2n})$ in column matrices.
		By definition, if $\mathcal{E}_z\neq\emptyset$, in particular, there is $h\in H(\Q_v)$ for some $v\in\mathbb{S}$ such that $h^{-1}z\in\mathrm{Supp}_\ell(r_\ell)$, so the matrix $(\langle z_i,z_j\rangle_{U})_{i,j}$ is non-degenerate. Therefore the order of $H(\Q)_{z}$ divides $2$: indeed, since the matrix $(\langle z_i,z_j\rangle_{U})_{i,j}$ is non-degenerate, the column vectors $z_1,\cdots,z_{2n}$ in the $\Q$-vector space $\Q^{2n+1}=\mathrm{M}_{2n+1,1}(\Q)$ are linearly independent. Choose $z_{2n+1}\in\Q^{2n+1}$ such that $(z_1,\cdots,z_{2n+1})$ forms a basis of $\Q^{2n+1}$. Then the matrix of an element $g\in H(\Q)_z$ under this basis is given by the following upper-triangular matrix ($a_1,\cdots,a_{2n+1}\in\Q$)
		\[
		A_g=A(a_1,\cdots,a_{2n+1}):=
		\begin{pmatrix}
			1 & 0 & & a_1 \\
			& \ddots & & \vdots \\
			& 0 & 1 & a_{2n} \\
			0 & \cdots & 0 & a_{2n+1}
		\end{pmatrix}.
		\]
		Since $\det(g)=\pm1$, we have $a_{2n+1}=\pm1$.
		\begin{enumerate}
			\item 
			If $a_{2n+1}=1$, since $g\in H(\Q)_z\subset H(\R)$ and $H(\R)$ is compact, we must have $a_1,\cdots,a_{2n}=0$ and so $A_g=1$, that is, $g=1$.

			\item 
			If $a_{2n+1}=-1$, one checks that $A_g^2=1$. Moreover, if there is another element $g'\in H(\Q)_z$ different from $1$, then we can write it as $A_{g'}=A(a_1',\cdots,a_{2n}',-1)$. Moreover,
			\[
			A_{g'}A_g
			=
			A(a_1-a_1',\cdots,a_{2n}-a_{2n}',1)
			\in
			H(\Q)_z\subset H(\R).
			\]
			By the case (1), we have $a_1-a_1'=\cdots=a_{2n}-a_{2n}'=0$ and thus $A_g=A_{g'}$.
			
		\end{enumerate}
		We conclude from the above two cases that there is at most one element $g\in H(\Q)_z$ different from $1$. So $\# H(\Q)_z$ is equal to $1$ or $2$, which clearly divides $2$.
		
		Since $p>2$, $w_{z,h_f}^{-1}
		=
		\sharp
		(H(\Q)_{z}\cap h_fK^{(\mathbb{S})}h_f^{-1})$ is a $p$-unit.
		Therefore every factor in the above expression for $a_{\Theta_{ \phi_{\lambda},\mathbf{f}}}(S_z)$
		lies in $\mathcal{O}$ or
		$\M_{\tau}(\mathcal{O})$,
		so
		$a_{\Theta_{ \phi_{\lambda},\mathbf{f}}}(S_z)
		\in
		\M_{\tau}(\mathcal{O})$.
	\end{proof}

	Now consider the $p$-primitivity of theta lifts.
	Choose a finite set of representatives
	$\widetilde{\mathcal{E}}_{T}$
	for the double coset
	$[T]_f/
	T(\widehat{\Z})$
	such that each
	$g\in\widetilde{\mathcal{E}}_{T}$
	has component
	$g_p=1$. So $\widetilde{\mathcal{E}}_T$ is a finite set.
	We relate the non-vanishing modulo $\mathfrak{P}$ of
	$a_{\Theta_{ \phi_{\lambda},\mathbf{f}}}(S)$
	to that of
	the toric period
	$\mathbf{B}_{S,\psi}
	(g_{\xi})$.

	\begin{theorem}\label{p-primitivity of theta lift of f}
		Assume conditions \ref{condition on lambda}, \ref{condition on p} and \ref{condition p not dividing the level}. If $\mathbf{f}\in\mathcal{A}_{\rho_{\lambda}}(H,K^{(\mathbb{S})})$ is	$p$-integral, $\xi=\xi(\underline{r})$ and $\psi\colon[T]\rightarrow	\mathbb{C}^\times$ a continuous character of finite order satisfying condition \ref{condition on the character psi}. Then
		\[
		\mathbf{B}_{S_\mathbf{z},\psi}(g_{\xi},\varepsilon)
		\in
		\M_{\tau}(\mathcal{O}).
		\]
		Moreover if $\overline{\mathbf{B}}_{S_\mathbf{z},\psi}(g_{\xi},\varepsilon)\neq0$ (recall $\overline{\mathbf{B}}_{S_\mathbf{z},\psi}(g_{\xi},\varepsilon)$ is reduction modulo $\mathfrak{P}$ of $\mathbf{B}_{S_\mathbf{z},\psi}(g_{\xi},\varepsilon)$), then there is
		$S\in\mathrm{Sym}_{2n}^\circ(\Z)$ such that $a_{\Theta_{\phi_{\lambda},\mathbf{f}}}(S)\not\equiv0	(\mathrm{mod}\,\mathfrak{P})$. In particular, the theta lift	$\Theta_{\phi_{\lambda},\mathbf{f}}$ is $p$-primitive.
	\end{theorem}
	\begin{proof}
		The argument is similar to \cite[Lemma 5.2]{HsiehNamikawa2017}.
		Since $\mathcal{T}\mathbf{z}\in W^+(\mathcal{O})$, the evaluation of $z_{\lambda;i,j}$ at $\mathbf{z}$ is equal to $\mathbf{E}_j(\mathcal{T}\mathbf{z})(\mathbf{e}_i^\vee)$, which lies in $\mathcal{O}$ (here we use the interpretation that $\mathbf{E}_j$ is a $\M_{\lambda}(\C)$-valued pluri-harmonic polynomial). Thus $\mathbf{P}_{\mathbf{f}}(\psi,\varsigma_{\underline{r},0})\in\M_{\tau}(\mathcal{O})$ and
		the first part follows from
		(\ref{toric period as a sum}).
		Now if
		$\overline{\mathbf{B}}_{S_\mathbf{z},\psi}(g_{\xi},\varepsilon)
		\neq0$,
		then
		again by
		(\ref{toric period as a sum}),
		we have
		\begin{align*}
			&
			(\chi_{U}(\mathrm{det}\,\xi)
			|\mathrm{det}\,\xi|^{(2n+1)/2}_{\mathbb{A}_f})^{-1}
			[ H(\widehat{\Z}):K^{(\mathbb{S})}]
			\mathbf{B}_{S_\mathbf{z},\psi}
			(g_{\xi},\varepsilon)
			\\
			=
			&
			\varepsilon N(\underline{r})\gamma(\det(\xi),q_U)
			\sum_{t\in[T]_f/
				T(\widehat{\Z})_{\underline{r}}}
			\,
			\int\limits_{ H(\Q)_{\mathbf{z}}
				\backslash
				\mathcal{E}_{\mathbf{z},\xi}}
			\langle
			\Delta_{\lambda}(\mathcal{T}\mathbf{z}),
			\mathbf{f}(th_f)
			\rangle_{W,U}
			\psi(t)
			dh_f
			\\
			=
			&
			\varepsilon N(\underline{r})\gamma(\det(\xi),q_U)
			\mathbf{e}_\infty(-i\mathrm{Tr}S_\mathbf{z})
			\sum_{t\in[T]_f/
				T(\widehat{\Z})_{\underline{r}}}
			\mathbf{a}_{\Theta_{\phi_{\lambda},\mathbf{f}},
				S_\mathbf{z}}
			(\mathbf{j}(t)g_{\xi})
			\psi(t).
		\end{align*}
		Recall
		$\mathbf{j}(t)g_{\xi}
		=
		\mathrm{diag}
		(\mathbf{j}'(t)\xi,
		(\mathbf{j}'(t)\xi)^{-\mathrm{t}}
		)
		$.
		By strong approximation for
		$\mathrm{SL}_{2n}$,
		we can write
		\[
		\mathrm{GL}_{2n}(\mathbb{A}_f)
		\ni
		\mathbf{j}'(t)\xi=r_tu_t
		\text{ with }
		r_t\in\mathrm{GL}_{2n}(\Q)
		\text{ and }
		u_t\in\mathrm{GL}_{2n}(\widehat{\Z}).
		\]
		By our choice of $\ell$,
		we actually have
		$\mathrm{diag}(u_t,u_t^{-\mathrm{t}})
		\in
		\Gamma_0(2,N_{\mathbb{S}})$.
		By
		(\ref{Fourier coefficient, transformation property}), (\ref{a_{f,S} and a_f(S)})
		and the fact that
		$\Theta_{\phi_{\lambda},\mathbf{f}}$
		is of level $\Gamma_0(2,N_{\mathbb{S}})$,
		we have
		\begin{align}\label{W and a}
			\begin{split}
				\mathbf{a}_{\Theta_{\phi_{\lambda},\mathbf{f}},
					S_\mathbf{z}}
				(\mathbf{j}(t)g_{\xi})
				&
				=
				\mathbf{a}_{\Theta_{\phi_{\lambda},\mathbf{f}},
					r_t^\mathrm{t}S_\mathbf{z}r_t}
				(\mathrm{diag}(r_t^{-1},r_t^{\mathrm{t}})_\infty u_t)
				=
				\mathbf{a}_{\Theta_{\phi_{\lambda},\mathbf{f}},
					r_t^\mathrm{t}S_\mathbf{z}r_t}
				(\mathrm{diag}(r_t^{-1},r_t^{\mathrm{t}})_\infty)
				\\
				&
				=
				\det(\xi_\infty)^{(2n+1)/2}
				\mathbf{e}_\infty(i\mathrm{Tr}S_\mathbf{z})
				\rho_{\tau}
				(r_t^{-\mathrm{t}})
				a_{\Theta_{ \phi_{\lambda},\mathbf{f}}}
				(r_t^\mathrm{t}S_\mathbf{z}r_t)
				\\
				&
				=
				\mathbf{e}_\infty(i\mathrm{Tr}S_\mathbf{z})
				\rho_{\tau}
				(r_t^{-\mathrm{t}})
				a_{\Theta_{ \phi_{\lambda},\mathbf{f}}}
				(r_t^\mathrm{t}S_\mathbf{z}r_t).
			\end{split}
		\end{align}
		Here $\mathrm{diag}(r_t^{-1},r_t^{\mathrm{t}})_\infty$ is the image of $\mathrm{diag}(r_t^{-1},r_t^{\mathrm{t}})\in\mathrm{GL}_{2n}(\Q)$ in $\mathrm{GL}_{2n}(\mathbb{A})$ via the embedding $\mathrm{GL}_{2n}(\Q)\hookrightarrow\mathrm{GL}_{2n}(\R)\hookrightarrow\mathrm{GL}_{2n}(\mathbb{A})$. Thus we get		
		\begin{equation}
			\frac{[ H(\widehat{\Z}):K^{(\mathbb{S})}]
				\mathbf{B}_{S_\mathbf{z},\psi}
				(g_{\xi},\varepsilon)}{\varepsilon N(\underline{r})\gamma(\mathrm{det}(\xi),q_U)\chi_{U}(\mathrm{det}\,\xi)
				|\mathrm{det}\,\xi|^{(2n+1)/2}_{\mathbb{A}_f}}
			=
			\sum_{t\in[T]_f/
				T(\widehat{\Z})_{\underline{r}}}
			\rho_{\tau}(r_t^{-\mathrm{t}})
			a_{\Theta_{ \phi_{\lambda},\mathbf{f}}}
			(r_t^\mathrm{t}S_\mathbf{z}r_t)			
			\psi(t).
		\end{equation}
		It follows that if
		$\overline{\mathbf{B}}_{S_\mathbf{z},\psi}(g_{\xi},\varepsilon)
		\neq0$,
		then the LHS of the above identity is
		$\not\equiv0(\mathrm{mod}\,\mathfrak{P})$
		(we have $p\nmid  N(\underline{r})$
		by
		(\ref{p not dividing t_r})
		and
		$p\nmid[H(\widehat{\Z}):K^{(\mathbb{S})}]$
		by assumption).

		Let
		$\widetilde{\mathcal{E}}_{T,\xi}$
		be a set of representatives for
		the quotient
		$(
		T(\Q)T(\widehat{\Z})
		)/
		(
		T(\Q)T(\widehat{\Z})_{\underline{r}}
		)$
		consisting of elements whose $p$-th component is $1$.
		Then
		the product set
		$\widetilde{\mathcal{E}}_{T}
		\cdot
		\widetilde{\mathcal{E}}_{T,\xi}$
		is a set of representatives for
		$[T]_f/
		T(\widehat{\Z})_{\underline{r}}$.
		For all
		$t\in
		\widetilde{\mathcal{E}}_T\cdot\widetilde{\mathcal{E}}_{T,\xi}$,
		one has
		$t_p=1$,
		thus
		\[
		r_t\in\mathrm{SL}_{2n}(\Z_p)
		\cap
		\mathrm{SL}_{2n}(\Q).
		\]
		So
		$\rho_{\tau}(r_t^{-\mathrm{t}})
		a_{\Theta_{ \phi_{\lambda},\mathbf{f}}}
		(r_t^{\mathrm{t}}S_{\mathbf{z}}r_t)
		\in\M_{\tau}(\mathcal{O})$.
		On the other hand, for any
		$t\in\widetilde{\mathcal{E}}_T
		\cdot\widetilde{\mathcal{E}}_{T,\xi}$, one has
		\[
		a_{\Theta_{ \phi_{\lambda},\mathbf{f}}}(r_t^\mathrm{t}S_{\mathbf{z}}r_t)
		=
		a_{\Theta_{ \phi_{\lambda},\mathbf{f}}}(S_{\mathbf{z}r_t}).
		\]
		Now if
		$\overline{\mathbf{B}}_{S_\mathbf{z},\psi}
		(g_{\xi},\varepsilon)
		\neq0$,
		then for at least one
		$t$,
		we have
		$a_{\Theta_{ \phi_{\lambda},\mathbf{f}}}(S_{\mathbf{z}r_t})\not\equiv0
		(\mathrm{mod}\,\mathfrak{P})$,
		which gives the $p$-primitivity of
		$\Theta_{\phi_{\lambda},\mathbf{f}}$.
	\end{proof}

	\begin{remark}
		It remains to show $\overline{\mathbf{B}}_{S_\mathbf{z},\psi}
		(g_{\xi},\varepsilon)
		\neq0$ for some $\xi(\underline{r})$ and $\psi$. This will be done in §\ref{non-vanishing of Besse mod p}.
	\end{remark}

	\section{Non-vanishing modulo $p$ of toric periods}\label{non-vanishing of Besse mod p}		
	\subsection{Preliminaries}
	\label{Equidistribution of CM points}
	For $1\leq i_1<i_2<\cdots<i_k\leq2n+1$, we write
	\begin{align*}
		\delta_{i_1,\cdots,i_k}
		&
		=\delta_{i_1}\cdots\delta_{i_k},
		\\
		\widetilde{E}_{i_1,\cdots,i_k}
		&
		=\widetilde{E}_{i_1}\otimes\cdots\otimes\widetilde{E}_{i_k}.
	\end{align*}
	We recollect some results from
	\S \ref{Odd spin groups}. Apart from the condition \ref{condition on ell} on the prime $\ell$,
	we consider the following condition on the quadratic space $(U,\langle-,-\rangle)$ in this section:
	\begin{taggedtheorem}{(C4)}\label{condition-delta are squares in Q}
		The products
		$\delta_{1,2,3,4},\delta_{1,2,5,6},\cdots,
		\delta_{1,2,2n-1,2n}$
		are all perfect squares in
		$\Q$; moreover, for each prime $q$, there exist a pair of integers $(i,j)$ such that $\delta_i\delta_j$ is a $q$-adic unit. Here we require $1\le i<j\le2n$ and $(i,j)\neq(2k_0-1,2k_0)$ for any integer $k_0$.\index{C@(C4)}
	\end{taggedtheorem}
	We write
	$\mathbf{H}=\mathrm{GSpin}(U)$
	for the general spin group associated to the quadratic space
	$U$.
	Recall
	$\mathbf{H}( \Q_\ell)$
	is given by the units
	$C(U_\ell)^\times$
	of the even part $C(U_\ell)$
	of the Clifford algebra
	\[
	\mathrm{Cliff}(U_\ell)=\bigoplus_{k=0}^\infty U_\ell^{\otimes k}/
	\langle
	v\otimes v-\langle v,v\rangle_U\mid v\in U_\ell
	\rangle
	\]
	associated to the quadratic space
	$(U_\ell,\langle-,-\rangle_{U_\ell})$.
	We have a central extension of group schemes over $\Z$:
	\[
	1
	\rightarrow
	\mathbb{G}_m
	\rightarrow
	\mathbf{H}
	\xrightarrow{\mathfrak{s}}
	H
	\rightarrow
	1,
	\]
	which is induced from the conjugate action of $\mathbf{H}$ on $U$.
	We will identify $\mathbb{G}_m$
	with the center of
	$\mathbf{H}$.
	The map
	$\mathfrak{s}$
	sends $g\in \mathbf{H}$ to
	the element in $H$
	which acts on
	$U$
	by mapping
	$v$ to
	$gvg^{-1}$.

	For distinct basis vectors
	$E_{i_1},\cdots,E_{i_k}$,
	write
	$C_{i_1,\cdots,i_k}$
	for the even part of the Clifford algebra associated to the quadratic submodule
	$(\Z(E_{i_1},\cdots,E_{i_k}))$
	of $U$. Recall we have defined $d_i=1/\sqrt{2(-1)^{i-1}\delta_i}$.
	We then write ($1\leq r,j,k\leq2n+1$)
	\begin{align*}
		e_{r,j,k}^\pm
		&
		=
		\frac{1}{4}
		\left(
		d_rd_k\mathbf{i}^{(r+k)/2}E_{r,j}\pm d_jd_k\mathbf{i}^{(j+k+1)/2}E_{j,k}
		\right),
		\\
		\tau_{j,k}(a)
		&
		=
		\frac{1}{a}
		\left(
		1+\frac{a^2-1}{8}
		\left(
		d_j\mathbf{i}^{(j+1)/2}E_j+d_k\mathbf{i}^{k/2}E_k
		\right)
		\left(
		d_j\mathbf{i}^{(j+1)/2}E_j-d_k\mathbf{i}^{k/2}E_k
		\right)
		\right),
		\quad
		a\in\Q_\ell^\times.
	\end{align*}
	Here $\mathbf{i}$ is a square root of $-1$. We then put
	\[
	\mathbf{e}_k
	=
	\frac{1}{2}
	\left(
	1-\frac{E_{1,2,2k-1,2k}}{\sqrt{\delta_{1,2,2k-1,2k}}}
	\right),
	\]
	an idempotent central element in
	$C_{1,2,2k-1,2k}$.
	Then one can check that
	$C_{1,2,2k-1,2k}\mathbf{e}_k$ and
	$C_{1,2,2k-1,2k}(1-\mathbf{e}_k)$
	are both central simple algebras over $\Z$, free of rank $4$.
	Then we put
	\[
	\mathbf{H}_j
	=
	\begin{cases*}
		(C_{1,2,2j+1,2j+2}\mathbf{e}_{j+1})^\times,
		&
		$j=1,\cdots,n-1$;
		\\
		C_{1,2,n}^\times,
		&
		$j=n$.
	\end{cases*}
	\]
	For $w\in\Q_\ell$ and $i=1,\cdots,n$, we write
	\begin{align*}
		v_i^+(w)
		&
		=
		\begin{cases*}
			1-\frac{w}{4}\widetilde{E}_{1,i},
			&
			$i=1,\cdots,n-1$;
			\\
			1+we_{1,2,n}^+,
			&
			$i=n$.
		\end{cases*}
		\\
		v_i^-(w)
		&
		=
		\begin{cases*}
			1-\frac{w}{4}\widetilde{E}_{n+1,i},
			&
			$i=1,\cdots,n-1$;
			\\
			1+we_{1,2,n}^-,
			&
			$i=n$.
		\end{cases*}
	\end{align*}
	The unipotent subgroups
	$\mathbf{U}_1^+,\cdots,\mathbf{U}_{n}^+,\mathbf{U}_1^-,\cdots,\mathbf{U}_{n}^-$
	in
	\S \ref{Odd spin groups}
	are given by
	(note that the enumerations are different, the $\mathbf{U}_1^\pm,\cdots,\mathbf{U}_{n}^\pm$
	in \S \ref{Odd spin groups}
	correspond to our $\mathbf{U}_{n}^\pm,
	\mathbf{U}_1^\pm,\cdots,\mathbf{U}_{n-1}^\pm$
	below
	respectively):
	\[
	\mathbf{U}_i^\pm=v_i^\pm(\Q_\ell).
	\]

	We also have torus subgroups ($i=1,2,\cdots,n$)
	\[
	\mathbf{T}_i=C_{2i-1,2i}^\times
	\subset \mathbf{H}.
	\]	
	These tori $\mathbf{T}_i$
	clearly commute with each other and each of these
	$\mathbf{T}_i$
	contains the center
	$\mathbb{G}_m\subset
	\mathbf{H}$. We put
	\[
	\mathbf{T}=\prod_{i=1}^{n}
	\mathbf{T}_i.\index{T@$\mathbf{T},\mathbf{T}_i$}
	\]
	
	Moreover the above central extension induces the following exact sequence
	of algebraic groups over $\Z$:
	\begin{equation*}
		1
		\rightarrow
		\mathbb{G}_m
		\rightarrow
		\mathbf{T}_i
		\xrightarrow{\mathfrak{s}}
		T_i
		\rightarrow
		1.
	\end{equation*}
	\begin{lemma}\label{surjectivity of s for torus}
		For any field extension
		$F/\Q$,
		we have a short exact sequence
		\[
		1
		\rightarrow
		\mathbb{G}_m(F)
		\rightarrow
		\mathbf{T}_i(F)
		\xrightarrow{\mathfrak{s}}
		T_i(F)
		\rightarrow
		1.
		\]
	\end{lemma}
	\begin{proof}
		It suffices to show the surjectivity of $\mathfrak{s}$. We first make $\mathfrak{s}$ explicit.
		For ease of notations,
		we write
		$\delta=\delta_{2i-1}/\delta_{2i}$
		and
		$\mathfrak{e}_1=E_{2i-1}$,
		$\mathfrak{e}_2=E_{2i}$.
		The group
		$T_i(F)$
		is isomorphic to a subgroup of
		$\mathrm{SL}_2(F)$
		consisting of matrices
		$g=\begin{pmatrix}
			a & b\delta \\
			-b & a
		\end{pmatrix}
		\in
		\mathrm{SL}_2(F)$
		with $a,b\in R$.
		In this case
		$\mathbf{T}_i(F)$
		consists of
		$x+y\mathfrak{e}_1\mathfrak{e}_2
		\in
		C(F\mathfrak{e}_1+F\mathfrak{e}_2)^\times$
		with
		$x,y\in F$.
		Moreover we have the following conjugate action
		\begin{align*}
			(x+y\mathfrak{e}_1\mathfrak{e}_2)
			\mathfrak{e}_1
			(x+y\mathfrak{e}_1\mathfrak{e}_2)^{-1}
			&
			=
			\frac{1}{x^2+y^2\delta_{2i-1,2i}}
			(x^2-y^2\delta_{2i-1,2i})\mathfrak{e}_1-2xy\delta_{2i-1}\mathfrak{e}_2,
			\\
			(x+y\mathfrak{e}_1\mathfrak{e}_2)
			\mathfrak{e}_2
			(x+y\mathfrak{e}_1\mathfrak{e}_2)^{-1}
			&
			=
			\frac{1}{x^2+y^2\delta_{2i-1,2i}}
			2xy\delta_{2i}\mathfrak{e}_1+(x^2-y^2\delta_{2i-1,2i})\mathfrak{e}_2.
		\end{align*}
		We put these two identities in a compact form
		\[
		(x+y\mathfrak{e}_1\mathfrak{e}_2)
		(\mathfrak{e}_1,\mathfrak{e}_2)
		(x+y\mathfrak{e}_1\mathfrak{e}_2)^{-1}=
		(\mathfrak{e}_1,\mathfrak{e}_2)
		\frac{1}{x^2+y^2\delta_{2i-1,2i}}
		\begin{pmatrix}
			x^2-y^2\delta_{2i-1,2i} & 2xy\delta_{2i} \\
			-2xy\delta_{2i-1} & x^2-y^2\delta_{2i-1,2i}
		\end{pmatrix}.
		\]
		Thus the map
		$\mathfrak{s}$
		is given by
		\[
		\mathfrak{s}
		\colon
		\mathbf{T}_i(F)
		\rightarrow
		T_i(F),
		\quad
		x+y\mathfrak{e}_1\mathfrak{e}_2
		\mapsto
		\frac{1}{x^2+y^2\delta_{2i-1,2i}}
		\begin{pmatrix}
			x^2-y^2\delta_{2i-1,2i} & -2xy\delta_{2i-1} \\
			2xy\delta_{2i} & x^2-y^2\delta_{2i-1,2i}
		\end{pmatrix}.
		\]
		For any
		$a,b\in F$
		such that
		$a^2+b^2\delta=1$,
		we can always find
		$x,y\in F$
		such that
		$a=(x^2-y^2\delta_{2i-1,2i})/(x^2+y^2\delta_{2i-1,2i})$
		and
		$b=-2xy\delta_{2i}$:
		if
		$a=1$,
		we take
		$x=1,y=0$;
		if
		$a\ne1$,
		we take any
		$y\ne0$ and set
		$x=y\delta b/(a-1)$. This shows that $\mathfrak{s}$ is surjective.
	\end{proof}
	A similar computation as in the above lemma shows that
	the image of these
	$\mathbf{T}_i( \Q_\ell)$ and $\mathbf{U}_i^+$
	under the map
	$\mathfrak{s}$
	is given as follows:
	\begin{align*}
		\mathfrak{s}(\mathbf{T}_i( \Q_\ell))
		&
		=		
		\mathcal{T}^{-1}
		\{
		\mathrm{diag}(1_{i-1},t,1_{n-i},1_{i-1},t^{-1},1_{n-i})
		\mid
		t\in
		\Q_\ell^{\times}
		\}
		\mathcal{T},
		\\
		\mathfrak{s}
		(v_1^+(w_1)
		\cdots
		v_{n}^+
		(w_n))
		&
		=
		u_1(w_1)\cdots u_{n}(w_n);
		\quad
		\forall
		w_1,\cdots,w_n\in\Q_\ell.
	\end{align*}
	In particular,
	\begin{equation}\label{unipotent in appendix related to main body of article}
		\mathfrak{s}
		(
		v_1^+(\ell^{-r_1})
		\cdots
		v_{n}^+(\ell^{-r_{n}})
		)
		=
		\varsigma_{(r_1,\cdots,r_{n}),0}.
	\end{equation}
	Here $\varsigma_{(r_1,\cdots,r_{n}),0}$ is given as in (\ref{varsigma}).
	Therefore, the unipotent subgroups
	$\mathbf{U}_1^+,\cdots,\mathbf{U}^+_{n}$
	are the root subgroups of
	$\mathfrak{s}(\mathbf{T}( \Q_\ell))$
	and their opposite root subgroups are given by
	$\mathbf{U}_1^-,\cdots,\mathbf{U}^-_{n}$.
	By \S \ref{Odd spin groups}, the torus
	$\mathbf{T}( \Q_\ell)$
	acts by conjugation on these unipotent subgroups
	$\mathbf{U}_i^\pm$
	by the following formula ($j=1,n+1$):
	\begin{align}\label{conjugate action of torus on unipotent group}
		\begin{split}
			\mathrm{Ad}_{\tau_{2i-1,2i}(a)}(1+v\widetilde{E}_j\widetilde{E}_{n+2})
			&
			=
			\begin{cases*}
				(1+a^{\pm2}v\widetilde{E}_j\widetilde{E}_{n+2}),
				&
				$i=1,2$;
				\\
				(1+v\widetilde{E}_j\widetilde{E}_{n+2}),
				&
				$i=3,\cdots,n$;
			\end{cases*}
			\\
			\mathrm{Ad}_{\tau_{1,2}(a)}(1+v\widetilde{E}_j\widetilde{E}_k)
			&
			=
			(1+a^{\pm2}v\widetilde{E}_j\widetilde{E}_k),
			\quad
			k=2,\cdots,n;
			\\
			\mathrm{Ad}_{\tau_{2i-1,2i}(a)}(1+v\widetilde{E}_j\widetilde{E}_k)
			&
			=
			\begin{cases*}
				(1+a^{\mp2}v\widetilde{E}_j\widetilde{E}_k),
				&
				$i=k\in\{2,\cdots,n\}$;
				\\
				(1+v\widetilde{E}_j\widetilde{E}_k),
				&
				$i\ne k\in\{2,\cdots,n\}$;
			\end{cases*}
			\\
			\mathrm{Ad}_{\tau_{2i-1,2i}(a)}(1+ve_{1,2,n}^\pm)
			&
			=
			\begin{cases*}
				1+a^{\pm2}ve_{1,2,n}^\pm,
				&
				$i=n$;
				\\
				1+ve_{1,2,n}^\pm,
				&
				$i\in\{1,2,\cdots,n-1\}$.
			\end{cases*}
		\end{split}
	\end{align}
	For $i=1,\cdots,n$,
	we define characters
	\[
	\chi_i
	\colon
	\mathbf{T}( \Q_\ell)
	\cap
	\mathbf{H}^1( \Q_\ell)
	\to
	\C^\times,
	\quad
	\prod_{i=1}^{n}\mathbf{T}_i( \Q_\ell)\ni
	(t_1,\cdots,t_{n})
	\mapsto
	\mu_i(\mathfrak{s}(t_i)).
	\]
	Then (\ref{conjugate action of torus on unipotent group}) shows that
	the action of
	$\mathbf{T}( \Q_\ell)$
	on
	$\mathbf{U}_i^\pm$
	correspond to the character
	$(\chi_i\chi_{n}^{-1})^{\pm1}$
	if
	$i=1,\cdots,n-1$
	and to the character
	$\chi_{n}^{\pm1}$
	if
	$i=n$.
	Therefore
	these
	$\mathbf{U}_1^+,\cdots,\mathbf{U}_{n}^+$
	correspond to a root basis for
	$\mathbf{T}( \Q_\ell)
	\cap
	\mathbf{H}^1( \Q_\ell)$
	(similarly for
	$\mathbf{U}_1^-,\cdots,\mathbf{U}_{n}^-$).

	\begin{lemma}\label{a variant of a theorem of Ihara}
		Fix two positive integers
		$k_1<k_2$
		and define a subgroup of $\mathrm{SL}_2(\Z_\ell)$:
		\[
		\mathfrak{G}
		=
		\left\{
		\begin{pmatrix}
			a & b \\
			c & d
		\end{pmatrix}
		\in
		\mathrm{SL}_2(\Z_\ell)
		|
		a,b,d\in\Z_\ell
		\text{ and }
		a-1,d-1,c\in\ell^{k_1}\Z_\ell
		\right\}.
		\]
		Then
		$\mathfrak{G}$
		and the matrix
		$\begin{pmatrix}
			1 & \ell^{-k_2} \\
			0 & 1
		\end{pmatrix}$
		generate the whole group
		$\mathrm{SL}_2( \Q_\ell)$.	    
	\end{lemma}
	\begin{proof}
		Write
		\[
		A=\begin{pmatrix}
			\ell^{k_1} & 0 \\ 0 & 1
		\end{pmatrix},
		\quad
		B(k)=\begin{pmatrix}
			1 & \ell^{-k} \\ 0 & 1
		\end{pmatrix}.
		\]
		Since $k_2>k_1$, it suffices to show
		$A\mathfrak{G}A^{-1}$
		and
		$AB(k_1+1)A^{-1}
		=B(1)=B(k_2-k_1)^{\ell^{k_2-k_1-1}}$
		generate $\mathrm{SL}_2( \Q_\ell)$.

		It is well known that $\mathrm{SL}_2(\Z_{\ell})$ and $B(1)$ generate $\mathrm{SL}_2(\Q_\ell)$. So it suffices to show that
		\[
		\mathrm{SL}_2(\Z_\ell)
		=
		\left\langle
		A\mathfrak{G}A^{-1},
		B(1)^\ell
		\right\rangle.
		\]
		We write the projection map
		\[
		\pi\colon\mathrm{SL}_2(\Z_\ell)\to\mathrm{SL}_2(\Z/\ell^{k_1}\Z).
		\]
		Note that $A\mathfrak{G}A^{-1}$ contains $\mathrm{Ker}(\pi)$. So to prove the lemma, it suffices to show that $\pi(A\mathfrak{G}A^{-1})$ and $\pi(B(1)^\ell)$ generate $\mathrm{SL}_2(\Z/\ell^{k_1}\Z)$. For this, note that $\pi(A\mathfrak{G}A^{-1})$ contains all matrices of the form
		$\begin{pmatrix}
			1 & 0 \\ \ast & 1
		\end{pmatrix}$ ($\ast\in\Z/\ell^{k_1}\Z$), while the group generated by $\pi(B(1)^{\ell})$ contains all matrices of the form $\begin{pmatrix}
			1 & \ast \\ 0 & 1
		\end{pmatrix}$ ($\ast\in\Z/\ell^{k_1}\Z$). Using the surjectivity of the reduction-mod-$\ell^{k_1}$ map $\mathrm{SL}_2(\Z)\to\mathrm{SL}_2(\Z/\ell^{k_1}\Z)$, one deduces that
		\[
		\mathrm{SL}_2(\Z/\ell^{k_1}\Z)
		=
		\left\langle
		\begin{pmatrix}
			1 & 0 \\ \ast & 1
		\end{pmatrix},
		\begin{pmatrix}
			1 & \ast \\ 0 & 1
		\end{pmatrix}
		\right\rangle.
		\]
		This finishes the proof of the lemma.

	\end{proof}

	\begin{lemma}\label{SL_2 and K generate Spin}
		Let $K$ be a compact open subgroup of $\mathbf{H}^1( \Q_\ell)$.
		For any $i=1,\cdots,n$,
		the two subgroups
		$\mathbf{H}_i^1(\Q_\ell)$
		and $K$
		generate
		$\mathbf{H}^1( \Q_\ell)$.
	\end{lemma}
	\begin{proof}
		We treat the case $i=1$, the other cases are similar. By (\ref{conjugate action of torus on unipotent group}),	for any positive integer $N>0$, $j=1,n+1$ and $k=2,\cdots,n$,	we have the following identity
		\[
		\left\{
		\mathrm{Ad}_{\tau_{1,2}(a)}(1+v\widetilde{E}_j\widetilde{E}_k)
		\mid
		a\in \Q_\ell^\times,
		v\in
		\ell^N\Z_\ell
		\right\}
		=
		1+ \Q_\ell\widetilde{E}_j\widetilde{E}_k.
		\]
		We choose $N\gg0$ such that $K$ contains $v_i^+(\ell^N\Z_{\ell})$ and $v_i^-(\ell^N\Z_{\ell})$ for all $i=2,\cdots,n$.
		Thus $\mathbf{T}( \Q_\ell)\cap\mathbf{H}^1( \Q_\ell)$ and $K$ generate a subgroup of	$\mathbf{H}^1( \Q_\ell)$ which contains	all the basic root subgroups $\mathbf{U}_2^\pm,\cdots,\mathbf{U}_{n}^\pm$ of $\mathbf{H}^1( \Q_\ell)$ (thus contains also the subgroups $\mathbf{H}_2^1(\Q_{\ell}),\cdots,\mathbf{H}_{n}^1(\Q_{\ell})$).
		Since $\mathbf{H}^1( \Q_\ell)$ is semi-simple and simply-connected,
		these subgroups $\mathbf{H}_1^1(\Q_{\ell}),\cdots,\mathbf{H}_n^1(\Q_{\ell})$
		generate $\mathbf{H}^1( \Q_\ell)$ and we are done.
	\end{proof}

	\subsection{Proof of non-vanishing modulo $p$ of toric periods}
	\label{non-vanishing of Besse mod p-subsection}	
	We assume condition \ref{condition-delta are squares in Q}.	We will show that
	there is $\xi=\xi(\underline{r})$ and a continuous character of finite order
	$\psi\colon[T]\rightarrow\mathbb{C}^\times$	satisfying condition \ref{condition on the character psi} such that $\overline{\mathbf{B}}_{S_\mathbf{z},\psi}
	(g_{\xi},\varepsilon)\neq0$. This is the most technical part of the article
	and we recommend the reader to start with the first paragraph in the proof of
	Theorem \ref{theorem on the non-vanishing mod p of toric periods}
	to get an idea of the plan of this section.

	By Lemma \ref{surjectivity of s for torus}, the map $\mathfrak{s}$ induces the following isomorphisms:
	\[
	\mathbf{T}( \Q_\ell)
	/
	\mathbf{T}(\Q)\mathbb{G}_m( \Q_\ell)
	\simeq
	T( \Q_\ell)
	/
	T(\Q)
	\simeq
	\mathcal{G}(\infty)
	/
	\mathcal{G}(\infty)'.
	\]
	Moreover, we have the following surjective maps
	\begin{equation}\label{T(A_f)-->G_0}
		\mathbf{T}(\mathbb{A}_f)
		\twoheadrightarrow
		\mathbf{T}( \Q_\ell)
		\twoheadrightarrow
		\mathbf{T}( \Q_\ell)
		/\mathbf{T}(\Q)\mathbb{G}_m( \Q_\ell)
		\xrightarrow{\sim}
		\mathcal{G}(\infty)/\mathcal{G}(\infty)'
		\twoheadrightarrow
		\mathcal{G}_0.
	\end{equation}

	We fix a finite subset
	\[
	\mathcal{R}=\{g_1,g_2,\cdots,g_{\sharp(\mathcal{G}_0)}\}
	\subset\mathbf{T}(\mathbb{A}_f)
	\]
	which is mapped bijectively to $\mathcal{G}_0$ by the map (\ref{T(A_f)-->G_0}).

	\begin{lemma}\label{R satisfies condition on g_k,g_l}
		For any two elements $g_i\neq g_k$ in $\mathcal{R}$, their $\ell$-th components satisfy
		\[
		(g_k)_\ell(g_i)_\ell^{-1}
		\notin
		\mathbf{T}(\Q_\ell)
		\cap
		\left(		
		\cup_{j=1}^{n}
		\mathbf{H}_j(\Q)
		Z(\mathbf{H}_j(\Q_\ell))
		\right).
		\]
	\end{lemma}
	\begin{proof}
		It suffices to note that
		$Z(\mathbf{H}_i(\Q_\ell))=Z(\mathbf{H}(\Q_\ell))
		\subset
		\mathbf{T}(\Q_\ell)$
		for any $i=1,\cdots,n$.
	\end{proof}
	For a ring $A$,
	write
	\[
	\widetilde{H}^1(A)
	:=
	\mathrm{Im}
	\left(
	\mathbf{H}^1(A)
	\hookrightarrow
	\mathbf{H}(A)
	\xrightarrow{\mathfrak{s}}
	H(A)
	\right),
	\]
	then the quotient $H(A)/\widetilde{H}^1(A)$ is an abelian group.
	We define the following map
	associated to $\mathbf{f}$ (recall $\mathbf{f}_p$ is the $p$-adic avatar of $\mathbf{f}$, see Definition \ref{automorphic forms on SO}):
	\begin{equation*}
		F_{\mathbf{f}}
		\colon
		H(\mathbb{A}_f)
		\rightarrow
		\M_{\tau}(\mathbb{C}_p),
		\quad
		g
		\mapsto
		\left\langle
		\Delta_{\lambda}(\mathcal{T}\mathbf{z}),
		\mathbf{f}_p
		(g)
		\right\rangle_{W,U}.
	\end{equation*}
	It follows that if
	$\mathbf{f}$ is $p$-integral,
	then
	$F_{\mathbf{f}}(g)\in\M_{\tau}(\mathcal{O})$
	for all $g\in H(\mathbb{A}_f)$.
	Under conditions \ref{condition on lambda} and \ref{condition on p}, $\M_{\lambda}(\kappa)$ is an irreducible representation of $H(\kappa)$ and it is of dimension $>1$, in particular, this representation does not factor through any abelian quotient of $H(\kappa)$. Thus if the image of $\overline{\mathbf{f}_p}$ generates the whole $\kappa$-vector space $\M_{\lambda}(\kappa)$, then $\overline{F}_{\mathbf{f}}$ is not $\widetilde{H}^1(\mathbb{A}_f)$-invariant (under right translation).

	\begin{theorem}\label{theorem on the non-vanishing mod p of toric periods}
		Assume condition \ref{condition on lambda}, \ref{condition on p}, \ref{condition p not dividing the level} and \ref{condition-delta are squares in Q}.
		Let
		$\mathbf{f}
		\in
		\mathcal{A}_{\rho_{\lambda}}
		( H,K^{(\mathbb{S})})$
		be a
		\emph{$p$-primitive} automorphic form.
		Then for any $r_1,\cdots,r_{n}
		\gg0$,
		there is a character
		$\psi
		\colon
		[T]
		\rightarrow
		\mathbb{C}^\times$
		satisfying condition
		\ref{condition on the character psi}
		such that
		\[
		\overline{\mathbf{B}}_{S_{\mathbf{z}},\psi}
		\neq0
		\]
	\end{theorem}
	\begin{proof}
		Since the whole proof is a bit technical,
		we give a summary of it:
		first we
		rewrite
		the toric period
		$\mathbf{P}_{\mathbf{f}}(\psi,\varsigma_{\xi})$
		in terms of a sum over
		$\mathcal{G}_0$
		of another function
		$\overline{F}_{\mathbf{f},\underline{r}',\underline{\alpha}'}$
		(this is (\ref{p=f})).
		Then
		we show that
		the latter is not $\widetilde{H}^1(\mathbb{A}_f)$-invariant provided that
		$\overline{F}_{\mathbf{f}}$
		is not $\widetilde{H}^1(\mathbb{A}_f)$-invariant,
		this is proved in
		Lemma \ref{widetilde(Q)_psi'' is not Spin-invariant},
		the main idea of the proof of this lemma is an argument by contradiction:
		we can show
		$\overline{F}_{\mathbf{f}}$
		is invariant under right translation of the subgroup $\mathbf{H}_n^1(\Q_{\ell})$ of
		$H(\Q_\ell)$.
		On the other hand,
		we know already
		this function is invariant under right translation of a certain compact open subgroup of
		$H(\Q_\ell)$
		and therefore one applies
		Lemma \ref{SL_2 and K generate Spin}
		to conclude that
		$\overline{F}_{\mathbf{f}}$
		is $\widetilde{H}^1(\mathbb{A}_f)$-invariant,
		a contradiction.
		Then we can apply
		Theorem
		\ref{automorphic application of appendix}
		and conclude that
		the above mentioned toric integral is non-zero modulo $\mathfrak{P}$
		(this is the last paragraph of this proof).

		Now we begin the proof of the theorem.
		Note that by (\ref{Q_Uz, computation}), $\Delta_\lambda(\mathbf{z})\neq0(\mathrm{mod}\,\mathfrak{P})$, thus by (\ref{toric period as a sum}), it suffices to show there exists $\psi$ satisfying condition \ref{condition on the character psi} such that
		\[
		\overline{\mathbf{P}_{\mathbf{f}}}(\psi,\varsigma_{\underline{r},0})
		\neq0.
		\]
		Recall that we have fixed an isomorphism
		\[
		\mathcal{G}(\infty)
		\simeq
		\mathcal{G}_0
		\times
		\Gamma(\infty),
		\]
		which induces the following two characters associated to $\psi$:
		\begin{align*}
			\psi'
			&
			\colon
			\mathcal{G}_0
			\hookrightarrow
			\mathcal{G}(\infty)
			\xrightarrow{\psi}
			\mathbb{C}^\times,
			\\
			\psi''
			&
			\colon
			\Gamma(\infty)
			\hookrightarrow
			\mathcal{G}(\infty)
			\xrightarrow{\psi}
			\mathbb{C}^\times.
		\end{align*}
		For $t_1,\cdots,t_{n}
		\in
		\Q_\ell$, we write
		\[
		u(t_1,\cdots,t_{n})
		=
		u_1(t_1)\cdots u_{n}(t_{n}).
		\]
		By
		(\ref{splitting of G(infty)})
		and
		(\ref{adjoint action of T on u, another basis}),
		one has the following
		\begin{align*}
			&
			\mathbf{P}_{\mathbf{f}}(\psi,\varsigma_{\underline{r},0})
			=
			\sum_{t\in\mathcal{G}(\infty)/T(\Z_\ell)_{\underline{r}}}
			F_{\mathbf{f}}(t\varsigma_{\underline{r},0})
			\psi(t)
			\\
			=
			&
			\sum_{\substack{t'\in\mathcal{G}_0 \\
					t''\in\prod_i^{n}\Gamma_i(r_i)}}
			F_{\mathbf{f}}
			(t't''u(\ell^{-r_1},\cdots,\ell^{-r_{n}})(t'')^{-1}t'')
			\psi'(t')
			\psi''(t'')
			\\
			=
			&
			\sum_{\substack{t'\in\mathcal{G}_0 \\
					t''=(t_1,\cdots,t_{n})\in\prod_i^{n}\Gamma_i(r_i)}}
			F_{\mathbf{f}}
			(t'u(\ell^{-r_1}\widetilde{\mu}_1(t_1),
			\cdots,\ell^{-r_{n}}\widetilde{\mu}_{n}(t_{n})))
			\psi'(t')
			\psi''(t'')
			\\
			=
			&
			\sum_{\substack{t'\in\mathcal{G}_0 \\
					t''=(t_1,\cdots,t_{n})\in\prod_i^{n}\Gamma_i(r_i)}}
			F_{\mathbf{f}}
			\left(
			t'
			u(\ell^{-r_1},\cdots,\ell^{-r_{n}})
			\cdot
			u(\ell^{-r_1}(\widetilde{\mu}_1(t_1)-1),
			\cdots,
			\ell^{-r_{n}}
			(\widetilde{\mu}_{n}(t_{n}-1)
			))
			\right)
			\psi'(t')
			\psi''(t'').
		\end{align*}
		For any $x\in\Z_\ell$,
		we define the following map
		\[
		(-)^x
		\colon
		1+\ell\Z_\ell
		\rightarrow
		1+\ell\Z_\ell,
		\quad
		1+y\mapsto (1+y)^x=\sum_{k=0}^\infty\binom{x}{k}y^k.
		\]
		Fix primitive characters
		$\psi_{i,r_i}
		\colon
		\Gamma_i(r_i)
		\rightarrow
		\mathbb{C}^\times$
		($i=1,\cdots,n$). For each $\underline{\alpha}=(\alpha_1,\cdots,\alpha_{n})\in\prod_{i=1}^n(\Z/\ell^{r_i}\Z)^\times$,
		we write
		$\psi_{i,r_i}^{\alpha_i}(-)=\psi_{i,r_i}((-)^{\alpha_i})=(\psi_{i,r_i}(-))^{\alpha_i}$ and define a primitive character
		\[
		\psi_{\underline{r},\underline{\alpha}}
		=
		\prod_{i=1}^{n}
		\psi_{i,r_i}^{\alpha_i}
		\colon
		\prod_{i=1}^n\Gamma_i(r_i)\to\C^\times.
		\]
		Conversely, any primitive character of $\prod_{i=1}^n\Gamma_i(r_i)$ is of the form $\psi_{\underline{r},\underline{\alpha}}$ for some $\underline{\alpha}\in\prod_{i=1}^n(\Z/\ell^{r_i}\Z)^\times$.
		
		We fix an $n$-tuple of positive integers
		$\underline{r}'=(r_1',\cdots,r'_{n})$
		such that
		$0
		<
		2r'_i
		<r_i$
		for all $i$ and an element
		$\underline{\alpha}'
		\in
		\prod_{i=1}^{n}
		(\Z_\ell/\ell^{r_i})^\times
		$,
		we then put
		\begin{align*}
			I(r'_i,r_i,\alpha_i')
			&
			=
			\alpha'_i+\ell^{r_i'}\Z_\ell/\ell^{r_i}\Z_\ell
			\\
			I(\underline{r}',\underline{r},\underline{\alpha'})
			&
			=
			\prod_{i=1}^{n}
			I(r_i',r_i,\alpha_i).
		\end{align*}

		We define the following map associated to $\overline{F}_{\mathbf{f}}$:
		\begin{align}\label{widetilde{F}}
			\begin{split}
				\overline{F}_{\mathbf{f},\underline{r}',\underline{\alpha}'}
				\colon
				&
				H(\mathbb{A}_f)
				\rightarrow
				\M_{\tau}(\kappa)
				\\
				g
				&
				\mapsto
				\ell^{-\sum_{i=1}^{n}(r_i-r_i')}
				\sum_{\substack{\underline{\alpha}\in I(\underline{r}',\underline{r},\underline{\alpha}') \\
						t''=(t_1,\cdots,t_{n})\in\prod_i^{n}\Gamma_i(r_i)}}
				\overline{F}_{\mathbf{f}}
				\left(
				g
				u
				\left(
				\ell^{-r_1}(\widetilde{\mu}_1(t_1)-1),
				\cdots,
				\ell^{-r_{n}}
				(\widetilde{\mu}_{n}(t_{n}-1)
				)
				\right)
				\right)
				\psi_{\underline{r},\underline{\alpha}}(t'').
			\end{split}			
		\end{align}

		For two integers
		$0<r'<r$,
		we fix an $\ell^r$-th root of unity $\zeta_{\ell^r}$,
		$\alpha',t\in\mathbb{Z}/\ell^e\mathbb{Z}$,
		$t_i\in\Gamma_i(r_i)$
		as in the above summation,
		we have the following identity:
		\begin{align*}
			\sum_{\alpha\in\alpha'+\ell^{r'}\mathbb{Z}/\ell^r\mathbb{Z}}
			\zeta_{\ell^r}^{t\alpha}
			&
			=
			\begin{cases*}	    	
				\ell^{r'}\cdot\zeta_{\ell^r}^{t\alpha'},
				&
				$\ell^{r-r'}\mid t$;
				\\
				0,
				&
				otherwise.
			\end{cases*},
			\\
			\sum_{
				\alpha_i\in I(r_i',r_i,\alpha_i')}
			\psi_{i,r_i}^{\alpha_i}(t_i)
			&
			=
			\begin{cases*}
				\ell^{r_i-r_i'}\cdot\psi_{i,r_i}^{\alpha_i'}(t_i),
				&
				$\ell^{r_i-r_i'}\mid(t_i-1)$;
				\\
				0,
				&
				otherwise.
			\end{cases*}
		\end{align*}
		We deduce easily
		\[
		\overline{F}_{\mathbf{f},\underline{r}',\underline{\alpha}'}
		(g)
		=		
		\sum_{t''=(t_1,\cdots,t_{n})\in\prod_{i=1}^{n}\Gamma_i(r_i)^{\ell^{r_i-r'_i}}}
		\overline{F}_{\mathbf{f}}
		\left(
		g
		u
		\left(
		\ell^{-r_1}(\widetilde{\mu}_1(t_1)-1),
		\cdots,
		\ell^{-r_{n}}
		(\widetilde{\mu}_{n}(t_{n}-1)
		)
		\right)
		\right)
		\psi_{\underline{r},\underline{\alpha}'}(t'').
		\]
		Note that $\overline{F}_{\mathbf{f},\underline{r}',\underline{\alpha}'}(g)$
		is independent of $\underline{r}$.
		Moreover the dependence of
		$\overline{F}_{\mathbf{f},\underline{r}',\underline{\alpha}'}$
		on
		$\underline{\alpha}'$
		factors through its image under the map
		$\prod_i
		(\Z_\ell/\ell^{r_i})^\times
		\rightarrow
		\prod_i
		(\Z_\ell/\ell^{r_i'})^\times$
		(we will identify $\underline{\alpha}'$ with its image).
		We denote
		\[
		\Gamma_i'(r_i')
		:=(1+\ell^{r_i-r_i'}\Z_\ell)/(1+\ell^{r_i}\Z_\ell).
		\]
		The bijection
		\[
		(1+\ell\Z_\ell)/(1+\ell^{r_i}\Z_\ell)
		\xrightarrow{\widetilde{\mu}_i^{-1}}
		\Gamma_i(r_i)
		\xrightarrow{t_i\mapsto\widetilde{\mu}_i(t_i)-1}
		\ell^{-r_i}\Z_\ell/\Z_\ell
		\xrightarrow{\times\ell^{r_i}}
		\Z_\ell/\ell^{r_i},
		\quad
		t\mapsto(t-1)/\ell
		\]
		induces 
		the following bijection
		\begin{equation}\label{logarithm becomes t to t-1}
			\mu_i'
			\colon
			\Gamma_i'(r_i')
			\xrightarrow{\widetilde{\mu}_i^{-1}}
			\Gamma_i(r_i)^{\ell^{r_i-r_i'}}
			\xrightarrow{t_i\mapsto\widetilde{\mu}_i(t_i)-1}
			\ell^{-(r_i-r'_i)}
			(\ell^{-r_i}\Z_\ell/\Z_\ell)
			\simeq
			\Z_\ell/\ell^{r_i'}.
		\end{equation}
		Recall we have
		$0<2r_i'<r_i$. So $\mu_i'$ is in fact a \emph{group homomorphism}. We can thus rewrite
		$\overline{F}_{\mathbf{f},\underline{r}',\underline{\alpha}'}$ as follows
		\begin{equation}\label{wildetilde{F}, another expression}
			\overline{F}_{\mathbf{f},\underline{r}',\underline{\alpha}'}(g)
			=
			\sum_{t''=(t_1,\cdots,t_{n})\in\prod_{i=1}^{n}\Gamma_i'(r_i')}
			\overline{F}_{\mathbf{f}}
			(
			g
			u
			(
			\ell^{-r_1'}\mu_1'(t_1),
			\cdots,
			\ell^{-r'_{n}}\mu'_{n}(t_{n}))
			)
			\psi_{\underline{r},\underline{\alpha}'}(t'').
		\end{equation}

		By definition, one has
		\begin{equation}\label{p=f}
			\sum_{\underline{\alpha}\in I(\underline{r}',\underline{r},\underline{\alpha}')}
			\overline{\mathbf{P}_{\mathbf{f}}}(\psi'\times
			\psi_{\underline{r},\underline{\alpha}},\varsigma_{\underline{r},0})			
			=
			\ell^{\sum_{i=1}^{n}(r_i-r_i')}
			\sum_{t'\in\mathcal{G}_0}
			\overline{F}_{\mathbf{f},\underline{r}',\underline{\alpha}'}
			(t'\varsigma_{\underline{r},0})
			\psi'(t').
		\end{equation}

		By Lemma \ref{widetilde(Q)_psi'' is not Spin-invariant} below,
		for any
		$0\ll r_1'\ll r_2'\ll\cdots\ll r_{n}'$,
		there exists
		$\underline{\alpha}'
		\in
		\prod_i
		(\Z_\ell/\ell^{r_i'})^\times$
		such that
		$\overline{F}_{\mathbf{f},\underline{r}',\underline{\alpha}'}$ is not invariant under $\widetilde{H}(\mathbb{A}_f)$.
		Next we want to apply Theorem \ref{automorphic application of appendix}, so we need to check that the assumption in this theorem is satisfied: indeed, we can take $\widetilde{\mathcal{G}}$ to be $\mathbf{T}(\mathbb{A}_f)$, $K$ to be the pre-image of $K^{(\mathbb{S})}$ under the map $\mathbf{H}(\mathbb{A}_f)\rightarrow H(\mathbb{A}_f)$. We have the spin norm map
		\[
		\mathrm{SN}\colon\mathbf{H}\rightarrow\mathbb{G}_{m},
		\]
		whose kernel is $\mathbf{H}^1$. By our choice of $\mathbb{S}$ and the condition \ref{condition-delta are squares in Q}, we know that $\mathrm{SN}(K)$ contains $\Z_q^\times$ for all primes $q$ such that $\mathrm{val}_q(\delta_1\cdots\delta_{2n})>0$. By Hasse-Minkowski theorem for the case $n\ge2$ and Hasse-Schilling theorem for $n=1$, we have
		\[
		\mathrm{SN}(\mathbf{H}(\Q))=\Q_{>0}.
		\]
		We then deduce
		\[
		\mathbb{A}_f^\times=\mathrm{SN}(\mathbf{H}(\Q)\mathbf{T}(\mathbb{A}_f)K),
		\]
		which implies in particular that the map $\widetilde{\mathcal{G}}\rightarrow\mathcal{Z}_K$ in Theorem \ref{automorphic application of appendix} is surjective.
		Taking into account of
		Lemma \ref{R satisfies condition on g_k,g_l},
		for any
		$0\ll r_1'\ll r_2'\ll\cdots\ll r_n'$,
		there is
		$h
		\in
		\mathbf{T}(\mathbb{A}_f)$
		such that
		\begin{align*}
			\sum_{t'\in\mathcal{G}_0}
			\overline{F}_{\mathbf{f},\underline{r}',\underline{\alpha}'}
			(
			t'h
			\varsigma_{\underline{r},0}
			)
			\psi'(t')
			=
			\sum_{t'\in\mathcal{R}}
			\overline{F}_{\mathbf{f},\underline{r}',\underline{\alpha}'}
			\left(
			\mathfrak{s}(t')h
			\mathfrak{s}(v_1(\ell^{-r'_1})\cdots v_{n}(\ell^{-r'_{n}}))
			\right)
			\psi'(\mathfrak{s}(t'))
			\neq0.
		\end{align*}
		Since
		$\mathcal{R}$
		and
		$\mathcal{R}h$
		are both mapped bijectively to $\mathcal{G}_0$ by the map (\ref{T(A_f)-->G_0}),
		it follows immediately that
		for such
		$\xi
		=\xi(\underline{r})$
		and $\psi'$,
		there is
		$\underline{\alpha}
		\in
		I(\underline{r}',\underline{r},\underline{\alpha}')$
		such that
		\[
		\overline{\mathbf{P}_{\mathbf{f}}}(\psi'\times\psi_{\underline{r},\underline{\alpha}},\varsigma_{\underline{r},0})
		\neq0.
		\]
		This finishes the proof of the theorem.
	\end{proof}

	\begin{remark}\label{no need of ramified places}
		\rm 
		Comparing with
		\cite[\S 5]{ChidaHsieh2016},
		there is no appearance of the group $\Delta^{\mathrm{alg}}$
		in our proof. Here is the explanation: in our case the group $\Delta^{\mathrm{alg}}$
		corresponds to the subgroup of
		$\mathbf{T}(\mathbb{A}_f)/\mathbb{G}_m(\mathbb{A}_f)
		T(\widehat{\Z})$
		generated by
		those
		$\sqrt{ q}$
		where
		$q$
		is a finite place of $\Q$
		divides exactly one of
		$\delta_{2i-1}$ and $\delta_{2i}$ (and not dividing both)
		for some $i=1,\cdots,n$.
		However under the map
		$\mathfrak{s}
		\colon
		\mathbf{T}(\mathbb{A}_f)
		\rightarrow
		T(\widehat{\Z})$,
		$\sqrt{ q}$
		is sent to
		an element of order $2$
		which is contained in $T(\widehat{\Z}^\ell)$.
		As a result
		the image of $\Delta^{\mathrm{alg}}$
		has trivial image in
		$\mathcal{G}_\infty$.
		That is why we do not need $\Delta^{\mathrm{alg}}$
		in our case.
	\end{remark}

	\begin{lemma}\label{widetilde(Q)_psi'' is not Spin-invariant}
		Let the notations be as in the proof of Theorem \ref{theorem on the non-vanishing mod p of toric periods}.
		For $0\ll r_1'\ll r_2'\ll\cdots\ll r_{n}'$,
		there is
		$\underline{\alpha}'\in
		\prod_i
		(\Z_\ell/\ell^{r_i'})^\times$
		such that the map
		$\overline{F}_{\mathbf{f},\underline{r}',\underline{\alpha}'}$
		defined in (\ref{widetilde{F}}) is not $\widetilde{H}^1(\mathbb{A}_f)$-invariant.
	\end{lemma}
	\begin{proof}
		Fix $\underline{r}=(r_1,\cdots,r_{n})$ as above
		such that
		$r_i>2r_i'$ for all $i$. Then we write
		\[
		\psi_{i,r_i'}=
		\psi_{i,r_i}^{\ell^{r_i-r_i'}}
		\colon
		\Gamma'_i(r_i')
		\rightarrow
		\mathbb{C}^\times
		\]
		where $\psi_{i,r_i}$
		is as in the proof of
		Theorem
		\ref{theorem on the non-vanishing mod p of toric periods}.
		For
		$i=1,\cdots,n$,
		we consider an element
		$\underline{\alpha}^{(i)}=(\alpha_1^{(i)},\cdots,\alpha_i^{(i)})
		\in
		\prod_{j=1}^i(\Z_\ell/\ell^{r_j'})^\times
		$ and set
		\[
		\psi_{\underline{r},\underline{\alpha}^{(i)}}
		=
		\prod_{j=1}^i
		\psi_{j,r_j'}^{\alpha_j^{(i)}}
		\]
		Then we define
		\begin{align*}
			\overline{F}_{\mathbf{f},\underline{r}',\underline{\alpha}^{(i)}}
			\colon
			H(\mathbb{A}_f)
			&
			\rightarrow
			\M_{\tau}(\kappa),
			\\
			g
			&
			\mapsto
			\sum_{t''=(t_1,\cdots,t_{i})\in\prod_{j=1}^{i}\Gamma'_j(r'_j)}
			\overline{F}_{\mathbf{f}}
			\left(
			g
			u_1(\ell^{-r'_1}\mu'_1(t_1))
			\cdots
			u_i(\ell^{-r'_i}\mu'_i(t_i))
			\right)
			\psi_{\underline{r},\underline{\alpha}^{(i)}}
			(t'').
		\end{align*}
		For $i=0$, we set
		$\overline{F}_{\mathbf{f},\underline{r}',\underline{\alpha}^{(0)}}
		=
		\overline{F}_{\mathbf{f}}$.
		Note that
		$\overline{F}_{\mathbf{f},\underline{r}',\underline{\alpha}^{(i)}}$
		is independent of
		$r_{i+1}',r_{i+2}',\cdots,r_{n}'$
		and
		$\overline{F}_{\mathbf{f},\underline{r}',\underline{\alpha}^{(n)}}$
		is the same as the one given in
		(\ref{widetilde{F}}) and (\ref{wildetilde{F}, another expression})
		for
		$\underline{\alpha}=\underline{\alpha}^{(n)}$.
		One verifies easily the following recursive relations among $\overline{F}_{\mathbf{f},\underline{r}',\underline{\alpha}^{(i)}}$:
		\begin{equation}\label{recursive relation}
			\overline{F}_{\mathbf{f},\underline{r}',\underline{\alpha}^{(i)}}
			(g)
			=
			\sum_{t_i\in\Gamma_i'(r_i')}
			\overline{F}_{\mathbf{f},\underline{r}',\underline{\alpha}^{(i-1)}}
			(gu_i(\ell^{-r_i'}\mu_i'(t_i)))
			\psi_{i,r_i'}^{\alpha_i^{(i)}}(t_i),
			\quad
			i=1,\cdots,n.
		\end{equation}

		We prove the lemma by contradiction and so assume that
		$\overline{F}_{\mathbf{f},\underline{r}',\underline{\alpha}^{(n)}}$
		is $\widetilde{H}^1(\mathbb{A}_f)$-invariant for all
		$\underline{\alpha}^{(n)}\in\prod_{j=1}^n(\Z/\ell^{r_j'})^\times$ for any $0\ll r_1'\ll r_2'\ll\cdots\ll r_n'$.
		We now use backward induction to deduce
		$\overline{F}_{\mathbf{f},\underline{r}',\underline{\alpha}^{(n-1)}}$,
		...,
		$\overline{F}_{\mathbf{f},\underline{r}',\underline{\alpha}^{(1)}}$
		and
		$\overline{F}_{\mathbf{f}}$		
		are all $\widetilde{H}(\mathbb{A}_f)$-invariant for all
		$\underline{\alpha}^{(n-1)},\cdots,\underline{\alpha}^{(1)}$ respectively.
		So by assumption
		$\overline{F}_{\mathbf{f},\underline{r}',\underline{\alpha}^{(n)}}$ is
		$\widetilde{H}^1(\mathbb{A}_f)$-invariant,
		in particular it is invariant under
		$\widetilde{H}^1( \Q_\ell)$.
		However
		the subgroup
		$\widetilde{H}^1(\Q_\ell)$
		has finite index in
		$ H( \Q_\ell)$
		(the index is $\leq[\Q_\ell^\times
		:
		(\Q_\ell^\times)^2]$).
		Thus
		we can take $r'_{n}\gg0$
		such that
		$u_{n}(\ell^{-r'_{n}})$
		lies in
		$\widetilde{H}^1(\Q_\ell)$.
		In particular,
		\[
		\overline{F}_{\mathbf{f},\underline{r}',\underline{\alpha}^{(n)}}
		(gu_{n}(\ell^{-r_{n}'}))
		=
		\overline{F}_{\mathbf{f},\underline{r}',\underline{\alpha}^{(n)}}(g),\,
		\forall
		g\in H(\mathbb{A}_f).
		\]
		However by definition
		(note that $\mu_{n}'$ is a group homomorphism!)
		\[
		\overline{F}_{\mathbf{f},\underline{r}',\underline{\alpha}^{(n)}}
		(gu_{n}(\ell^{-r_{n}'}))
		=
		\psi_{n,r_{n}'}^{\alpha_{n}^{(n)}}
		(1+\ell^{r_{n}-r_{n}'})
		\overline{F}_{\mathbf{f},\underline{r}',\underline{\alpha}^{(n)}}(g),
		\quad
		\forall
		g\in
		H(\mathbb{A}_f).
		\]
		The primitivity of the character
		$(\psi_{n,r_{n}'})^{\alpha_{n}^{(n)}}$
		implies that
		$\overline{F}_{\mathbf{f},\underline{r}',\underline{\alpha}^{(n)}}
		(g)
		=
		0$
		for all
		$g\in
		H(\mathbb{A}_f)$
		and
		$\underline{\alpha}^{(n)}
		\in
		\prod_{i=1}^{n}(\Z_\ell/\ell^{r'_i})^\times$.
		Using
		(\ref{recursive relation}),
		we claim
		that
		the map
		\[
		f\colon
		\Gamma_{n}'(r_{n}')
		\rightarrow
		\kappa,
		\quad
		t
		\mapsto
		\overline{F}_{\mathbf{f},\underline{r}',\underline{\alpha}^{(n-1)}}
		(gu_{n}(\ell^{-r_{n}'}\mu_{n}'(t))
		\]
		factors through the quotient
		$\Gamma_{n}'(r'_{n})
		\rightarrow
		\Gamma_{n}'(r'_{n}-1)$: indeed, write $\mathcal{S}$ for 
		$\kappa$-vector space
		of
		maps
		$f'
		\colon
		\Gamma_{n}'(r_{n}')
		\rightarrow
		\kappa$
		such that
		$\sum_t
		f'(t)
		\psi_{n,r_{n}'}^{\alpha_{n}^{(n)}}(t)
		=0$
		for all
		$\alpha_{n}^{(n)}
		\in
		(\Z_\ell/\ell^{r_{n}'})^\times$. Its dimension
		is equal to
		$\ell^{r_{n}'-1}$;
		on the other hand,
		any map
		$f'
		\colon
		\Gamma_{n}'(r_{n}')
		\rightarrow
		\kappa$
		induced from a map
		$f''\colon
		\Gamma_{n}'(r_{n}'-1)
		\rightarrow
		\kappa$
		via the above quotient map
		lies in this space $\mathcal{S}$,
		and this subspace has $\kappa$-dimension also equal to
		$\ell^{r_{n}'-1}$. This proves the claim.

		From the above argument, we deduce the following
		\begin{equation}\label{widetilde F is invariant under unipotent subgroup}
			\overline{F}_{\mathbf{f},\underline{r}',\underline{\alpha}^{(n-1)}}
			(g)
			=
			\overline{F}_{\mathbf{f},\underline{r}',\underline{\alpha}^{(n-1)}}
			(gu_{n}(\ell^{1-r_{n}'})),
			\quad
			\forall
			g\in
			H(\mathbb{A}_f).
		\end{equation}

		Fix an isomorphism
		$\mathbf{H}^1_{n}( \Q_\ell)
		\simeq
		\mathrm{SL}_2( \Q_\ell)$
		sending
		$\mathbf{H}^1_{n}(\Z_\ell)$
		onto
		$\mathrm{SL}_2(\Z_\ell)$,
		sending the torus
		$\mathbf{T}( \Q_\ell)
		\cap
		\mathbf{H}^1_{n}( \Q_\ell)$
		onto
		the diagonal subgroup
		and
		$u_{n}(\Z_\ell)$
		onto
		the subgroup
		$\begin{pmatrix}
			1 & \Z_\ell \\
			0 & 1
		\end{pmatrix}$.
		By definition,
		$\overline{F}_{\mathbf{f},\underline{r}',\underline{\alpha}^{(n-1)}}$
		is invariant under
		a compact open subgroup of
		$\mathbf{H}^1_{n}(\Z_\ell)$
		containing those elements whose image in
		$\mathrm{SL}_2(\Z_\ell)$
		are of the form
		$\begin{pmatrix}
			a & b \\
			c  & d
		\end{pmatrix}$
		with
		$a-1,d-1,c\in\ell^{e'}\Z_\ell$
		for some positive integer
		$e'$,
		depending only on
		$r_1',\cdots,r'_{n-1}$.
		Moreover
		(\ref{widetilde F is invariant under unipotent subgroup})
		implies that
		$\overline{F}_{\mathbf{f},\underline{r}',\underline{\alpha}^{(n-1)}}$
		is also invariant by
		$u_{n}
		(\ell^{-r'_{n}+1})$.
		Thus as long as we take
		$r'_{n}-1>e'$,
		Lemma
		\ref{a variant of a theorem of Ihara}
		shows that
		$\overline{F}_{\mathbf{f},\underline{r}',\underline{\alpha}^{(n-1)}}$
		is invariant under the subgroup
		$\mathbf{H}^1_{n}( \Q_\ell)$
		of
		$\mathbf{H}^1( \Q_\ell)$.
		Now apply
		Lemma
		\ref{SL_2 and K generate Spin},
		and we see that
		$\overline{F}_{\mathbf{f},\underline{r}',\underline{\alpha}^{(n-1)}}$
		is invariant under
		$\mathfrak{s}(\mathbf{H}^1( \Q_\ell))
		=
		\widetilde{H}^1(\Q_\ell)$
		for all
		$\underline{\alpha}^{(n-1)}
		\in
		\prod_{i=1}^{n-1}
		(\Z_\ell/\ell^{r_i'})^\times$,
		in particular, it is $\widetilde{H}^1(\mathbb{A}_f)$-invariant
		by the strong approximation property of
		$\mathbf{H}^1$.

		Now we repeat the above argument for the pair
		$(\overline{F}_{\mathbf{f},\underline{r}',\underline{\alpha}^{(n-1)}},
		\overline{F}_{\mathbf{f},\underline{r}',\underline{\alpha}^{(n-2)}})$
		instead of the pair
		$(\overline{F}_{\mathbf{f},\underline{r}',\underline{\alpha}^{(n)}},
		\overline{F}_{\mathbf{f},\underline{r}',\underline{\alpha}^{(n-1)}})$
		and we get that
		$\overline{F}_{\mathbf{f},\underline{r}',\underline{\alpha}^{(n-2)}}$
		is $\widetilde{H}^1(\mathbb{A}_f)$-invariant
		for all
		$\underline{\alpha}^{(n-2)}$, etc. until we get that
		$\overline{F}_{\mathbf{f},\underline{r}',\underline{\alpha}^{(0)}}
		=
		\overline{F}_{\mathbf{f}}$
		is $\widetilde{H}^1(\mathbb{A}_f)$-invariant, which is a contradiction.		
	\end{proof}

	Combining
	Theorems
	\ref{p-primitivity of theta lift of f}
	and
	\ref{theorem on the non-vanishing mod p of toric periods},
	we get the main result of this article
	\begin{theorem}\label{main result of the article}
		Assume conditions
		\ref{condition on lambda}, \ref{condition on p}, \ref{condition p not dividing the level} and
		\ref{condition-delta are squares in Q}.
		Let
		$\mathbf{f}
		\in
		\mathcal{A}_{\rho_{\lambda}}( H,K^{(\mathbb{S})})$ be $p$-primitive
		such that the image of
		$\overline{\mathbf{f}_p}$ generates $\M_{\lambda}(\kappa)$.
		Then the theta lift
		$\Theta_{\phi_{\lambda},\mathbf{f}}
		\in
		\mathcal{A}_{\rho_{\tau^\circ}}
		(\widetilde{G},\Gamma_0(2,N_{\mathbb{S}}),\chi_{U}^\circ)$
		is also
		$p$-primitive. 	
	\end{theorem}

	\appendix

	\section{Toric orbits of unipotent elements}\label{appendix}

	\subsection{Basic set-up}\label{basic set-up}
	Throughout this appendix we fix a prime number
	$\ell$.
	Write
	$\mathrm{U}_2$\index{U@$\mathrm{U}_2$}
	for the subgroup
	of
	$\mathrm{SL}_2$
	consisting of
	unipotent upper triangular matrices
	\[
	\mathrm{U}_2
	=
	\left\{
	\begin{pmatrix}
		1 & b \\ 0 & 1
	\end{pmatrix}\in\mathrm{SL}_2
	\right\}.
	\]

	We fix
	a connected linear algebraic group $\mathbf{H}$ over
	$\Q$
	and a maximal torus $\mathbf{T}$
	(over $\Q$)
	of $\mathbf{H}$,
	both of which are split at $\ell$.
	We assume that the derived subgroup
	$\mathbf{H}^1$ is simply connected and that
	$\mathbf{H}^1(\Q)$ is discrete in
	$\mathbf{H}^1(\mathbb{A}_f)$.
	We assume moreover that there are $n$ algebraic subgroups 
	$\mathbf{H}_1,\cdots,\mathbf{H}_n$
	of $\mathbf{H}$
	such that the following conditions are satisfied:

	\begin{taggedtheorem}{(A3)}\label{condition on mathbf{H}}
		\begin{enumerate}
			\item
			Each $\mathbf{H}_j$ is split at $\ell$, isomorphic to either $\mathbf{B}[j]^\times$ for some quaternion algebra $\mathbf{B}[j]$ over $ \Q$ or is a unitary group $\mathbf{U}[j]$ of rank $2$ over an imaginary quadratic field. We fix once and for all an isomorphism
			\[
			\iota_j\colon\mathrm{GL}_2(\Q_\ell)\simeq\mathbf{H}_j(\Q_\ell).
			\]		
			Then these subgroups $\mathbf{H}_j^1(\Q_\ell)$ ($j=1,\cdots,n$)			generate the whole $\mathbf{H}^{1}( \Q_\ell)$.

			\item
			The unipotent subgroups
			\[
			\mathbf{U}_j:=
			\iota_j(\mathrm{U}_2( \Q_\ell))
			\]
			all commute with each other. Moreover $\mathbf{T}(\Q_\ell)$ normalizes $\mathbf{H}_j(\Q_\ell)$ and $\mathbf{U}_j$ for $j=1,\cdots,n$ and $\mathbf{T}(\Q_\ell)$ acts by conjugation on $\mathbf{U}_j$ via an algebraic character $\chi_j$ and that these characters $\chi_j$ are linearly independent in the character group $X^\ast(\mathbf{T}_{/\Q})$ of $\mathbf{T}_{/\Q}$.\index{A@(A3)}
		\end{enumerate}		
	\end{taggedtheorem}

	\begin{remark}\rm
		\begin{enumerate}
			\item
			Note that
			the simple-connectedness of $\mathbf{H}^1$ implies that 
			$\mathbf{H}^{1}(\mathbb{A}_f)$
			satisfies the strong approximation property
			(with respect to the place $\ell$).

			\item 
			By
			\cite[Proposition 1.4]{Gross1999},
			we know that
			$\mathbf{H}^1(\Q)$ is discrete in
			$\mathbf{H}^1(\mathbb{A}_f)$
			if and only if
			$\mathbf{H}^1(\Q)$ is discrete and co-compact in
			$\mathbf{H}^1(\mathbb{A}_f)$,
			if and only if
			$\mathbf{H}^1(\R)$
			is compact. By assumption, these
			$\mathbf{H}_j^{1}(\R)$
			are all compact,
			$\mathbf{H}_j^{1}$
			are all simply-connected and
			$\mathbf{H}_j^{1}(\mathbb{A}_f)$
			satisfy the strong approximation property
			(with respect to the place $\ell$).

			\item 
			The two groups
			$\mathbf{B}[j]^\times$ and
			$\mathbf{U}[j]$
			are closely related
			as follows:
			for each quaternion algebra $\mathbf{B}$ over $\Q$
			with a principal involution $\tau$,
			we can define a dimension $2$ Hermitian space
			over a quadratic extension $K$ of $\Q$
			(which is a maximal commutative subalgebra of $\mathbf{B}$
			and is stable under the involution $\tau$)
			and write
			$\mathbf{SU}$
			for the corresponding
			special unitary group
			associated to this space.
			Write
			$\mathbf{B}^1$
			for the set of elements of
			$\mathbf{B}$
			whose reduced norm is $1$.
			Then
			$\mathbf{B}^1\simeq\mathbf{SU}$
			as algebraic groups over $\Q$.
			So the groups
			$\mathbf{H}_j$
			can all be seen as extensions of
			a rank one torus (over $\Q$)
			by some
			$\mathbf{B}^1$.
			From (2) we have an isomorphism
			$\prod_{j=1}^n\mathbf{U}_j\simeq
			\Q_\ell^n$
			of $\ell$-adic groups.
			Moreover
			$\mathbf{T}(\Q_\ell)$
			normalizes each
			$\mathbf{H}_j^1(\Q_\ell)$.
		\end{enumerate}    	
	\end{remark}

	The following are two examples of $\mathbf{H}$ satisfying the above assumptions.
	We refer to
	§\ref{Example}
	for other more examples.
	\begin{examples}\label{examples}
		\begin{enumerate}
			\item
			$\mathbf{H}$ is $\mathbf{B}^\times$
			where
			$\mathbf{B}$ is
			a definite quaternion algebra over $\Q$;

			\item 
			$\mathbf{H}$ is
			a unitary group $\mathbf{U}$ of rank $2$ over an imaginary quadratic field.
		\end{enumerate}
		
		We assume that $\mathbf{H}$ is split at $\ell$.
		In this case we fix a maximal torus 
		$\mathbf{T}$ of $\mathbf{H}$
		which is also split at $\ell$.
		Then we take
		$n=1$,
		$\mathbf{H}_1=\mathbf{H}$,
		$\mathbf{U}_1$
		is a unipotent subgroup of $\mathbf{H}^1(\Q_\ell)$
		normalized by
		$\mathbf{T}(\Q_\ell)$.
	\end{examples}

	\subsection{Ratner's theorems and commensurability}

	\subsubsection{SL(2) case}
	In this subsection
	we recall Ratner's orbit closure theorem
	and uniform distribution theorem
	on unipotent flows.
	Write
	$\mathbf{m}$\index{m@$\mathbf{m}$}
	to be the normalized Haar measure on
	$ \Q_\ell$.
	For a subset
	$\kappa$ of $\Q_\ell$ and an integer $N$,
	we write
	\[
	\kappa_N=\kappa\ell^{-N}
	=
	\{k\ell^{-N}|k\in\kappa\}.
	\index{k@$\kappa_N$}
	\]

	Let
	$\mathbf{H}$ be as in the preceding section and put
	$\mathbf{G}=\mathbf{H}^1(\Q_\ell)$. A lattice in $\mathbf{G}$ is a discrete and cocompact subgroup.
	Then Ratner's theorems give
	(\cite[Theorems 2\&3]{Ratner1995}):
	\begin{theorem}\label{Ratner}		
		Let
		$\Gamma$ be a lattice in
		$\mathbf{G}$
		and
		$\mathbf{U}$ be
		a subgroup of $\mathbf{G}$
		generated by one-parameter unipotent subgroups of
		$\mathbf{G}$.
		\begin{enumerate}
			\item
			(\textbf{Orbit closure theorem})
			For any
			$x\in\Gamma\backslash\mathbf{G}$,
			the closure
			$x\mathbf{U}$ in
			$\Gamma\backslash\mathbf{G}$
			is of the form
			$x\mathbf{L}$
			for some closed subgroup $\mathbf{L}$ of $\mathbf{G}$
			containing $\mathbf{U}$.

			\item 
			(\textbf{Uniform distribution theorem})
			Let $\mathbf{U}=\{\mathbf{u}(t)|t\in \Q_\ell\}$
			be a one-parameter unipotent subgroup of $\mathbf{G}$.
			Write
			$\mu_{\mathbf{L}}$
			for the unique Borel measure on
			$\Gamma\backslash\mathbf{G}$
			invariant under the action of $\mathbf{L}$
			and supported on
			$x\mathbf{L}$.
			Then
			for any locally constant function
			$f\colon
			\Gamma\backslash\mathbf{G}
			\rightarrow
			\mathbb{C}$
			and any compact open subset
			$\kappa\subset \Q_\ell$,
			one has
			\[
			\lim\limits_{N\rightarrow+\infty}
			\frac{1}{\mathbf{m}(\kappa_N)}
			\int_{\kappa_N}
			f(x\mathbf{u}(t))dt
			=
			\int_{\Gamma\backslash\mathbf{G}}
			fd\mu_{\mathbf{L}}.
			\]
		\end{enumerate}
	\end{theorem}
	
	\begin{remark}\rm
		In \textit{loc.cit}, $f$ is assumed to be continuous. However, since
		$\Gamma\backslash\mathbf{G}$ is compact, $f$ can be uniformly approximated by locally constant functions on $\Gamma\backslash\mathbf{G}$. Thus the last conclusion in the above theorem is equivalent to the one given in
		\textit{loc.cit}.
	\end{remark}

	The proofs of the above results rely on a careful study for the case
	(Theorem 6 of \textit{loc.cit}),
	\[
	\mathbf{G}=\mathrm{SL}_2(\Q_\ell)^r,
	\]
	which we will also need in the following. Let $\Gamma_1,\cdots,\Gamma_r$
	be lattices in $\mathrm{SL}_2( \Q_\ell)$ and write
	\[
	\Gamma=\Gamma_1\times\cdots\times\Gamma_r.
	\]
	We write $\Delta$ the diagonal embedding
	\[
	\Delta
	\colon
	\mathrm{SL}_2( \Q_\ell)
	\rightarrow
	\mathrm{SL}_2( \Q_\ell)^r
	\]
	and $V$ a one-parameter unipotent subgroup of $\mathrm{SL}_2( \Q_\ell)$
	\[
	V=\{v(t)|t\in \Q_\ell\}.
	\]
	So $\Delta(V)$ is a one-parameter unipotent subgroup of $\mathrm{SL}_2(\Q_\ell)^r$.
	\begin{theorem}
		For any $g=(g_1,\cdots,g_r)\in\mathrm{SL}_2( \Q_\ell)^r$, the closure
		$\Gamma g\Delta(V)$ inside $\mathrm{SL}_2( \Q_\ell)^r$ is of the form
		$\Gamma g\mathbf{L}$ for a closed subgroup $\mathbf{L}$ of	$\mathrm{SL}_2( \Q_\ell)^r$ containing $\Delta(V)$ and there is an element
		$c\in V^r$ such that $c\mathbf{L} c^{-1}\supset\Delta(\mathrm{SL}_2(\Q_\ell))$. Moreover there is a unique
		$\mathbf{L}$-invariant Borel measure $\mu_\mathbf{L}$ on	$\Gamma\backslash\mathrm{SL}_2( \Q_\ell)^r$ supported on
		$\Gamma g\mathbf{L}$, and the measure $\mu_\mathbf{L}$
		is ergodic for $\mathbf{L}$: for any locally constant function $f\colon
		\Gamma\backslash\mathrm{SL}_2( \Q_\ell)^r\rightarrow\mathbb{C}$
		and any compact open subset $\kappa\subset\Q_\ell$,
		\[
		\lim\limits_{N\rightarrow+\infty}
		\frac{1}{\mathbf{m}(\kappa_N)}
		\int_{\kappa_N}
		f(gv(t))dt
		=
		\int_{\Gamma\backslash\mathrm{SL}_2( \Q_\ell)^r}
		fd\mu_\mathbf{L}.
		\]
	\end{theorem}
	The part $c\mathbf{L}c^{-1}\supset\Delta(\mathrm{SL}_2(\Q_\ell))$ for some
	$c\in V^r$ comes from Theorem 6 of \cite{Ratner1995} and Theorem 1.1 of \cite{Shah2009}. We can say something more for the closed subgroup
	$\mathbf{L}$ under certain conditions. For this we need the notion of
	$V$-commensurability: we say two lattices $\Gamma,\Gamma'$ of
	$\mathrm{SL}_2(\Q_\ell)$ are \textit{$V$-commensurable}, if there is an element $v\in V$ such that $\Gamma$ and $v\Gamma' v^{-1}$ are commensurable (that is, $\Gamma\cap v\Gamma' v^{-1}$ has finite index in both $\Gamma$ and $v\Gamma' v^{-1}$). Then one has (\cite[Proposition 2.35]{CornutVatsal2005})
	\begin{theorem}\label{not V-commensurable}
		Maintain the notations of the preceding Theorem. If for any $i\ne j\in\{1,\cdots,r\}$, the lattices $g_i^{-1}\Gamma_ig_i$ and $g_j^{-1}\Gamma_jg_j$ are \emph{not} $V$-commensurable, then
		$\mathbf{L}=\mathrm{SL}_2( \Q_\ell)^r$.
	\end{theorem}
	\begin{remark}\rm
		In fact the converse is also true (see \textit{loc.cit}). But we will not need this in the following.
	\end{remark}
	
	\subsubsection{Adelic reformulation}\label{adelic formulation}
	For later applications,
	it is useful to give an adelic point of view of the preceding result.
	Let $\mathbf{H}$ be as in
	Examples
	\ref{examples}.

	Fix a compact open subgroup
	$K^\ell$ of
	$\mathbf{H}(\mathbb{A}_f^\ell)$
	and let
	$\mathbf{H}^1(\Q_\ell)$
	act on the right on
	$\mathbf{H}(\Q)\backslash \mathbf{H}(\mathbb{A}_f)/K^\ell$.
	For an element
	$g\in \mathbf{H}(\mathbb{A}_f)$,
	write
	$\Gamma_{K^\ell}(g)$\index{G@$\Gamma_{K^\ell}(g)$}
	for the stabilizer inside
	$\mathbf{H}^{1}( \Q_\ell)$
	of the double coset
	$\mathbf{H}(\Q)gK^\ell
	\in
	\mathbf{H}(\Q)\backslash
	\mathbf{H}(\mathbb{A}_f)/K^{\ell}$.
	\begin{lemma}\label{stabilizer for SL(2)}
		The stabilizer
		$\Gamma_{K^\ell}(g)$
		is a lattice of
		$\mathbf{H}^{1}( \Q_\ell)$.
	\end{lemma}
	\begin{proof}
		Since $\mathbf{H}^1\subset \mathbf{H}$ are linear algebraic groups over $\Q$,
		\[
		\mathbf{H}(\Q)\cap \mathbf{H}^1(\Q_\ell)=\mathbf{H}^1(\Q).
		\]
		Write $\Gamma$ to be the image of the projection map $\mathbf{H}^1(\mathbb{A}_f)\to\mathbf{H}^1(\Q_\ell)$ of the following subgroup
		\[
		\mathbf{H}(\Q)
		\cap
		(gK^\ell g^{-1}\times \mathbf{H}^1(\Q_\ell))
		=
		\mathbf{H}^1(\Q)
		\cap
		\left(
		gK^{\ell}g^{-1}\times \mathbf{H}^1(\Q_\ell)
		\right)
		\subset
		\mathbf{H}^1(\mathbb{A}_f).
		\]
		Then it is easy to see that
		\[
		\Gamma_{K^\ell}(g)
		=
		g_\ell^{-1}\Gamma g_\ell.
		\]

		Write
		$W=
		gK^{\ell}g^{-1}\times \mathbf{H}^1(\Q_\ell)
		\subset
		\mathbf{H}^1(\mathbb{A}_f)$.
		Then the continuous injective map
		\[
		(\mathbf{H}^1(\Q)\cap W)\backslash W
		\rightarrow
		\mathbf{H}^1(\Q)\backslash \mathbf{H}^1(\mathbb{A}_f)
		\]
		is open since
		$W$ is open in
		$\mathbf{H}^1(\mathbb{A}_f)$
		and is also
		surjective since
		$\mathbf{H}^1(\mathbb{A}_f)$
		satisfies the strong approximation property.
		Thus it is a homeomorphism.
		By assumption
		$\mathbf{H}^1(\Q)$ is a lattice inside
		$\mathbf{H}^1(\mathbb{A}_f)$.
		This implies that
		$\mathbf{H}^1(\Q)\cap W$
		is a lattice in
		$W$ and since the factor
		$gK^{\ell}g^{-1}$
		is compact,
		$\Gamma$
		is also a lattice in
		$\mathbf{H}^1(\Q_\ell)$
		by
		\cite[p.105, Lemme 1.2]{Vigneras1980}.
		Thus
		$\Gamma_{K^\ell}(g)$
		is a lattice in
		$\mathbf{H}^1(\Q_\ell)$.
	\end{proof}
	Since
	$\mathbf{H}^{1}(\mathbb{A}_f)$
	satisfies strong approximation property
	with respect to $\ell$, one has
	\[
	g^{-1}\mathbf{H}(\Q)g
	\mathbf{H}^{1}( \Q_\ell)K^\ell
	=
	g^{-1}\mathbf{H}(\Q)g
	\mathbf{H}^{1}(\mathbb{A}_f)
	K^{\ell}.
	\]
	One deduces that
	the following natural map is a homeomorphism of topological spaces
	\begin{align*}
		\Gamma_{K^\ell}(g)
		\backslash
		\mathbf{H}^{1}( \Q_\ell)
		&
		\simeq
		\mathbf{H}(\Q)
		\backslash
		\mathbf{H}(\Q)g\mathbf{H}^{1}(\mathbb{A}_f)
		K^{\ell}
		/K^\ell,
		\\
		\Gamma_{K^\ell}(g)h
		&
		\mapsto
		\mathbf{H}(\Q)
		ghK^\ell
		=
		\mathbf{H}(\Q)
		gK^\ell h.
	\end{align*}
	Note that this homeomorphism is $\mathbf{H}^{1}( \Q_\ell)$-equivariant.
	Thus the measure
	$\mu_{\mathbf{L}}$on
	$\Gamma_{K^\ell}(g)\backslash \mathbf{H}^{1}( \Q_\ell)$
	as in Theorem \ref{Ratner}
	corresponds to a measure
	on the RHS
	which is $\mathbf{L}$-invariant. We denote this measure by $\mu_{g,V}$.
	So
	Theorem \ref{Ratner} gives
	\begin{corollary}
		For any locally constant function
		$f\colon
		\mathbf{H}(\Q)
		\backslash
		\mathbf{H}(\mathbb{A}_f)
		/K^\ell
		\rightarrow
		\mathbb{C}$
		and any compact open subset
		$\kappa$ of
		$ \Q_\ell$,
		\[
		\lim\limits_{N\rightarrow+\infty}
		\frac{1}{\mathbf{m}(\kappa_N)}
		\int_{\kappa_N}
		f(gv(t))dt
		=
		\int_{\mathbf{H}(\Q)
			\backslash
			\mathbf{H}(\Q)g\mathbf{H}^{1}(\mathbb{A}_f)
			K^{\ell}
			/K^\ell}
		fd\mu_{g,V}.
		\]
	\end{corollary}
	We want to put some conditions on $g$
	to ensure that
	$\mathbf{L}$
	is as large as possible,
	as in
	Theorem \ref{not V-commensurable}.
	We proceed this in two steps
	in the next two subsections.

	\subsubsection{Single-copy case}
	We now consider a higher-dimensional generalization of the above results.
	In this subsection we treat the case $r=1$.

	Let
	$\mathbf{H}$ be as in
	§\ref{basic set-up}.
	We fix isomorphisms
	\[
	u_j\colon
	\Q_\ell
	\simeq
	\mathbf{U}_j,
	\quad
	t\mapsto u_j(t),
	\quad
	\forall
	j=1,\cdots,n.
	\]
	Let
	$g\in
	\mathbf{T}(\mathbb{A}_f)$
	and
	$K^\ell$ be a compact open subgroup of
	$\mathbf{H}(\mathbb{A}_f^\ell)$.
	We have
	\begin{lemma}
		The stabilizer
		$\Gamma_{K^\ell}(g)$
		is a lattice in
		$\mathbf{H}^{1}( \Q_\ell)$.
		Similarly
		the intersection
		\[
		\Gamma_{K^\ell}^{(j)}(g)
		:=
		\Gamma_{K^\ell}(g)
		\cap
		\mathbf{H}_j^{1}
		( \Q_\ell)
		\]
		is a lattice in
		$\mathbf{H}_j^{1}
		( \Q_\ell)$.
	\end{lemma}
	\begin{proof}
		Write $\Gamma$ for the image under the projection map
		$\mathbf{H}^1(\mathbb{A}_f)\to\mathbf{H}^1(\Q_\ell)$ of the following subgroup
		\[
		\mathbf{H}(\Q)
		\cap
		\left(
		gK^\ell g^{-1}\times \mathbf{H}^{1}( \Q_\ell)
		\right)
		\subset
		\mathbf{H}^1(\mathbb{A}_f).
		\]
		Then we have
		\[
		\Gamma_{K^\ell}(g)
		=
		g_\ell^{-1}
		\Gamma
		g_\ell.
		\]

		Write $\Gamma_j
		=
		\Gamma
		\cap
		\mathbf{H}_j^1(\Q_\ell)$
		for the projection image to
		$\mathbf{H}_j^{1}( \Q_\ell)$
		of the following subgroup
		\[
		\mathbf{H}_j(\Q)
		\cap
		\left(
		gK^\ell g^{-1}
		\times
		\mathbf{H}_j^{1}( \Q_\ell)
		\right)
		\subset
		\mathbf{H}_j^{1}(\mathbb{A}_f).
		\]
		Since
		$g_\ell\in \mathbf{T}(\Q_\ell)$
		normalizes
		$\mathbf{H}_j(\Q_\ell)$
		and
		$\mathbf{H}_j^1(\Q_\ell)$, we have
		\[
		\Gamma_{K^\ell}^{(j)}(g)
		=
		g_\ell^{-1}\Gamma_jg_\ell.
		\]
		
		Now the same proof as for
		Lemma
		\ref{stabilizer for SL(2)}
		shows that
		$\Gamma$,
		resp.
		$\Gamma_j$
		is a lattice inside
		$\mathbf{H}^{1}( \Q_\ell)$,
		resp.
		$\mathbf{H}_j^{1}( \Q_\ell)$
		(using the fact that $\mathbf{H}^{1}(\mathbb{A}_f)$,
		resp.
		$\mathbf{H}_j^{1}(\mathbb{A}_f)$
		satisfies the strong approximation property
		and
		$\mathbf{H}^1(\Q)$,
		resp.
		$\mathbf{H}_j^1(\Q)$
		is a lattice in
		$\mathbf{H}^1(\mathbb{A}_f)$,
		resp.
		$\mathbf{H}_j^1(\mathbb{A}_f)$).
	\end{proof}
	
	For a lattice $\Gamma'$ in $\mathbf{H}_j^1(\Q_\ell)$, we define the commensurator of $\Gamma'$ inside $\mathbf{H}_j(\Q_\ell)$ to be
	\[
	\mathcal{C}_j(\Gamma')
	:=
	\left\{
	h\in
	\mathbf{H}_j( \Q_\ell)|
	\Gamma'
	\text{ and }h\Gamma'h^{-1}
	\text{are commensurable}
	\right\}.
	\]
	\begin{remark}\label{commensurable set}\rm
		By \cite[Lemma 2.19]{CornutVatsal2005}, we know
		\begin{equation*}
			\mathcal{C}_j(\Gamma_{K^\ell}^{(j)}(g))
			=
			\mathbf{H}_j(\Q)
			Z(\mathbf{H}_j( \Q_\ell)),
			\quad
			j=1,\cdots,n
		\end{equation*}		
		Since
		all compact open subgroups $K^{\ell}$ of
		$\mathbf{H}^1(\mathbb{A}_f^{\ell})$
		are commensurable,
		$\mathcal{C}_j(\Gamma_{K^\ell}^{(j)}(g))$
		does not depend on
		$K^{\ell}$.
	\end{remark}
	
	For a compact open subset
	$\kappa$ of
	$ \Q_\ell$, an $n$-tuple of integers
	$\underline{N}=(N_1,N_2,\cdots,N_n)
	\in\mathbb{N}^n$
	and an element
	$\underline{t}=(t_1,\cdots,t_n)
	\in
	\Q_\ell^n$,
	write
	\begin{align*}
		\kappa_{\underline{N}}
		&
		=\kappa_{N_1}\times\cdots\times\kappa_{N_n}
		\subset
		\Q_\ell^n,\index{k@$\kappa_{\underline{N}}$}
		\\
		\mathbf{m}(\kappa_{\underline{N}})
		&
		=
		\mathbf{m}(\kappa_{N_1})
		\times
		\cdots
		\times
		\mathbf{m}(\kappa_{N_n}),
		\\
		u(\underline{t})
		&
		=
		u_1(t_1)
		\times
		\cdots
		\times
		u_n(t_n)
		\in\prod_{j=1}^n\mathbf{U}_j.\index{u@$u(\underline{t})$}
	\end{align*}
	Moreover we write
	$\underline{N}\rightarrow+\infty$
	to mean
	$N_1,N_2,\cdots,N_n\rightarrow
	+\infty$.	
	Then
	\begin{theorem}\label{single-copy}
		The subset
		$\Gamma_{K^\ell}(g)
		\prod_{j=1}^n\mathbf{U}_j$
		is dense in $\mathbf{H}^{1}(\Q_\ell)$. Moreover, if we write $\mu_g$
		for the unique Borel measure on
		$\mathbf{H}(\Q)
		\backslash
		\mathbf{H}(\Q)
		g
		\mathbf{H}^{1}(\mathbb{A}_f)K^{\ell}
		/K^{\ell}$
		invariant under
		$\mathbf{H}^{1}(\Q_\ell)$, then for any locally constant complex-valued function
		$f$ on 
		$\mathbf{H}(\Q)
		\backslash
		\mathbf{H}(\mathbb{A}_f)
		/K^\ell$
		and any compact open subset
		$\kappa$ of $ \Q_\ell$, we have the following equidistribution result
		\begin{equation}\label{equidistribution result}
			\lim\limits_{\underline{N}\rightarrow+\infty}
			\frac{1}{\mathbf{m}(\kappa_{\underline{N}})}
			\int_{\kappa_{\underline{N}}}
			f(gu(\underline{t}))
			d\underline{t}
			=
			\int_{\mathbf{H}(\Q)\backslash \mathbf{H}(\Q)g\mathbf{H}^{1}(\mathbb{A}_f)K^{\ell}/K^\ell}
			fd\mu_g.
		\end{equation}
	\end{theorem}

	\begin{proof}
		To ease notations, we write
		\[
		X=\mathbf{H}( \Q)\backslash \mathbf{H}(\Q)g\mathbf{H}^{1}(\mathbb{A}_f)K^{\ell}/K^\ell.
		\]
		
		For the density,
		by
		Theorem
		\ref{Ratner},
		the closure of
		$\Gamma_{K^\ell}(g)\prod_j\mathbf{U}_j$
		is of the form
		$\Gamma_{K^\ell}(g)\mathbf{L}$
		for a closed subgroup $\mathbf{L}$ of
		$\mathbf{H}^{1}( \Q_\ell)$
		containing
		$\prod_j\mathbf{U}_j$.
		So for any
		$h_j\in
		\mathbf{H}_j^{1}( \Q_\ell)$
		and any open subset
		$W_j$ of
		$\mathbf{H}_j^{1}( \Q_\ell)$
		containing $h_j$, we have
		\[
		\left(
		\Gamma_{K^\ell}^{(j)}(g)\mathbf{U}_j
		\right)
		\bigcap
		\left(
		\Gamma_{K^\ell}^{(j)}(g)h_jW_j
		\right)
		\ne
		\emptyset.
		\]
		Thus for any open subset
		$W$ of
		$\mathbf{H}^{1}( \Q_\ell)$
		containing $h_j$,
		the intersection
		$W\cap \mathbf{H}_j^1(\Q_\ell)$
		is open in
		$\mathbf{H}_j^1(\Q_\ell)$
		for any $j=1,\cdots,n$
		and
		one deduces
		from the previous non-empty intersection
		\[
		\Big(
		\Gamma_{K^\ell}(g)\prod_{j'}\mathbf{U}_{j'}
		\Big)
		\bigcap
		\Big(
		\Gamma_{K^\ell}(g)h_jW
		\Big)
		\ne\emptyset.
		\]
		Therefore
		$\mathbf{L}$
		contains all these
		$\mathbf{H}_j^{1}( \Q_\ell)$,
		which 
		generate the whole
		$\mathbf{H}^{1}( \Q_\ell)$,
		thus
		$\mathbf{L}
		=\mathbf{H}^{1}( \Q_\ell)$.

		One deduces that the measure
		$\mu_{\mathbf{L}}$
		is the unique Borel measure on
		$\Gamma_{K^\ell}(g)\backslash \mathbf{H}^1(\Q_\ell)$
		which is $\mathbf{H}^1(\Q_\ell)$-invariant.
		This measure is transferred
		to the unique Borel measure
		on
		$X$,
		via \emph{transport de structure} by the following natural
		$\mathbf{H}^1(\Q_\ell)$-equivariant
		homeomorphism	
		\[
		\Gamma_{K^\ell}(g)
		\backslash
		\mathbf{H}^1(\Q_\ell)
		\simeq
		X.
		\]
		This is the measure $\mu_g$ in the theorem.

		To prove (\ref{equidistribution result}),
		we argue as follows:
		for any
		$\underline{N}$
		as above,
		we define a Borel probability measure
		$\mu_{\underline{N}}$ on
		$X$
		by the formula
		\[
		\int_Xf(x)d\mu_{\underline{N}}(x)
		=
		\frac{1}{\mathbf{m}(\kappa_{\underline{N}})}
		\int_{\underline{t}\in\kappa_{\underline{N}}}
		f(u(\underline{t}))d\underline{t}
		\quad
		\text{for any continuous/locally constant}
		f\colon
		X
		\rightarrow
		\mathbb{C}.
		\]
		Since $X$ is compact,
		for any sequence
		$\underline{N}_{(k)}=(N_{1,k},\cdots,N_{n,k})$
		with $N_{i,k}\rightarrow+\infty$
		for $k\rightarrow+\infty$,
		there is a subsequence
		$\{\underline{N}_{(k_s)}\}_{s\in\mathbb{N}}$
		such that
		$\mu_{\underline{N}_{(k_s)}}$
		converges
		(under the weak topology) to
		a Borel probability measure
		$\mu'$ on $X$.
		We claim that
		$\mu'$ is $\mathbf{H}^1(\Q_{\ell})$-invariant
		and thus we necessarily have $\mu'=\mu_g$,
		which finishes the proof of the theorem.

		To prove the claim, note that
		$\prod_{j=1}^n\mathbf{U}_j$
		preserves $\mu'$
		(because $\kappa$ is a compact \textit{open} subset of $\Q_{\ell}$).
		For $h\in \mathbf{H}^1(\Q_{\ell})$,
		write
		$\mu'\circ h$
		for the right translation of $h$ on $\mu'$,
		that is,
		\[
		\int_Xf(x)d(\mu'\circ h)(x)
		:=
		\int_Xf(xh)d\mu'(x).
		\]
		Any locally constant function
		$f\colon
		X\rightarrow
		\mathbb{C}$
		is invariant under the right translation by a compact open subgroup
		$K'$ of $\mathbf{H}^1(\Q_{\ell})$.
		Thus we have
		\[
		\int_Xf(x)d(\mu'\circ k)(x)
		=
		\int_Xf(x)d\mu'(x),
		\quad
		\forall
		\,
		k\in K'.
		\]
		Fix $j=1,\cdots,n$, since $\mathbf{H}_j^1(Q_\ell)\cap K'$ is an open subgroup in $\mathbf{H}_j^1(\Q_\ell)$,
		$\mathbf{U}_j$ and $\mathbf{H}_j^1(\Q_{\ell})\cap K'$
		generate $\mathbf{H}_j^1(\Q_{\ell})$,
		then one deduces
		\[
		\int_X
		f(x)d(\mu'\circ h_j)(x)
		=
		\int_Xf(x)d\mu'(x),
		\quad
		\forall
		\,
		h_j\in
		\mathbf{H}_j^1(\Q_{\ell}).
		\]
		Since this is true for arbitrary locally constant function $f$,
		we deduce that
		$\mu'$ is invariant under
		$\mathbf{H}_j^1(\Q_{\ell})$
		for any $j$,
		and therefore also invariant under
		$\mathbf{H}^1(\Q_{\ell})$
		(it is generated by these $\mathbf{H}_j^1(\Q_{\ell})$).
		This proves our claim.
	\end{proof}

	\subsubsection{Multi-copy case}
	In this subsection we treat the case
	$r>1$.

	Again let $\mathbf{H}$ be as in the preceding section,
	$K^\ell$
	be a compact open subgroup of
	$\mathbf{H}(\mathbb{A}_f^\ell)$
	and
	$g=(g_1,\cdots,g_r)$
	with
	$g_1,\cdots,g_r\in \mathbf{T}(\mathbb{A}_f)$
	such that their $\ell$-th components satisfy
	\begin{equation}\label{pairwise distinct-1}
		(g_k)_\ell(g_i)_\ell^{-1}
		\notin
		\mathbf{T}(\Q_\ell)
		\cap
		\big(
		\cup_{j=1}^n
		\mathbf{H}_j(\Q)Z(\mathbf{H}_j(\Q_\ell))
		\big),
		\quad
		\forall
		k\ne i
		\in
		\{
		1,\cdots,r
		\}.
	\end{equation}
	As in the preceding subsection,
	$\Gamma_{K^\ell}^{(j)}(g_i)
	$
	is a lattice inside
	$\mathbf{H}_j^{1}( \Q_\ell)$.
	We write
	\begin{align*}
		\Gamma_{K^\ell}^{(j)}(g)
		&
		=
		\Gamma_{K^\ell}^{(j)}(g_1)
		\times
		\cdots
		\times
		\Gamma_{K^\ell}^{(j)}(g_r),
		\\
		\Gamma_{K^\ell}(g)		
		&
		=
		\Gamma_{K^\ell}(g_1)
		\times
		\cdots
		\times
		\Gamma_{K^\ell}(g_r),
	\end{align*}
	which is a lattice in
	$\mathbf{H}_j^{1}( \Q_\ell)^r$,
	resp.
	$\mathbf{H}^{1}( \Q_\ell)^r$.
	Then 
	we have
	\begin{lemma}
		Fix an element
		$h\in \mathbf{T}(\mathbb{A}_f)$.
		For any $j=1,\cdots,n$,
		$\Gamma_{K^\ell}^{(j)}(\Delta(h)g)
		\Delta(\mathbf{U}_j)$
		is dense in
		$\mathbf{H}_j^{1}( \Q_\ell)^r$.
		Similarly
		$\Gamma_{K^\ell}(\Delta(h)g)
		\Delta(\prod_{j=1}^n\mathbf{U}_j)$
		is dense in
		$\mathbf{H}^{1}( \Q_\ell)^r$.
	\end{lemma}
	\begin{proof}		
		\sloppy
		Fix $i\ne k\in\{1,\cdots,r\}$.
		We claim there is no
		$w\in \mathbf{U}_j$
		such that
		$\Gamma_{K^\ell}^{(j)}(h g_i)$
		and
		$w\Gamma_{K^\ell}^{(j)}(h g_k)w^{-1}$
		are commensurable:
		otherwise,
		recall the expression for
		$\Gamma_{K^\ell}^{(j)}(h g_i)$,
		using
		Remark
		\ref{commensurable set}, we have the following
		\begin{align*}
			(h g_i)_\ell^{-1}
			\mathcal{C}_j(\Gamma_{K^\ell}^{(j)}(1))
			(h g_i)_\ell
			&
			=
			\mathcal{C}_j(\Gamma_{K^\ell}^{(j)}(h g_i))
			=
			\mathcal{C}_j(w\Gamma_{K^\ell}^{(j)}(h g_k)w^{-1})
			\\
			&
			=
			w^{-1}
			\mathcal{C}_j(\Gamma_{K^\ell}^{(j)}(h g_k))
			w
			\\
			&
			=
			(h_\ell(g_k)_\ell w)^{-1}
			\mathcal{C}_j(\Gamma_{K^\ell}^{(j)}(1))
			(h_\ell(g_k)_\ell w).
		\end{align*}
		We write $b=h_\ell
		(g_k)_\ell
		w
		(g_i)_\ell^{-1}
		h_\ell^{-1}$. Then one deduces
		\[
		b
		\in
		\mathbf{H}_j(\Q)
		Z(\mathbf{H}_j( \Q_\ell)).
		\]
		Note that
		$h_\ell,(g_i)_\ell,(g_k)_\ell$
		normalize
		$\mathbf{U}_j$,
		so there is another
		$w'\in \mathbf{U}_j$
		such that
		$b=(g_k)_\ell(g_i)_\ell^{-1}w'$.
		Recall
		$\mathbf{T}(\Q_\ell)$
		normalizes
		$\mathbf{H}_j(\Q_\ell)$,
		thus
		for any
		$s\in \mathbf{T}(\Q_\ell)$
		and
		$z\in Z(\mathbf{H}_j(\Q_\ell))$,
		$szs^{-1}
		\in
		Z(\mathbf{H}_j(\Q_\ell))$.
		So we have a continuous morphism of topological groups
		\[
		\varphi
		\colon
		\mathbf{T}(\Q_\ell)
		\rightarrow
		\mathrm{Auto}(Z(\mathbf{H}_j(\Q_\ell)))
		\]
		induced by the conjugate action.
		Thus for any
		$s\in
		\mathrm{Ker}(\varphi)\cap
		\mathbf{T}(\Q)$
		(which is dense in $\mathrm{Ker}(\varphi)$),
		since
		$\mathbf{T}(\Q_\ell)$
		also normalizes
		$\mathbf{U}_j$,
		the commutator
		$sbs^{-1}b^{-1}
		\in
		\mathbf{H}(\Q)$
		is a unipotent element in
		$\mathbf{H}^1(\Q_\ell)$.
		Since
		$\mathbf{H}^1(\R)$
		is compact,
		$\mathbf{H}^1(\Q)$
		can not contain non-trivial unipotent elements,
		thus one must have
		$sbs^{-1}b^{-1}=1$
		for any
		$s\in
		\mathrm{Ker}(\varphi)
		\cap
		\mathbf{T}(\Q)$. Moreover, $\mathrm{Ker}(\varphi)
		\cap
		\mathbf{T}(\Q)$ is dense in
		$\mathrm{Ker}(\varphi)$,
		so
		$b$ commutes with all
		$\mathrm{Ker}(\varphi)$,
		an open subgroup of $\mathbf{T}(\Q_\ell)$.
		Thus
		$b$ commutes with $\mathbf{T}(\Q_\ell)$.
		This implies that
		$b\in
		\mathbf{T}( \Q_\ell)$
		since $\mathbf{T}(\Q_\ell)$ is a maximal torus in $\mathbf{H}(\Q_\ell)$. On the other hand, $b=(g_k)_\ell(g_i)_\ell^{-1}w'$.
		So we must have
		$w'=1$
		and therefore
		\[
		(g_k)_\ell(g_i)_\ell^{-1}
		\in
		\mathbf{T}(\Q_\ell)
		\cap
		\left(
		\mathbf{H}_j(\Q)Z(\mathbf{H}_j(\Q_\ell))
		\right).
		\]
		This contradicts our assumption
		(\ref{pairwise distinct-1}) on
		$g_i,g_k$.
		So
		$\Gamma_{K^\ell}^{(j)}(h g_i)$
		and
		$\Gamma_{K^\ell}^{(j)}(h g_k)$
		are
		\textit{not}
		$\mathbf{U}_j$-commensurable
		for any
		$i\ne k$.
		
		Now apply
		Theorem
		\ref{not V-commensurable}
		and we see
		that
		$\Gamma_{K^\ell}^{(j)}(\Delta(h)g)
		\Delta(\mathbf{U}_j)$
		is dense in
		$\mathbf{H}_j^{1}( \Q_\ell)^r$
		for any $j=1,\cdots,n$.
		Apply Theorem \ref{Ratner}(1)
		and we see that
		$\Gamma_{K^\ell}(\Delta(h)g)
		\Delta(\prod_j\mathbf{U}_j)$
		is dense in
		$\mathbf{H}^{1}( \Q_\ell)^r$
		since
		these
		$\mathbf{H}_j^{1}( \Q_\ell)^r$
		generate
		$\mathbf{H}^{1}( \Q_\ell)^r$.
	\end{proof}

	Analogue to 
	Theorem
	\ref{single-copy},
	we have
	\begin{theorem}\label{multi-copy}
		Fix $h\in\mathbf{T}(\mathbb{A}_f)$.
		For 
		any locally constant function $f$ on $(\mathbf{H}(\Q)\backslash \mathbf{H}(\mathbb{A}_f)/K^\ell)^r$
		and any compact open subset
		$\kappa$ of $ \Q_\ell$,
		we have
		\begin{equation}\label{equidistribution result-2}
			\lim\limits_{\underline{N}\rightarrow+\infty}
			\frac{1}{\mathbf{m}(\kappa_{\underline{N}})}
			\int_{\kappa_{\underline{N}}}
			f(\Delta(h)g\Delta(u(\underline{t})))
			d\underline{t}
			=
			\int_{\prod_i
				G(\Q)\backslash G(\Q)h
				g_i\mathbf{H}^{1}(\mathbb{A}_f)K^{\ell}/K^\ell}
			fd\mu_{\Delta(h)g}.
		\end{equation}
		Here
		$\mu_{\Delta(h)g}$
		is the product of
		the measures
		$\mu_{h g_i}$
		on
		$G(\Q)
		\backslash
		G(\Q)
		h
		g_i
		\mathbf{H}^{1}(\mathbb{A}_f)K^{\ell}/K^\ell$
		for
		$i=1,\cdots,r$.
	\end{theorem}

	\subsection{Main result}
	Let
	$\mathbf{H},\mathbf{T},\mathbf{H}_j,\mathbf{U}_j$ be as in
	§\ref{basic set-up}. We fix also the following data
	\begin{enumerate}
		\item
		$r$ elements $g_1,\cdots,g_r\in\mathbf{T}(\mathbb{A}_f)$
		such that
		(\ref{pairwise distinct-1})
		is satisfied. In particular, these $g_1,\cdots,g_r$ are all distinct. We write
		\[
		g=(g_1,g_2,\cdots,g_r).
		\]

		\item 
		a compact open subgroup $K^{\ell}$ of
		$\mathbf{H}(\mathbb{A}_f^{\ell})$.

		\item 
		an open subgroup $\mathcal{G}$ of
		$\mathbf{T}(\Q)\backslash
		\mathbf{T}(\mathbb{A}_f)$ and a Haar measure $\mu_{\mathcal{G}}$ on $\mathcal{G}$.
	\end{enumerate}

	Then integrating both sides (\ref{equidistribution result-2})
	with respect to the variable
	$h\in\mathcal{G}$,
	we get
	\begin{equation}\label{fundamental identity}
		\lim\limits_{\underline{N}\rightarrow+\infty}
		\frac{1}{\mathbf{m}(\kappa_{\underline{N}})}
		\int_{\kappa_{\underline{N}}}
		d\underline{t}
		\int_{\mathcal{G}}
		f(\Delta(h)g\Delta(u(\underline{t})))
		d\mu_{\mathcal{G}}(h)
		=
		\int_{\mathcal{G}}
		d\mu_{\mathcal{G}}(h)
		\int_{
			\prod_i
			\mathbf{H}(\Q)
			\backslash
			\mathbf{H}(\Q)h g_i\mathbf{H}^{1}(\mathbb{A}_f)K^\ell/K^\ell}
		fd\mu_{\Delta(h)g}.
	\end{equation}

	Under condition
	\ref{condition on mathbf{H}},
	we can choose $n$ torus subgroups
	$T_1\simeq \Q_\ell^\times,
	\cdots,
	T_n
	\simeq \Q_\ell^\times$
	of $\mathbf{T}( \Q_\ell)$
	whose conjugate action on
	$\prod_{j=1}^n\mathbf{U}_j$
	are pairwise distinct.

	\begin{proposition}
		Fix a compact open subgroup
		$K_\ell
		\subset
		\mathbf{H}( \Q_\ell)$.
		Then for any $m_1,\cdots,m_n\gg0$, we have
		\[
		\mathcal{G}
		\prod_{i=1}^nu_i
		\left(
		(1+\ell^{m_i}\Z_\ell)\ell^{-k_i}
		\right)
		K_\ell
		=
		\mathcal{G}
		\prod_{i=1}^n
		u_i(\ell^{-k_i})
		K_\ell,
		\quad
		\forall
		\,
		k_1,\cdots,k_n\ge0.
		\]
	\end{proposition}    
	\begin{proof}
		Recall we fixed isomorphisms
		\begin{align*}
			\sigma_i
			&
			\colon T_i\simeq\Q_\ell^\times,
			\\
			u_j
			&
			\colon\Q_\ell\simeq \mathbf{U}_j.
		\end{align*}
		thus
		there is an integer $r_{i,j}$
		such that the conjugate action
		$\varphi$
		of
		$T_i$ on
		$\mathbf{U}_j$
		is given by
		\[
		\varphi(\sigma_i(x_i),u_j(t))
		=
		u_j(x_i^{r_{i,j}}t),
		\quad
		\forall
		\,
		x_i\in \Q_\ell^\times,
		\,
		t\in   \Q_\ell.
		\]
		Under condition \ref{condition on mathbf{H}},
		the matrix
		$(r_{i,j})_{i,j=1}^n$
		is non-singular.
		For any
		$x_i\in   \Q_\ell^\times$
		viewed as an element in
		$T_i$
		and any
		$t_i\in   \Q_\ell$,
		one has
		\[
		\prod_{j=1}\sigma_j(x_j)
		\prod_{j=1}u_j(t_j\ell^{-k_j})
		\prod_{j=1}\sigma_j(x_j)^{-1}
		=
		\prod_{j=1}u_j\left(\prod_{i=1}^nx_i^{r_{i,j}}t_j\ell^{-k_j}\right).
		\]
		For any $m_1,\cdots,m_n\gg0$, we can choose $m_1',\cdots,m_n'>0$
		such that the system of equations on the variables
		$ x_1,\cdots, x_n$:
		\begin{equation*}
			\begin{cases*}
				x_1^{r_{1,1}}\cdots x_n^{r_{n,1}}=t_1^{-1}
				\in
				1+\ell^{m_1} \Z_\ell;
				\\
				x_1^{r_{1,2}}\cdots x_n^{r_{n,2}}=t_2^{-1}
				\in
				1+\ell^{m_2} \Z_\ell;
				\\
				\cdots
				\\
				x_1^{r_{1,n}}\cdots x_n^{r_{n,n}}=t_n^{-1}
				\in
				1+\ell^{m_n} \Z_\ell.
			\end{cases*}
		\end{equation*}
		always has a solution in
		$ x_i
		\in
		1+\ell^{m_i'} \Z_\ell$ ($i=1,\cdots,n$).
		We write $T_i(m_i'):=\sigma_i(1+\ell^{m_i'}\Z_\ell)$. Then we have the following identity
		\[
		\prod_{i=1}^n
		T_i(m_i')
		\prod_{i=1}^n
		u_i(t_i\ell^{-k_i})
		\prod_{i=1}^n
		T_i(m_i')
		=
		\prod_{i=1}^n
		T_i(m_i')
		\prod_{i=1}^n
		u_i(\ell^{-k_i})
		\prod_{i=1}^n
		T_i(m_i'),
		\quad
		\forall
		t_i\in
		1+\ell^{m_i} \Z_\ell.
		\]
		Moreover for $m_1,\cdots,m_n\gg0$, we can choose $m_1',\cdots,m_n'>0$ such that
		these compact open subgroups
		$T_i(m_i')$
		are all contained in
		$\mathcal{G}$
		and
		$K_\ell$ and their images in $\mathbf{T}(\Q)\backslash\mathbf{T}(\mathbb{A}_f)$ are all contained in $\mathcal{G}$.
		Then we get the identity as in the proposition.
	\end{proof}
	
	For simplicity, in the following we choose
	$m_1=m_2=\cdots=m_n=m\gg0$
	and write
	\[
	\kappa=1+\ell^m\Z_\ell.
	\]
	We also put
	\[
	K=K^\ell K_\ell
	\subset
	\mathbf{H}(\mathbb{A}_f).
	\]	
	For
	$\underline{N}=(N_1,\cdots,N_n)
	\in
	\mathbb{N}^n$,
	write
	\[
	u(\ell^{-\underline{N}})
	=
	u_1(\ell^{-N_1})\cdots u_n(\ell^{-N_n}).
	\]
	
	For any map
	$f\colon
	(\mathbf{H}(\Q)\backslash \mathbf{H}(\mathbb{A}_f)/K)^r
	\rightarrow
	\mathbb{C}$,
	define
	\[
	A(f,u(\underline{t}))
	:=
	\int_{\mathcal{G}}
	f
	\left(
	\Delta(h)g\Delta(u(\underline{t}))
	\right)
	d\mu_{\mathcal{G}}(h),
	\quad
	\forall
	\,
	\underline{t}=(t_1,\cdots,t_n)
	\in
	\Q_\ell^n.
	\]
	\begin{corollary}
		Fix $h\in\mathbf{T}(\mathbb{A}_f)$.
		We have
		\[
		A(f,u(\underline{t}))
		=
		\int_{\mathcal{G}}
		f
		\left(
		\Delta(h)g\Delta(u(\ell^{-\underline{N}}))
		\right)
		d\mu_{\mathcal{G}}(h),
		\quad
		\forall
		\,
		\underline{t}
		\in
		\kappa_{\underline{N}}.
		\]
		In other words, the function
		$A(f,u(\underline{t}))$
		is constant on the variable
		$\underline{t}$
		in the above specified domain
		$\kappa_{\underline{N}}$.
	\end{corollary}
	
	We can rewrite (\ref{fundamental identity}) as
	\[
	\lim\limits_{\underline{N}\rightarrow+\infty}
	A(f,u(\ell^{-\underline{N}}))
	=
	B(f)
	:=
	\text{RHS of (\ref{fundamental identity})}.
	\]

	We define the following objects\footnote{The terminology $\mathrm{SP}$ comes from the fact that if $(\mathbf{H},X)$ is a Shimura datum, then $\mathrm{SP}$ are exactly those \emph{special points} on the Shimura variety $Sh(\mathbf{H},X)_{\C}$ whose Mumford-Tate group is $\mathbf{T}$.}
	\begin{align*}
		\mathrm{SP}
		&
		:=
		\mathbf{T}(\Q)
		\backslash
		\mathbf{H}(\mathbb{A}_f),
		\\
		\mathcal{X}
		&
		:=
		\mathbf{H}(\Q)
		\backslash
		\mathbf{H}(\mathbb{A}_f),
		\\
		\mathcal{Z}
		&
		:=
		\mathbf{H}(\Q)
		\backslash
		\mathbf{H}(\mathbb{A}_f)/
		\mathbf{H}^{1}(\mathbb{A}_f).
	\end{align*}
	Moreover we have the following natural projection maps
	\[
	\mathrm{SP}
	\xrightarrow{\mathfrak{R}}
	\mathcal{X}
	\xrightarrow{\mathfrak{A}}
	\mathcal{Z}.
	\]
	Note that
	$\mathcal{X}$
	and
	$\mathcal{Z}$
	are both compact.
	Let
	$\mathbf{H}(\mathbb{A}_f)$
	act on the right on
	$\mathrm{SP},\mathcal{X},\mathcal{Z}$. Then these maps $\mathfrak{R}$ and $\mathfrak{A}$ are 
	$\mathbf{H}(\mathbb{A}_f)$-equivariant. Similarly we define the following objects and maps
	\begin{align*}
		\mathrm{SP}_K
		&
		:=
		\mathrm{SP}/K,
		\\
		\mathcal{X}_K
		&
		:=
		\mathcal{X}/K,
		\\
		\mathcal{Z}_K
		&
		:=
		\mathcal{Z}/K.
	\end{align*}
	\[
	\mathrm{SP}_K
	\xrightarrow{\mathfrak{R}_K}
	\mathcal{X}_K
	\xrightarrow{\mathfrak{A}_K}
	\mathcal{Z}_K.
	\]
	Here
	$\mathcal{X}_K,\mathcal{Z}_K$
	are finite sets.
	We define furthermore the following maps:
	\begin{align*}
		\mathfrak{R}_K^r
		\colon
		&    	
		\mathrm{SP}_K
		\rightarrow
		\mathcal{X}_K^r,
		\quad
		x
		\mapsto
		(\mathfrak{R}_K(g_ix))_{i=1}^r,
		\\
		\mathfrak{A}_K^r
		\colon
		&
		\mathcal{X}_K^r
		\rightarrow
		\mathcal{Z}_K^r,
		\quad
		(y_i)_{i=1}^r
		\mapsto
		(\mathfrak{A}_K(y_i))_{i=1}^r,
		\\
		\mathfrak{A}^r
		\colon
		&
		\mathcal{X}^r
		\rightarrow
		\mathcal{Z}^r,
		\quad
		(y_i)_{i=1}^r
		\mapsto
		(\mathfrak{A}(y_i))_{i=1}^r.
	\end{align*}

	Since
	$\mathbf{H}^1(\mathbb{A}_f)$
	acts transitively on each fiber of the map
	$\mathfrak{A}$,
	for any
	$z\in\mathcal{Z}$,
	there is a unique Borel measure
	$\mu_{z}$
	on
	$\mathfrak{A}^{-1}(z)$
	invariant under
	$\mathbf{H}^1(\mathbb{A}_f)$.
	
	Here is the main result of this appendix:
	\begin{theorem}\label{density of torus orbits}
		For
		$N_1,\cdots,N_n\gg0$,
		one has
		\[
		\mathfrak{R}_K^r(\mathcal{G}u(\ell^{-\underline{N}}))
		=
		(\mathfrak{A}_K^r)^{-1}
		\left(
		\mathfrak{A}_K^r\circ
		\mathfrak{R}_K^r(\mathcal{G}u(\ell^{-\underline{N}}))
		\right).
		\]
	\end{theorem}
	\begin{proof}
		We prove that RHS is contained in LHS, the other direction being trivial.
		Let
		$s
		\in
		(\mathfrak{A}_K^r)^{-1}
		\left(
		\mathfrak{A}_K^r
		(\mathfrak{R}_K^r(\mathcal{G}u(\ell^{-\underline{N}})))
		\right)$
		and
		$\mathbb{I}_s
		\colon
		\mathcal{X}^r
		\rightarrow
		\mathbb{C}$
		be the characteristic function
		of the pre-image of $s$ under the projection map
		$\mathcal{X}^r
		\rightarrow
		\mathcal{X}^r_K$.
		We compute both sides of
		(\ref{fundamental identity})
		for $f=\mathbb{I}_s$.
		It is easy to see that $A(\mathbb{I}_s,u(\ell^{-\underline{N}}))$ is equal to the measure of the subset $\{
		h\in\mathcal{G}|		
		\mathfrak{R}_K^r(h u(\ell^{-\underline{N}}))
		=s
		\}$ of $\mathcal{G}$:
		\[
		A(\mathbb{I}_s,u(\ell^{-\underline{N}}))
		=
		\mu_{\mathcal{G}}
		\{
		h\in\mathcal{G}|		
		\mathfrak{R}_K^r(h u(\ell^{-\underline{N}}))
		=s
		\}.
		\]

		For any $z=(z_1,\cdots,z_r)\in\mathcal{Z}^r$, write $\mu_z$
		for the product of the measures
		$\mu_{z_i}$
		on
		$\mathfrak{A}^{-1}(z_i)$ ($i=1,\cdots,r$). Then we define
		\[
		I(\mathbb{I}_s,z)
		:=
		\int_{(\mathfrak{A}^r)^{-1}(z)}
		\mathbb{I}_sd\mu_{z}.
		\]

		It is easy to see that
		$I(\mathbb{I}_s,\cdot)$
		factors through
		the quotient
		$\mathcal{Z}^r
		\rightarrow
		\mathcal{Z}_K^r$.
		Write
		$\Omega(\mathcal{G})$
		for the common cardinal of $\mathcal{G}$-orbits
		on
		$\mathcal{Z}_K^r$.
		Then we have
		\[
		B(\mathbb{I}_s)
		=
		\frac{I(\mathbb{I}_s,\mathfrak{A}_K^r(s))}{\Omega(\mathcal{G})}.
		\]
		
		Lemma \ref{I>0}
		shows that
		$I(\mathbb{I}_s,\mathfrak{A}_K^r(s))
		\ne0$.
		Moreover
		$\mathcal{X}_K^r$
		is a finite set,
		thus
		there are only finitely many values for
		$I(\mathbb{I}_s,\mathfrak{A}_K^r(s))$
		(none of which is zero).
		So for 
		$N_1,\cdots,N_n\gg0$,
		we have
		\[
		|A(\mathbb{I}_s,u(\ell^{-\underline{N}}))-B(\mathbb{I}_s)|\le|B(\mathbb{I}_s)|/2.
		\]
		In particular,
		$A(\mathbb{I}_s,u(\ell^{-\underline{N}}))\ne0$,
		so
		for any
		$s
		\in
		(\mathfrak{A}_K^r)^{-1}
		\left(
		\mathfrak{A}_K^r
		(\mathfrak{R}_K^r(\mathcal{G}u(\ell^{-\underline{N}})))
		\right)$,
		there exists
		$h\in\mathcal{G}$
		such that
		$\mathfrak{R}_K^r
		(h u(\ell^{-\underline{N}}))=s$.
	\end{proof}
	\begin{lemma}\label{I>0}
		For
		$s\in
		(\mathfrak{A}_K^r)^{-1}
		\left(
		\mathfrak{A}_K^r
		(\mathfrak{R}_K^r(\mathcal{G}u(\ell^{-\underline{N}})))
		\right)$,
		we have
		\[
		I(\mathbb{I}_s,\mathfrak{A}_K^r(s))>0.
		\]
	\end{lemma}
	\begin{proof}
		We follow the proof in
		\cite[Proposition 2.14]{CornutVatsal2005}.

		We write $z_i=\mathfrak{A}_K(s_i)$ and set $\mathbb{I}_{i}
		\colon
		\mathcal{X}
		\rightarrow
		\{0,1\}$
		the characteristic function of
		the pre-image of
		$s_i$ under the projection map
		$\mathcal{X}
		\rightarrow
		\mathcal{X}_K$. Then we put
		\[
		I(s_i)=\int_{\mathfrak{A}^{-1}(z_i)}\mathbb{I}_id\mu_{z_i}.
		\]
		So one has
		\[
		I(s):=I(\mathbb{I}_s,\mathfrak{A}_K^r(s))=\prod_{i=1}^rI(s_i).
		\]
		Therefore to prove the lemma, it suffices to show that
		$I(s_i)>0$ for any $i=1,\cdots,n$.

		Let
		$z\in\mathcal{Z}$.
		The normalized Haar measure
		$\mu$
		on
		$\mathbf{H}^1(\mathbb{A}_f)$
		induces
		the measure
		$\mu_z$ on
		$\mathfrak{A}^{-1}(z)\subset\mathcal{X}$.
		Note that
		$\mu_z$ is the unique measure on
		$\mathfrak{A}^{-1}(z)$
		such that for any compact open subgroup
		$\mathbf{K}$ of
		$\mathbf{H}^1(\mathbb{A}_f)$
		and any
		$x\in\mathfrak{A}^{-1}(z)$,
		one has
		\[
		\mu_z(x\mathbf{K})
		=
		\frac{\mu(\mathbf{K})}{\sharp\mathbf{K}_x}.
		\]
		Note that $\mathbf{K}_x$
		is compact and discrete, thus is a finite set. So $\mu_z(x\mathbf{K})>0$.
		Moreover for any
		$g\in \mathbf{H}^1(\mathbb{A}_f)$,
		the measure
		$\mu_{zg}(\cdot g)=(\mu_{zg}\circ g)(\cdot)$ on
		$\mathfrak{A}^{-1}(z)$
		is equal to
		$\mu_z$.
		In other words,
		for any
		$z_1,z_2\in\mathcal{Z}$,
		$\mu_{z_1}(\mathfrak{A}^{-1}(z_1))
		=
		\mu_{z_2}(\mathfrak{A}^{-1}(z_2))$.
		Now for any
		$x\in\mathcal{X}$
		and
		$z\in\mathcal{Z}$,
		we write
		\[
		\phi_z(x)
		=
		\mu_z(xK\cap\mathfrak{A}^{-1}(z)).
		\]
		Then the map
		$\mathcal{Z}
		\rightarrow
		\mathbb{C}$
		sending $z$ to
		$\phi_z(x)$
		factors through
		$\mathcal{Z}\rightarrow\mathcal{Z}_K$.
		Similarly
		the map
		$\mathcal{X}
		\rightarrow
		\mathbb{C}$
		sending
		$x$ to $\phi_z(x)$
		factors through
		$\mathcal{X}
		\rightarrow
		\mathcal{X}_K$.
		By definition we have
		\[
		\phi_z(x)
		\begin{cases*}
			=0,
			&
			$z\notin\mathfrak{A}(xK)$;
			\\
			>0,
			&
			$z\in\mathfrak{A}(xK)$.
		\end{cases*}
		\]
		Thus for any $i=1,\cdots,r$, we have
		\[
		I(s_i)
		=
		\phi_{z_i}(s_i)
		>0.
		\]
		
	\end{proof}

	\subsection{Application to automorphic forms}\label{applications to auto forms}
	Let
	$\mathbf{H},\mathbf{T},\mathbf{H}_j,\mathbf{U}_j$
	be as in
	§\ref{basic set-up}.
	Let $K$
	be a compact open subgroup of
	$\mathbf{H}(\mathbb{A}_f)$.
	We fix a ring
	$A$,
	an $A$-module $M$.
	Let
	$g=(g_1,\cdots,g_r)\in\mathbf{T}(\mathbb{A}_f)$ satisfy
	(\ref{pairwise distinct-1}).
	Let
	$\mathcal{G}$
	be an open subgroup of
	$\mathbf{T}(\Q)\backslash \mathbf{T}(\mathbb{A}_f)$
	and write $\widetilde{\mathcal{G}}$
	its pre-image by the projection map
	$\mathbf{T}(\mathbb{A}_f)
	\rightarrow
	\mathbf{T}(\Q)\backslash \mathbf{T}(\mathbb{A}_f)$.
	Then a similar argument as in
	\cite[Cor.5.2]{ChidaHsieh2016}
	gives
	\begin{theorem}\label{automorphic application of appendix}
		Let
		$\{\beta_i\}_{i=1}^r$
		be a finite set of elements in $A$ with
		$\beta_1\in A^\times$.
		We assume that
		the composition map
		\[
		\widetilde{\mathcal{G}}
		\rightarrow
		\mathbf{H}(\mathbb{A}_f)
		\rightarrow
		\mathcal{Z}_K
		\]
		is surjective.
		Consider a map
		$f\colon\mathbf{H}(\Q)\backslash\mathbf{H}(\mathbb{A}_f)/K\rightarrow M$ which is not $\mathbf{H}^1(\mathbb{A}_f)$-invariant (under right translation).	
		Then for any
		$N_1,\cdots,N_n\gg0$,
		there is an element
		$h=h_{\underline{N}}
		\in
		\widetilde{\mathcal{G}}$
		such that
		\[
		\sum_{i=1}^r
		\beta_if(h g_iu(\ell^{-\underline{N}}))
		\ne0.
		\]
	\end{theorem}
	\begin{proof}
		By Theorem \ref{density of torus orbits},
		for any
		$N_1,\cdots,N_n\gg0$,
		\[
		\mathfrak{R}_K^r(\mathcal{G}u(\ell^{-\underline{N}}))
		=
		(\mathfrak{A}_K^r)^{-1}
		\left(
		\mathfrak{A}_K^r
		(
		\mathfrak{R}_K^r(\mathcal{G}u(\ell^{-\underline{N}}))
		)
		\right).
		\]
		We fix one such
		$\underline{N}=(N_1,\cdots,N_n)$.
		By assumption,
		there are elements
		$y_1\ne y_2\in\mathcal{X}$
		such that
		$\mathfrak{A}_K(y_1)=\mathfrak{A}_K(y_2)$
		and
		$f(y_1)\ne f(y_2)$.
		Since
		$\widetilde{\mathcal{G}}
		\rightarrow
		\mathcal{Z}_K$
		is surjective,
		we can choose
		$x\in\widetilde{\mathcal{G}}$
		such that
		\[
		\mathfrak{A}_K(y_1)=\mathfrak{A}_K(y_2)
		=
		\mathfrak{A}_K(\mathfrak{R}_K(xu(\ell^{-\underline{N}}))).
		\]
		Therefore for any
		finite subset
		$\{x_i\}_{i=1}^r$
		of
		$\mathfrak{A}_K^{-1}(\mathfrak{A}_K(y_1))$,
		we can find
		$h_1,h_2
		\in \widetilde{\mathcal{G}}$
		such that
		\[
		\mathfrak{R}_K^r(h_ixu(\ell^{-\underline{N}}))
		=
		(y_i,x_2,\cdots,x_r),
		\quad
		i=1,2.
		\]
		Thus
		one has
		\[
		\sum_{i=1}^r\beta_if(h_1g_ixu(\ell^{-\underline{N}}))
		-
		\sum_{i=1}^r\beta_if(h_2g_ixu(\ell^{-\underline{N}}))
		=
		\beta_1(f(y_1)-f(y_2))
		\ne0.
		\]
		So we can take
		$h=h_{\underline{N}}$ to be
		$h_1x$ or $h_2x$.    	
	\end{proof}

	\subsection{Examples}\label{Example}
	In this section we give some examples of $(\mathbf{H},\mathbf{T},\mathbf{H}_j,\mathbf{U}_j)$
	satisfying condition
	\ref{condition on mathbf{H}}.
	We give these explicit examples with the purpose in mind that
	they may be used directly in the theory of theta lifts.
	
	\subsubsection{Unitary groups}
	Let
	$K$
	be an imaginary quadratic number field
	with an embedding
	$K\hookrightarrow
	\Q_\ell$.
	Let
	$V=K^{n+1}$
	be a vector space over $K$ of dimension $n+1\ge2$
	equipped with a Hermitian form
	$Q$,
	which is represented,
	under the standard $K$-basis
	$\{E_1,\cdots,E_{n+1}\}$
	of $V$,
	by a diagonal matrix
	$Q=\mathrm{diag}(\delta_1,\cdots,\delta_{n+1})$.
	Suppose that $\delta_1,\cdots,\delta_{n+1}>0$.
	We take
	$\mathbf{H}=\mathrm{U}(V,Q)$
	the
	unitary group associated to $(V,Q)$.
	Then
	$\mathbf{H}_{/K}
	\simeq
	\mathrm{GL}_{n+1/K}$
	as algebraic groups over $K$.

	For any distinct basis elements
	$E_i,E_j$,
	we write
	$U_{i,j}$
	for the
	unitary group associated to
	the Hermitian subspace
	$(K(E_i,E_j),Q)$
	of
	$(V,Q)$.
	Now we put
	\[
	\mathbf{H}_j=U_{1,j+1},
	\quad
	j=1,\cdots,n.
	\]
	We fix an isomorphism
	\[\iota\colon
	\mathrm{GL}_{n+1}( \Q_\ell)
	\simeq
	\mathbf{H}( \Q_\ell),
	\]
	compatible with
	$\mathrm{GL}_{n+1}(K)\simeq\mathbf{H}(K)$
	under the embedding
	$K\hookrightarrow \Q_\ell$
	such that
	$\mathbf{H}_j( \Q_\ell)$
	is mapped isomorphically to
	the subgroup of
	$\mathrm{GL}_{n+1}( \Q_\ell)$
	consisting of matrices of the following form
	\[
	\begin{pmatrix}
		a & & b & \\
		& 1_{j-1} & & \\
		c & & d & \\
		& & & 1_{n-j}
	\end{pmatrix}
	\quad
	\text{with }
	\begin{pmatrix}
		a & b \\ c & d
	\end{pmatrix}\in\mathrm{GL}_2(\Q_\ell).
	\]
	Fix isomorphisms
	\[
	\iota_j
	\colon
	\mathrm{GL}_2( \Q_\ell)
	\simeq
	\mathbf{H}_j(\Q_\ell)
	\]
	such that
	$(\iota^{-1}\circ\iota_j)(\mathrm{U}_2( \Q_\ell))$ consists of unipotent upper triangular matrices.
	We then take
	the unipotent subgroup
	$\mathbf{U}_j\subset \mathbf{H}_j(\Q_\ell)$
	to be the image of $\mathrm{U}_2(\Q_\ell)$
	by $\iota_j$.
	It is easy to see that these groups
	$\mathbf{H}_j^{1}( \Q_\ell)$
	generate
	$\mathbf{H}^{1}( \Q_\ell)$.
	Finally we take
	\[
	\mathbf{T}=\prod_{i=1}^{n+1}
	\mathrm{U}(KE_i,Q).
	\]
	One verifies that
	$(\mathbf{H},\mathbf{T},\mathbf{H}_1,\cdots,\mathbf{H}_n,\mathbf{U}_1,\cdots,\mathbf{U}_n)$
	satisfies condition
	\ref{condition on mathbf{H}}.

	\subsubsection{Even spin groups}\label{Even spin groups}
	Suppose that $-1$ has a square root in $\Q_\ell$, which we denote by $\mathbf{i}$.

	Let
	$(V,Q)$
	be a vector space over $\Q$
	of even dimension
	$2n\ge4$ such that under the standard basis
	$\{E_1,E_2,\cdots,E_{2n}\}$
	of $V$,
	$Q$ is represented by a diagonal matrix
	$Q=\mathrm{diag}(\delta_1,\cdots,\delta_{2n})$. In the following, for $1\leq i_1<i_2<\cdots<i_k\leq2n$, we write
	\begin{align*}
		\delta_{i_1,\cdots,i_k}
		&
		=\delta_{i_1}\cdots\delta_{i_k},
		\\
		E_{i_1,\cdots,i_k}
		&
		=E_{i_1}\otimes\cdots\otimes E_{i_k}.
	\end{align*}
	We put the following conditions on $\delta_1,\cdots,\delta_{2n}$:
	\begin{enumerate}
		\item 
		$\delta_1,\cdots,\delta_{2n}>0$,
		
		\item 
		$\delta_1,\cdots,\delta_{2n}
		\in
		(\Q_\ell^\times)^2$,
		
		\item 
		$\delta_{1,2,3,4},
		\delta_{1,2,5,6},
		\cdots,
		\delta_{1,2,2n-1,2n}
		\in
		(\Q^\times)^2$.
	\end{enumerate}
	We fix $h_i\in\Q_\ell^\times$ ($i=1,\cdots,2n$) such that
	\[
	h_i^2=\delta_i^{-1}.
	\]
	Then we take
	$\mathbf{H}=\mathrm{GSpin}(V,Q)$.
	More precisely,
	write
	$C(V)$
	for the even degree part of the Clifford algebra
	\[
	\mathrm{Cliff}(V)=\bigoplus_{k=0}^\infty
	V^{\otimes k}/
	\langle
	v\otimes v-Q(v)|v\in V
	\rangle
	\]
	associated to $(V,Q)$.
	Then
	$\mathbf{H}$
	consists of those units
	$v\in C(V)^\times$
	such that
	$vV=Vv$.
	For any distinct basis elements
	$E_{i_1},E_{i_2},\cdots,E_{i_k}$,
	we write
	$C_{i_1,\cdots,i_k}$
	for
	the even degree part of the Clifford algebra associated to
	the quadratic subspace
	$(\Q(E_{i_1},\cdots,E_{i_k}),Q)$
	of $(V,Q)$.
	Then we have an isomorphism of
	$\Q_\ell$-algebras:
	\begin{align*}		    	
		\mathrm{M}_{2}( \Q_\ell)
		&
		\simeq
		C_{i,j,k}\otimes_\Q \Q_\ell,
		\\
		\begin{pmatrix}
			0 & 1 \\
			0 & 0
		\end{pmatrix}
		&
		\mapsto
		\frac{(h_iE_i+h_j\mathbf{i}E_j)(h_kE_k)}{2}
		=:
		e_{i,j,k}^+,
		\\
		\begin{pmatrix}
			0 & 0 \\
			1 & 0
		\end{pmatrix}
		&
		\mapsto
		-\frac{(h_iE_i-h_j\mathbf{i}E_j)(h_kE_k)}{2}
		=:
		e_{i,j,k}^-,
		\\
		\begin{pmatrix}
			a & 0 \\
			0 & 1
		\end{pmatrix}
		&
		\mapsto
		1+\frac{a-1}{4}
		(h_iE_i+h_j\mathbf{i}E_j)(h_i-h_j\mathbf{i}E_j),
		\quad
		a\in\Q_\ell,
		\\
		\begin{pmatrix}
			a & 0 \\
			0 & \frac{1}{a}
		\end{pmatrix}
		&
		\mapsto
		\frac{1}{a}
		\left(
		1+
		\frac{a^2-1}{4}
		(h_iE_i+h_j\mathbf{i}E_j)
		(h_iE_i-h_j\mathbf{i}E_j)
		\right)
		=:
		\tau_{i,j}(a),\quad a\in\Q_\ell^\times.
	\end{align*}
	We have the following identities
	\begin{align*}
		(1+ve_{i_1,i_2,k}^+)(1\pm ve_{i_1,i_2,l}^+)
		&
		=
		1+v(e_{i_1,i_2,k}^+\pm e_{i_1,i_2,l}^+),
		\quad
		i_1,i_2,k,l
		\in
		\{
		1,\cdots,2n
		\} \text{ distinct},
		\\
		\tau_{i_1,i_2}(a)(1+ve_{i_1,i_2,k}^\pm)\tau_{i_1,i_2}(a)^{-1}
		&
		=
		1+a^{\pm2}ve_{i_1,i_2,k}^\pm,
		\quad
		i_1,i_2,k\in
		\{1,\cdots,2n\} \text{ distinct}.
	\end{align*}
	For
	$k=2,3,\cdots,n$,
	we fix a positive square root
	$\sqrt{\delta_{1,2,2k-1,2k}}$
	of
	$\delta_{1,2,2k-1,2k}$
	in
	$\Q$
	and
	we write
	\[
	\mathbf{e}_k
	=
	\frac{1}{2}
	\left(
	1
	-
	\frac{E_{1,2,2k-1,2k}}{\sqrt{\delta_{1,2,2k-1,2k}}}
	\right),
	\]
	\sloppy
	which is a central idempotent element in
	$C_{1,2,2k-1,2k}$.
	One verifies that both
	$C_{1,2,2k-1,2k}\mathbf{e}_k$
	and
	$C_{1,2,2k-1,2k}(1-\mathbf{e}_k)$
	are central simple algebras over $\Q$ of dimension $4$.
	Then we put
	\[
	\mathbf{H}_j
	=
	\begin{cases*}
		(C_{1,2,3,4}(1-\mathbf{e}_2))^\times
		+
		\mathbf{e}_2,
		&
		$j=1$;
		\\
		(C_{1,2,2j-1,2j}\mathbf{e}_j)^\times
		+
		(1-\mathbf{e}_j),
		&
		$j=2,\cdots,n$.
	\end{cases*}
	\]
	The unipotent subgroups $\mathbf{U}_j$
	are given as follows:
	\[
	\mathbf{U}_j
	=
	\begin{cases*}
		(1+\Q_\ell e_{1,2,3}^+)(1-\mathbf{e}_2)
		+\mathbf{e}_2
		=
		1+\Q_\ell(e_{1,2,3}^+-\mathbf{i}e_{1,2,4}^+),
		&
		$j=1$;
		\\
		(1+\Q_\ell e_{1,2,2j-1}^+)\mathbf{e}_j
		+
		(1-\mathbf{e}_j)
		=
		1+\Q_\ell(e_{1,2,2j-1}^++\mathbf{i}e_{1,2,2j}^+),
		&
		$j=2,\cdots,n$.
	\end{cases*}
	\]
	Moreover it is easy to see
	these
	$\mathbf{U}_j$ commute with each other.
	The maximal torus
	is given by
	\[
	\mathbf{T}
	=
	\prod_{k=1}^n
	C_{2k-1,2k}^\times.
	\]
	Thus
	the conjugate
	action of
	$\tau_{2i-1,2i}(\Q_\ell)$
	on
	$\mathbf{U}_k$
	is given as follows
	\begin{align*}
		\mathrm{Ad}_{\tau_{2i-1,2i}(a)}
		\left(
		(1+ve_{1,2,3}^+)(1-\mathbf{e}_2)+\mathbf{e}_2
		\right)
		&
		=
		\begin{cases*}
			(1+a^2ve_{1,2,3}^+)(1-\mathbf{e}_2)+\mathbf{e}_2,
			&
			$i\in\{1,2\}$;
			\\
			(1+ve_{1,2,3}^+)(1-\mathbf{e}_2)+\mathbf{e}_2,
			&
			$i\in\{3,\cdots,n\}$.
		\end{cases*}
		\\
		\mathrm{Ad}_{\tau_{1,2}(a)}
		\left(
		(1+ve_{1,2,2k-1}^+)\mathbf{e}_k+(1-\mathbf{e}_k)
		\right)
		&
		=
		(1+a^2ve_{1,2,2k-1}^+)\mathbf{e}_k+(1-\mathbf{e}_k),
		\quad
		k\in\{2,\cdots,n\}.
		\\
		\mathrm{Ad}_{\tau_{2i-1,2i}(a)}
		\left(
		(1+ve_{1,2,2k-1}^+)\mathbf{e}_k+(1-\mathbf{e}_k)
		\right)
		&
		=
		\begin{cases*}
			(1+a^{-2}ve_{1,2,2k-1}^+)\mathbf{e}_k+(1-\mathbf{e}_k),
			&
			$i=k\in\{2,\cdots,n\}$;
			\\
			(1+ve_{1,2,2k-1}^+)\mathbf{e}_k+(1-\mathbf{e}_k),
			&
			$i\ne k\in\{2,\cdots,n\}$.
		\end{cases*}
	\end{align*}
	One verifies easily that
	these $\mathbf{H}_1^{1}( \Q_\ell),\cdots,\mathbf{H}_n^1(\Q_\ell)$
	generate
	$\mathbf{H}^{1}( \Q_\ell)$:
	indeed,
	the above formula shows that
	the opposite root subgroups
	$\mathbf{U}_1^-,\cdots,\mathbf{U}_n^-$
	of
	$\mathbf{U}_1,\cdots,\mathbf{U}_n$
	(with respect to $\mathbf{T}(\Q_\ell)$)
	are contained in
	the subgroup of $\mathbf{H}^1(\Q_\ell)$
	generated by
	$\mathbf{H}_1^1(\Q_\ell),\cdots,\mathbf{H}_n^1(\Q_\ell)$.
	Since
	$\mathbf{H}^1(\Q_\ell)$
	is simply-connected,
	we deduce that
	$\mathbf{H}_1^1(\Q_\ell),\cdots,\mathbf{H}_n^1(\Q_\ell)$
	generate $\mathbf{H}^1(\Q_\ell)$.
	One checks easily that the remaining part of
	condition
	\ref{condition on mathbf{H}} are also satisfied.

	\subsubsection{Odd spin groups}\label{Odd spin groups}
	As in §\ref{Even spin groups}, we suppose that $-1$ has a square root in $\Q_\ell$, which we denote by $\mathbf{i}$.

	Let
	$(V,Q)$
	be a vector space of odd dimension $2n+1\ge3$ such that under the standard basis
	$\{E_0,\cdots,E_{2n}\}$
	of $V$, the quadratic form
	$Q$ is represented by the diagonal matrix
	$\mathrm{diag}(\delta_0,\delta_1,\cdots,\delta_{2n})$.

	As in the preceding example, we put the following conditions on $\delta_0,\delta_1,\cdots,\delta_{2n}$:
	\begin{enumerate}
		\item 
		$\delta_0,\cdots,\delta_{2n}>0$,

		\item 
		$\delta_0,\cdots,\delta_{2n}
		\in
		(\Q_\ell^\times)^2$,
		
		\item 
		$\delta_{1,2,3,4},
		\delta_{1,2,5,6},
		\cdots,
		\delta_{1,2,2n-1,2n}
		\in
		(\Q^\times)^2$.
	\end{enumerate}
	Using the notations from the preceding example,
	we take
	$\mathbf{H}
	=
	\mathrm{GSpin}(V,Q)$
	and
	\begin{align*}
		\mathbf{H}_j
		=
		\begin{cases*}
			C_{0,1,2}^\times,
			&
			$j=1$;
			\\
			(C_{1,2,2j-1,2j}\mathbf{e}_i)^\times+(1-\mathbf{e}_j),
			&
			$j=2,\cdots,n$.
		\end{cases*}
	\end{align*}
	The unipotent subgroups $\mathbf{U}_j$
	are given as follows:
	\begin{align*}
		\mathbf{U}_j
		=
		\begin{cases*}
			1+\Q_\ell e_{1,2,0}^+,
			&
			$j=1$;
			\\
			(1+\Q_\ell e_{1,2,2j-1}^+)\mathbf{e}_j,
			&
			$j=2,\cdots,n$.
		\end{cases*}
	\end{align*}
	One checks that these
	$\mathbf{U}_j$ commute with each other.
	The maximal torus is given by
	\[
	\mathbf{T}
	=
	\prod_{k=1}^n
	C_{2k-1,2k}^\times.
	\]
	Thus
	the conjugate
	action of
	$\tau_{2i-1,2i}(\Q_\ell)$
	on
	$\mathbf{U}_1$
	is given as follows
	\[
	\mathrm{Ad}_{\tau_{2i-1,2i}(a)}(1+ve_{1,2,0}^+)
	=
	\begin{cases*}
		1+a^2ve_{1,2,0}^+,
		&
		$i=1$;
		\\
		1+ve_{1,2,0}^+
		&
		$i\in\{2,\cdots,n\}$.
	\end{cases*}
	\]
	and their actions on
	$\mathbf{U}_2,\cdots,\mathbf{U}_n$
	are the same as in the preceding example.
	The same reasoning as above shows that
	$\mathbf{H}_1^{1}( \Q_\ell),
	\cdots,
	\mathbf{H}_n^{1}( \Q_\ell)$
	generate
	$\mathbf{H}^{1}( \Q_\ell)$.

	\subsubsection{Symplectic groups}
	We write $J_1=
	\begin{pmatrix}
		0 & 1 \\
		-1 & 0
	\end{pmatrix}$ and for a positive integer $n$, we write
	\[
	J_{2n}=\mathrm{diag}(J_1,\cdots,J_1)
	\in
	\mathrm{GL}_{2n}.
	\]
	We define
	$\mathrm{GSp}_{2n}$ to be the group scheme over $\mathbb{Z}$
	consisting of matrices
	$X\in\mathrm{GL}_{2n}$
	such that
	$XJ_{2n}X^\mathrm{t}
	=
	\mu(X)
	J_{2n}$ with $\mu(X)\in\mathbb{G}_m$.

	Fix an integer $n>1$ and $\mathbf{B}$ a definite quaternion algebra over $\Q$,
	which is split at the place $\ell$ of $\Q$
	as in
	Example
	\ref{examples}
	and $n>1$ an integer.
	Let
	$\ast
	\colon
	\mathbf{B}
	\rightarrow
	\mathbf{B}$
	be a main involution and
	for an $n\times n$-matrix
	$g=(g_{i,j})_{i,j=1}^n
	$
	with entries in
	$\mathbf{B}$,
	write
	$g^\ast$
	to be
	$(g_{j,i}^\ast)_{i,j=1}^n$.
	Then we define a quaternionic unitary group
	\[
	\mathbf{H}
	=
	\{
	g\in\mathrm{GL}_n(B)|
	gg^\ast=\mu(g)\cdot1_n,
	\,
	\text{for some }
	\mu(g)\in \Q^\times
	\}.
	\]
	We view $\mathbb{H}$ as an algebraic group over $\Q$.
	Since $\mathbf{B}$ is split at $\ell$, we have an isomorphism of $\Q_\ell$-algebras
	\[
	\iota\colon
	\mathrm{M}_2(\Q_\ell)
	\simeq
	\mathbf{B}\otimes\Q_\ell,
	\]
	which induces an isomorphism
	\[
	\iota\colon
	\mathrm{GSp}_{2n}(\Q_\ell)
	\simeq
	\mathbf{H}(\Q_\ell).
	\]
	
	We fix then a maximal commutative subalgebra
	$K$ of $\mathbf{B}$,
	which is a quadratic field extension of $\Q$. We fix a $\Q$-basis of $\mathbf{B}$
	\[
	\mathbf{B}=\Q(1,\mathbf{i},\mathbf{j},\mathbf{k}).
	\]
	We assume that
	$\mathbf{k}=\mathbf{i}\mathbf{j}=-\mathbf{j}\mathbf{i}$,
	$\mathbf{i}^2,\mathbf{j}^2,\mathbf{k}^2\in\Q^\times$ and
	$K=\Q[\mathbf{i}]$ splits at $\ell$\footnote{Note that here $\mathbf{i}$ is not necessarily a square root of $-1$, as in the main body of this article}.	
	We fix embeddings $\iota_j\colon\mathbf{B}^\times\hookrightarrow\mathbf{H}$ as follows ($j=1,\cdots,n$):
	\[
	\iota_j
	\colon
	\mathbf{B}^\times
	\rightarrow
	\mathbf{H},
	\quad
	x=
	a+b\mathbf{i}+c\mathbf{j}+d\mathbf{k}
	\mapsto
	\begin{cases*}
		\begin{pmatrix}
			x & \\
			& 1_{n-1}
		\end{pmatrix},
		&
		$j=1$;
		\\
		\begin{pmatrix}
			a+b\mathbf{i} & & c\mathbf{j}+d\mathbf{k} & \\
			& 1_{j-2} & & \\
			-c\mathbf{j}-d\mathbf{k} & & a+b\mathbf{i} & \\
			& & & 1_{n-j}
		\end{pmatrix},
		&
		$j=2,\cdots,n$.
	\end{cases*}
	\]
	Then set
	$\mathbf{H}_j=\iota_j(\mathbf{B}^\times)$.
	We fix an isomorphism of $\Q_\ell$-algebras
	\[
	\mathbf{B}
	\otimes
	\Q_\ell
	\rightarrow
	\mathrm{M}_2(\Q_\ell),
	\quad
	\mathbf{i}\otimes1
	\mapsto
	\begin{pmatrix}
		\mathbf{i} & 0 \\
		0 & -\mathbf{i}
	\end{pmatrix},
	\,
	\mathbf{j}\otimes1
	\mapsto
	\begin{pmatrix}
		0 & \mathbf{j}^2 \\
		1  & 0
	\end{pmatrix},
	\,
	\mathbf{k}\otimes1
	\mapsto
	\begin{pmatrix}
		0 & \mathbf{i}\mathbf{j}^2 \\
		-\mathbf{i} & 0
	\end{pmatrix}.
	\]
	The unipotent subgroups
	$\mathbf{U}_j$ are given by
	\[
	\mathbf{U}_j=(\iota\circ\iota_j\circ\iota^{-1})(\mathrm{U}_2).
	\]
	Clearly these
	unipotent subgroups commute with each other.
	The maximal torus $\mathbf{T}$ of $\mathbf{H}$ is given by
	\[
	\mathbf{T}=
	\left\{
	g=\mathrm{diag}(g_1,\cdots,g_n)\in \mathbf{H}|
	g_1,\cdots,g_n\in K^\times
	\text{ with }
	g_1g_1^\ast
	=
	g_2g_2^\ast
	=
	\cdots
	=
	g_ng_n^\ast
	=
	\mu(g)\in\Q^\times
	\right\}.
	\]	
	So
	in the induced isomorphism
	$\mathbf{H}(\Q_\ell)
	\simeq
	\mathrm{GSp}_{2n}(\Q_\ell)$,
	$\mathbf{T}(\Q_\ell)$
	is mapped to the subgroup of diagonal matrices.
	For any
	$a\in \Q_\ell^\times$
	and any
	$j=1,\cdots,n$,
	we write
	\[
	\tau_j(a)
	=
	\mathrm{diag}
	(1_{2(j-1)},
	a,
	1/a,
	1_{2(n-j)})
	\in
	\mathrm{GSp}_{2n}(\Q_\ell).
	\]
	For any $j=1,\cdots,n$,
	we write
	\[
	u_j(t)
	=
	(\iota\circ\iota_j\circ\iota^{-1})\left(\begin{pmatrix}
		1 & t \\ 0 & 1
	\end{pmatrix}\right)\in\mathbf{H}(\Q_\ell).
	\]
	Then one checks easily the conjugate action $\mathrm{Ad}$ of $\tau_j(a)$ on $\mathbf{U}_i$ is given by
	\[
	\mathrm{Ad}_{\tau_j(a)}(u_i(v))
	=
	\begin{cases*}
		u_i(a^2v),
		&
		$i=j=1$;
		\\
		u_i(av),
		&
		$i\ne j=1$;
		\\
		u_i(a^{-1}v),
		&
		$i=j\ne1$;
		\\
		u_i(v),
		&
		$i\ne j\ne1$.
	\end{cases*}
	\]
	Using the above formulas
	and the argument similar to the preceding examples,
	one verifies easily that $\mathbf{H}_1^1(\Q_\ell),\cdots,\mathbf{H}_n^1(\Q_\ell)$
	generate
	$\mathbf{H}^1(\Q_\ell)$
	and
	condition \ref{condition on mathbf{H}}
	is satisfied.

	\printindex

\end{document}